\documentclass[11pt,reqno,twoside,a4paper]{amsart}

\usepackage{amsmath, amssymb, amscd, amsthm}

\theoremstyle{plain}
\newtheorem{thm}{Theorem}[section]

\newtheorem{lem}[thm]{Lemma}
\newtheorem{prop}[thm]{Proposition}

\theoremstyle{remark}
\newtheorem{rem}[thm]{Remark}

\newcounter{sspar}[subsection]
\renewcommand\thesspar{(\thesubsection.\arabic{sspar})}
    {\par\ \newline
     \vskip-\baselineskip\vskip.1truecm
     \noindent\refstepcounter{sspar}
     \noindent\textbf{\thesspar} \ignorespaces}
    {\vskip-\baselineskip
    \ignorespaces}
    {\refstepcounter{sspar}
     \textup{\textbf{\thesspar}} \ignorespaces}
    {\vskip-\baselineskip
    \ignorespaces}

\newcommand{\R}{\mathbb R}

\newcommand{\C}{\mathbb C}
\newcommand{\Z}{\mathbb Z}
\newcommand{\T}{\mathbb T}
\newcommand{\al}{\alpha}
\newcommand{\be}{\beta}
\newcommand{\ga}{\gamma}
\newcommand{\Ga}{\Gamma}
\newcommand{\de}{\delta}
\newcommand{\De}{\Delta}
\newcommand{\eps}{\varepsilon}
\newcommand{\si}{\sigma}
\newcommand{\te}{\theta}

\newcommand{\la}{\lambda}

\newcommand{\om}{\omega}
\newcommand{\Om}{\Omega}
\newcommand{\inprod}[2]{\langle #1,#2 \rangle}
\newcommand{\tensor}{\otimes}
\newcommand{\dirint}{\sideset{}{^\oplus}\int\limits_{0}^{\infty}}
\newcommand{\qdirint}{\sideset{}{^\oplus}\int\limits_{0}^{-\frac{\pi}{2\ln
q}}}
\newcommand{\F}[5]{\,_{#1}F_{#2} \left( \genfrac{.}{.}{0pt}{}{#3}{#4}
\ ;#5 \right)}
\newcommand{\ph}[5]{\,_{#1}\varphi_{#2}\!\left( \genfrac{.}{.}{0pt}{}{#3}{#4}
\,;#5 \right)}
\newcommand{\hf}{\frac{1}{2}}
\newcommand{\su}{\mathfrak{su}}
\newcommand{\U}{\mathcal U}
\newcommand{\Res}[1]{\underset{#1}{\mathrm{Res}}\,}
\newcommand{\Csp}{\mathrm{span}_\C }

\numberwithin{equation}{section}

\begin{document}
\title{Wilson function transforms related to Racah coefficients}
\author{Wolter Groenevelt}
\address{Department of Mathematics\\
Chalmers University of Technology and G\"oteborg University\\
SE-412 96 G\"oteborg\\
Sweden}
\email{wolter@math.chalmers.se}
\thanks{The author is supported by a NWO-TALENT stipendium of the Netherlands Organization for Scientific Research (NWO)}

\begin{abstract}
The irreducible $*$-representations of the Lie algebra $\su(1,1)$ consist of discrete series representations, principal unitary series and complementary series. We calculate Racah coefficients for tensor product representations that consist of at least two discrete series representations. We use the explicit expressions for the Clebsch-Gordan coefficients as hypergeometric functions to find explicit expressions for the Racah coefficients. The Racah coefficients are Wilson polynomials and Wilson functions. This leads to natural interpretations of the Wilson function transforms. As an application several sum and integral identities are obtained involving Wilson polynomials and Wilson functions. We also compute Racah coefficients for $\U_q(\su(1,1))$, which turn out to be Askey-Wilson functions and Askey-Wilson polynomials.
\end{abstract}
\maketitle


\section{Introduction}
In this paper we study Racah coefficients, or $6j$-symbols, for representations of the Lie algebra $\su(1,1)$. The Racah coefficients for certain tensor product representations of $\su(1,1)$ lead to unitary integral transforms with a very-well-poised $_7F_6$-function, called a Wilson function, as a kernel. These Wilson functions and the corresponding integral transforms were recently introduced by the author in \cite{Gr03} with the applications in this paper in mind. We also consider a $q$-analogue of $\su(1,1)$, namely the quantized universal enveloping algebra $\U_q(\su(1,1))$. From the Racah coefficients for certain tensor product representations of $\mathcal U_q(\su(1,1))$ we obtain a new interpretation of the Askey-Wilson functions transform introduced by Koelink and Stokman \cite{KSt01A}. As an application we obtain several identities for special functions involving (Askey-) Wilson functions and polynomials.

Racah coefficients for $\su(2)$ play an important role in the theory of angular momentum in quantum physics \cite{BL81}. They were first studied in the 1940's by Racah \cite{Ra42}, who also obtained an explicit expression as a finite single sum for these coefficients. Only much later it was realized that the Racah coefficients are multiples of polynomials of hypergeometric type, so that the orthogonality relations for the Racah coefficients lead to discrete orthogonality relations for certain polynomials. These polynomials are nowadays called the Racah polynomials, which can be defined by
\[
R_n(x)=R_n(x;\al,\be,\ga,\de) = \F{4}{3}{ -n, n+\al+\be+1, -x, x+ \ga+\de+1}{\al+1, \be+\de+1, \ga+1}{1},
\]
where one of the lower parameters is equal to $-N$, $N \in \Z_{\geq 0}$, and $0 \leq n \leq N$. These are polynomials of degree $n$ in the variable $x(x+ \ga + \de +1)$, and they are orthogonal on the set $\{0,1,\ldots,N\}$;
\[
\sum_{x=0}^N R_m(x) R_n(x) W(x) = 0, \qquad m \neq n, 
\] 
where 
\[
W(x) = W(x;\al,\be,\ga,\de) =\frac{ 2\ga +2\de+2x+1}{ 2\ga+2\de+1}\, \frac{ (\al+1)_x (\be+\de+1)_x (\ga+1)_x (\ga+\de+1)_x }{(\ga+\de-\al+1)_x (\ga-\be+1)_x (\de+1)_x\, x!}.
\] 
From their definition it is clear that the Racah polynomials satisfy the \emph{duality property}
\[
R_n(x;\al,\be,\ga,\de) = R_x(n;\ga,\de,\al,\be),
\]
so the Racah polynomials are also polynomials of degree $x$ in $n(n+\al+\be+1)$. This shows that the dual orthogonality relation for the Racah polynomials are essentially the same as the orthogonality relations. In section \ref{sec:2} it is shown how the orthogonality relations for the Racah polynomials can be obtained directly from their interpretation as Racah coefficients. We refer to the lecture notes by Van der Jeugt \cite{VdJ03} for a detailed treatement of how properties of Racah polynomials can be obtained from their interpretation as Racah coefficients, see also Vilenkin and Klimyk \cite{VK91}.

The Wilson polynomials, introduced by Wilson \cite{Wil80}, can be defined by
\[
P_n(x)= \F{4}{3}{-n ,n+a+b+c+d-1,a+ix, a-ix}{a+b,a+c,a+d}{1}.
\]
Assuming that $a,b,c,d>0$ these polynomials satisfy continuous orthogonality relations on $[0,\infty)$, with the notation $f(a \pm b)=f(a+b)f(a-b)$, 
\[
\frac{1}{2\pi}\int_0^\infty P_m(x) P_n(x) \frac{ \Ga( a \pm ix) \Ga(b \pm ix) \Ga(c \pm ix) \Ga(d \pm ix) }{ \Ga(\pm 2ix) }dx =0, \qquad m \neq n,
\]
so we may consider the Wilson polynomials as continuous extensions of the Racah polynomials. 

Together with the Racah polynomials the Wilson polynomials form the top level of the Askey-scheme of hypergeometric polynomials \cite{KS98}. This is a large scheme of orthogonal polynomials in one variable that can be written as hypergeometric series. All the polynomials below the top level in the Askey-scheme, and many of their properties, can (at least formally) be obtained by taking appropriate limits of the Wilson or Racah polynomials. Therefore these polynomials are fundamental objects in the theory of hypergeometric orthogonal polynomials.

In \cite{Gr03} the Wilson functions were introduced. The Wilson function can be considered as an analytic continuation of the Wilson polynomial in its degree, or equivalently as continuous extension of the Racah polynomials in both the variable $x$ and the dual variable (the degree) $n$. The Wilson funtions satisfy generalized orthogonality relations, in the sense that they are the kernel in two different unitary integral transforms, called the Wilson transforms I and II. The dual integral transforms are essentially the same as the integral transforms, which can be considered as continuous analogues of the fact that the orthogonality relations and the dual relations for the Racah polynomials are essentially the same.

The main goal of this paper is to show that the Wilson functions have a natural interpretation as Racah coefficients for the Lie algebra $\su(1,1)$. Both Wilson function transforms appear in this context, and unitarity of both integral transforms can be obtained in this way. Moreover, we show that the Wilson polynomials also appear as Racah coefficients for $\su(1,1)$. As an application we obtain several sum or integral identities for Wilson polynomials and Wilson functions.
 
Let us also mention some other generalizations of the Racah polynomials. Askey and Wilson \cite{AW79} considered basic hypergeometric extensions of the Racah polynomials, and called these the $q$-Racah polynomials. The $q$-Racah polynomials satisfy orthogonality relations on a finite discrete set, similar to the Racah polynomials. Later Kirillov and Reshetikhin gave an interpretation of the $q$-Racah polynomials as Racah coefficients for the quantized universal enveloping algebra $\U_q(\su(2))$. Askey and Wilson \cite{AW85} showed that the $q$-Racah polynomials also satisfy continuous orthogonality relations, and considered as a polynomial in a continuous variable, these polynomials are the now famous Askey-Wilson polynomials. Koelink and Stokman \cite{KSt01A} introduced the Askey-Wilson functions, which can be considered as the analytic continuation of the Askey-Wilson polynomials in their degree, and they studied the corresponding integral transform. In section \ref{sec:Uq} we show that the Askey-Wilson polynomials and the Askey-Wilson funtions have an interpretation as Racah coefficients for $\U_q(\su(1,1))$.

Let us now describe the contents of this paper. In section \ref{sec:2} we recall the irreducible $*$-representations of the Lie algebra $\su(1,1)$. The irreducible representations consist of four classes: positive and negative discrete series $\pi^\pm$, principle unitary series $\pi^P$ and complementary series $\pi^C$. Then we calculate Clebsch-Gordan coefficients and Racah coefficients related to tensor product representations of positive discrete series. The results from section \ref{sec:2} are well known in the literature. This section is mainly included to introduce notations and to explain the ideas we use later on. Also we use the results from section \ref{sec:2} in following sections. In section \ref{sec:CGC} we first decompose $2$-fold tensor product representations, of which one is a discrete series, into irreducible representations, and we calculate the corresponding Clebsch-Gordan coefficients. The results from this section are used later on in sections \ref{sec:+-+}, \ref{sec:+P-}, \ref{sec:more} to compute Racah coefficients related to various $3$-fold tensor product representation of $\su(1,1)$. In section \ref{sec:+-+}  we compute the Racah coefficients for the tensor product $\pi^+ \tensor \pi^-\tensor \pi^+$. This leads to a natural interpretation of the Wilson transform I. In section \ref{sec:+P-} we consider Racah coefficients for the tensor product $\pi^+ \tensor \pi^P \tensor \pi^-$, and in this way we obtain a natural interpretation of the Wilson transform II. In section \ref{sec:more} we give the Racah coefficients for some other tensor products involving at least two discrete series representation. In section \ref{sec:identities} we use natural identities for Racah coefficients to obtain sum and integral identities involving Wilson polynomials and Wilson funtions. Finally in section \ref{sec:Uq} we sketch how the Askey-Wilson polynomials and the Askey-Wilson functions can be obtained as Racah coefficients for $\U_q(\su(1,1))$. 

We use the standard notations for (basic) hypergeometric functions as in the book by Gasper and Rahman \cite{GRa04}.\\

\emph{Acknowledgement.} I would like to thank Erik Koelink for many useful discussions and suggestions.


\section{Representations of $\su(1,1)$ and Racah polynomials} \label{sec:2}
This section is meant to explain the ideas we use later on in the paper to find Racah coefficients, and to introduce some notations. 

\subsection{The Lie algebra $\su(1,1)$} 
The Lie algebra $\su(1,1)$ is generated by $H$, $E$ and $F$, which satisfy the commutation relations
\[
[H,E]=2E, \quad [H,F]=-2F, \quad [E,F]=H.
\]
There is a $\ast$-structure defined by $H^*=H$ and $E^*=-F$. The center of the universal enveloping algebra $\U(\su(1,1))$ is generated by the Casimir element $\Omega$, which is given by
\begin{equation} \label{Casimir}
\Omega = -\frac{1}{4}(H^2+2H+4FE)= -\frac14(H^2-2H+4EF).
\end{equation}
The irreducible $*$-representations of $\su(1,1)$ are determined by Bargmann \cite{Ba47}, see e.g.~ also \cite[\S6.4]{VK91}. There are, besides the trivial representation, four classes of irreducible $*$-representations of $\su(1,1)$:\\

\emph{Positive discrete series.}
The positive discrete series representations $\pi_k^+$ are labelled
by $k>0$. The representation space is $\ell^2(\mathbb{Z}_+)$ with
orthonormal basis $\{e_n\}_{n \in \Z_{\geq 0}}$. The action is given by
\begin{equation} \label{def:pos}
\begin{split}
\pi_k^+(H)\, e_n =&\ 2(k+n)\, e_n, \\
\pi_k^+(E)\, e_n =&\ \sqrt{(n+1)(2k+n)}\, e_{n+1}, \\
\pi_k^+(F)\, e_n =&\ -\sqrt{n(2k+n-1)} \,e_{n-1}, \\
\pi_k^+(\Omega)\, e_n =&\ k(1-k)\, e_n.
\end{split}
\end{equation}

\emph{Negative discrete series.}
The negative discrete series representations $\pi_{k}^-$ are labelled by
$k>0$. The representation space is $\ell^2(\Z_{\geq 0})$ with orthonormal
basis $\{e_n\}_{n \in \Z_{\geq 0}}$. The action is given by
\begin{equation} \label{def:neg}
\begin{split}
\pi_{k}^-(H)\, e_n =&\ -2(k+n)\, e_n, \\
\pi_{k}^-(E)\, e_n =&\ -\sqrt{n(2k+n-1)}\, e_{n-1},  \\
\pi_{k}^-(F)\, e_n =&\ \sqrt{(n+1)(2k+n)}\, e_{n+1},  \\
\pi_{k}^-(\Omega)\, e_n =&\ k(1-k)\, e_n.
\end{split}
\end{equation}

\emph{Principal unitary series.}
The principal unitary series representations $\pi^P_{\rho,\eps}$ are labelled
by  $\eps \in [0,1)$ and $\rho \geq 0$, where $(\rho,\eps)
\neq (0,\frac{1}{2})$. The representation space
is $\ell^2(\mathbb{Z})$ with orthonormal basis $\{e_n\}_{n \in \Z}$.
The action is given by
\begin{equation} \label{def:princ}
\begin{split}
\pi_{\rho, \eps}^P(H)\, e_n =&\ 2(\eps+n)\, e_n, \\
\pi_{\rho, \eps}^P(E)\, e_n =&\ \sqrt{(n+\eps+\frac{1}{2}-i\rho)
(n+\eps+\frac{1}{2}+i\rho)}\, e_{n+1}, \\
\pi_{\rho, \eps}^P(F)\, e_n =&\ -\sqrt{(n+\eps-\frac{1}{2}-i\rho)
(n+\eps-\frac{1}{2}+i\rho)}\, e_{n-1}, \\
\pi_{\rho, \eps}^P(\Omega)\, e_n =&\ (\rho^2+\frac{1}{4})\, e_n.
\end{split}
\end{equation}
For $(\rho,\eps)= (0,\hf)$ the representation $\pi^P_{0,\hf}$ splits into a direct sum of a positive and a negative discrete series representation: $\pi^P_{0,\hf}= \pi^+_\hf \oplus \pi^-_\hf$. The representation space splits into two invariant subspaces: $\{e_n \, |\, n<0 \} \oplus \{ e_n\, | \, n \geq 0 \}$. \\

\emph{Complementary series.}
The complementary series representations $\pi^C_{\lambda,\eps}$ are
labelled by $\eps$ and $\lambda$, where $\eps \in [0,\frac{1}{2})$ and
$\lambda \in (-\frac{1}{2},-\eps)$  or $\eps \in (\frac{1}{2},1)$ and
$\lambda \in (-\frac{1}{2},\eps-1)$. The representation space is
$\ell^2(\Z)$ with orthonormal basis $\{e_n\}_{n \in \Z}$. The actions of the generators $H$, $E$, $F$ are the same as in the principal unitary series, with $i\rho$ replaced by $\la+\hf$. \\

We remark that in case $\eps \not\in [0,1)$ the principle unitary series $\pi^P_{\rho,\eps}$ are unitarily equivalent to the $\pi^P_{\rho,\eps+p}$, where $p \in \Z$ is such that $\eps+p \in [0,1)$. This follows from the fact that for $\eps\in [p,p+1)$, $p \in \Z$, and for $X \in \su(1,1)$,  the action of $\pi^P_{\rho,\eps}(X)$ on the basis $\{e_n \}_{n \in \Z}$ is the same as the action of $\pi^P_{\rho,\eps-p}$ acting on the basis $\{e_{n+p}\}_{n \in \Z}$. A similar observation can be made for the complementary series.

The operators defined above are unbounded closable operators, with domain the finite linear combinations of the basis vectors. The representations are $*$-representations in the sense of Schm\"udgen \cite[Ch.8]{Sc90}.

\subsection{Clebsch-Gordan coefficients for ${\pi^+ \tensor \pi^+}$} \label{ssec:CGC}
We consider the tensor product representation $\pi_{k_1}^+ \tensor \pi_{k_2}^+$. Recall that for the definition of the a tensor product representation we need the coproduct. For a Lie algebra element $X$ the coproduct $\De$ is defined by $\De(X) = 1 \tensor X + X \tensor 1$. The coproduct can be extended to the universal enveloping algebra as an algebra morphism, and then we obtain from \eqref{Casimir}  for the Casimir element
\begin{equation} \label{eq De Om}
\Delta(\Omega) = 1 \tensor \Omega + \Omega \tensor 1 - \frac{1}{2} H
\tensor H - F \tensor E - E \tensor F.
\end{equation}

In order to decompose the tensor product $\pi^+_{k_1} \tensor \pi^+_{k_2}$ into irreducible representation we diagonalize the operator $\pi^+_{k_1} \tensor \pi^+_{k_2}(\De(\Om))$. Applying this operator to $e_{n} \tensor e_{p-n} \in \ell^2(\Z_{\geq 0})^{\tensor 2}$ and using \eqref{def:pos} and \eqref{eq De Om}, we find an expression of the form
\[
\pi^+_{k_1} \tensor \pi^+_{k_2}(\De(\Om)) \, e_{n} \tensor e_{p-n}= a_n\, e_{n+1} \tensor e_{p-n-1} + b_n\ e_{n} \tensor e_{p-n} + a_{n-1} e_{n-1} \tensor e_{p-n+1},
\]
for some coefficients $a_n$ and $b_n$. From writing out $a_n$ explicitly we see that $a_{-1}=a_{p+1}=0$. So for $p \in \Z_{\geq 0}$ the space
\[
\mathcal H^p =  \Csp \{ e_n \tensor e_{p-n} \ | \ n =0,\ldots,p \} \cong \C^{p+1}
\]
is invariant under the action of the Casimir operator. Let $v^p(\la)=\sum_{n=0}^p v_n^p(\la) e_n \tensor e_{p-n} \in \mathcal H^p$ be an eigenvector of $\pi^+_{k_1} \tensor \pi^+_{k_2}(\De(\Om))$ for eigenvalue $\la$, then the coefficients $v_n^p$ are solutions of
\[
\la v_n^p(\la) = a_n v_{n+1}^p(\la) + b_n v_n^p(\la) + a_{n-1}v_{n-1}^p(\la).
\]
This shows that the coefficients $v_n^p$ can be generated from $v_0^p$. To give an explicit expression for the coefficients $v_n^p$ we introduce the dual Hahn polynomials.

The dual Hahn polynomials are defined by
\[
T_n(x;\ga,\de,N) = \F{3}{2}{-n, -x, x+ \ga+\de+1}{\ga+1,-N}{1},\qquad n=0,1,\ldots,N,
\]
where $N$ is a non-negative integer. These are polynomials in $x(x+\ga+\de+1)$. Note that the dual Hahn polynomial can be obtained as a limit case of a Racah polynomial; $T_n(x;\ga,\de,N)= \lim_{\be \rightarrow \infty} R_n(x;-N-1,\be,\ga,\de)$. The dual Hahn polynomial is not symmetric in $\ga$ and $\de$, but we can interchange $\ga$ and $\de$ using 
\begin{equation} \label{eq:symT}
T_n(x;\ga,\de,N) = (-1)^x \frac{ (\de+1)_x }{(\ga+1)_x} T_{N-n}(x;\de,\ga,N),
\end{equation}
where we assume that $x \in \{0,1,\ldots,N\}$. The orthonormal dual Hahn polynomials are defined by
\[
t_n(x) = t_n(x;\ga,\de,N) = (-1)^n\sqrt{ \frac{ (\ga+1)_n (-N)_n }{n!\, (-\de-N)_n} }\, T_n(x;\ga,\de,N),
\]
and they satisfy, for $\ga>-1$ and $\de>-1$ or for $\ga<-N$ and $\de<-N$, 
\[
\begin{split}
&\sum_{x=0}^N t_m(x) t_n(x) W(x;\ga,\de,N) = \de_{mn},\\
&W(x;\ga,\de,N) = \frac{ (2x+\ga+\de+1) (\ga+1)_x (-N)_x (\de+1)_N}{ (-1)^x (x+ \ga + \de +1)_{N+1} (\de+1)_x x! }.
\end{split}
\]
The dual Hahn polynomials satisfy the following three-term recurrence relation,
\[
-x(x+\ga+\de+1)\, t_n(x)= a_n t_{n+1}(x) + b_n t_n(x) + a_{n-1} t_{n-1}(x),
\]
where
\[
\begin{split}
a_n &= \sqrt{ (n+1) (n+ \ga+1) (n- \de-N) (n-N) },\\
b_n &= (n+\ga+1)(n-N) + n(n-\de-N-1).\\
\end{split}
\]
Note that for $x\in \{0, \ldots, N\}$ the polynomial $T_n(x)$ is also a polynomial in $n$ of degree $x$. Considered as a polynomial in $n$, $T_n(x)$ is called a Hahn polynomial.

We see that the coefficients $v^p_n$ are multiples of dual Hahn polynomials with parameters $(\ga,\de,N)=(2k_1-1,2k_2-1,p)$, and the eigenvalue $\la$ is given by $\la(x)=(k_1+k_2+x)(1-k_1-k_2-x)$. We normalize the coefficients $v_n^p$ as follows. We split the weight function for the dual Hahn polynomials in a $p$-dependent and a $p$-independent part, such that the $p$-dependent part equals $1$ for $p=0$, and we glue the square root of the $p$-dependent part to the dual Hahn polynomial. So for $p \in \Z_{\geq 0}$ we set
\[
v_n^p(x;k_1,k_2)=\sqrt{ \frac{(-1)^x (-p)_x (2k_2)_p }{(2k_1+2k_2+x)_p } }\ t_n(x;2k_1-1, 2k_2-1, p), \qquad x=0,1,\ldots,p.
\]
Let $W(x;k_1,k_2)$ be the weight function given by
\[
W(x;k_1,k_2) = \frac{ (2x+2k_1+2k_2-1) (2k_1)_x }{ (x+2k_1+2k_2-1) (2k_2)_x x! },\qquad x \in \Z_{\geq 0}.
\]
We denote by $\C^{p+1}(W(x;k_1,k_2))$ the $(p+1)$-dimensional Hilbert space with the inner product associated to this weight function at the points $x \in \{0,1,\ldots,p\}$. Then, for fixed $p$, the functions $x\mapsto v_n^p(x;k_1,k_2)$ form an orthonormal basis for $\C^{p+1}(W(x;k_1,k_2))$. We remark that it is customary to multiply the coefficients $v_n^p$ by the square root of the weight function $W$, since then they are orthogonal in $\C^{p+1}$ with respect to the usual inner product. We will work instead with weighted inner products.  We can now formulate our finding as follows.
\begin{prop} \label{prop I}
For $p\in \Z_{\geq 0}$, the operator $I: \mathcal H^p \rightarrow \C^{p+1}(W(x;k_1,k_2))$ defined by
\[
I : e_n \tensor e_{p-n} \mapsto \Big(x \mapsto v_n^p(x;k_1,k_2)\Big),
\]
is unitary and intertwines $\pi^+_{k_1} \tensor \pi^+_{k_2}(\De(\Om))$ with $M_{(k_1+k_2+x)(1-k_1-k_2-x)}$. 
\end{prop}
Here $M_a$ is denotes multiplication by $a$. We have found that 
\[
v^p(x)=\sum_{n=0}^p v_n^p(x;k_1,k_2)\, e_{n} \tensor e_{n-p} 
\]
is an eigenvector of $\pi^+_{k_1} \tensor \pi^+_{k_2}(\De(\Om))$ for eigenvalue $(k_1+k_2+x)(1-k_1-k_2-x)$. The squared norm of $v^p(x)$ is $W(x;k_1,k_2)^{-1}$. We write $v^p(x)=v^p(x;k_1,k_2)$ if we need to stress the dependence on the representation labels $k_1$ and $k_2$. A comparison of the multiplication in Proposition \ref{prop I} with the action of $\Om$ in the discrete series representation \eqref{def:pos} and \eqref{def:neg}, 
suggests that for all $x \in \Z_{\geq 0}$ the discrete series representation $\pi_{k_1+k_2+x}$ appears in the decomposition of $\pi^+_{k_1} \tensor \pi^+_{k_2}$. To find the exact decomposition we determine the action of the generators $H$, $E$ and $F$ on the eigenvector $v^p(x)$. The action of $H$ is easy to find, and the actions of $E$ and $F$ follow from the following contiguous relations:
\begin{equation} \label{eq:cont.rel}
\begin{split}
(d-1)&\F{3}{2}{a,b,c}{d-1,e}{1} = a \F{3}{2}{a+1,b,c}{d,e}{1} + (a-d+1) \F{3}{2}{a,b,c}{d,e}{1}, \\
(d-b)(d-c)& \F{3}{2}{a,b,c}{d+1,e}{1} = \\
&d(e-a)\F{3}{2}{a-1,b,c}{d,e}{1} +d(a+b+c-d-e) \F{3}{2}{a,b,c}{d,e}{1}.\\
\end{split}
\end{equation}
The first relation is obtained from combining (3.7.9) and (3.7.11) from \cite{AAR99}, the second relation from combining (3.7.10), (3.7.12) and (3.7.14) from \cite{AAR99}. This gives us the following actions
\[
\begin{split}
\pi^+_{k_1} \tensor \pi^+_{k_2}(\De(H)) \, v^p(x) &= 2(k_1+k_2+p)\, v^p(x),\\
\pi^+_{k_1} \tensor \pi^+_{k_2}(\De(E)) \, v^p(x)& = \sqrt{(p-x+1)(2k_1+2k_2+x+p)}\, v^{p+1}(x),\\
\pi^+_{k_1} \tensor \pi^+_{k_2}(\De(F)) \, v^p(x) &= \sqrt{(p-x)(2k_1+2k_2+x+p-1)}\, v^{p-1}(x).
\end{split}
\]
Observe that these are the same as the actions of the generators in the representation $\pi^+_{k_1+k_2+x}$ on the basis vectors $e_{p-x}$. Therefore we define 
\[
\begin{split}
I': \ell^2(\Z_{\geq 0})^{\tensor 2} &\rightarrow \bigoplus_{x=0}^\infty \ell^2(\Z_{\geq 0})\, W(x;k_1,k_2),\\
e_{n_1} \tensor e_{n_2} &\mapsto \sum_{x=0}^{n_1+n_2} v_{n_1}^{n_1+n_2}(x;k_1,k_2)\,e_{n_1+n_2-x} ,
\end{split} 
\]
which corresponds to inverting $v^p(x) = \sum_n v_n^p(x;k_1,k_2)\, e_n \tensor e_{p-n}$ using the orthogonality relations for the dual Hahn polynomials. Here the weighted direct sum can be considered as a direct integral with a measure $\mu$ supported on $\Z_{\geq 0}$, satisfying $\mu(\{x\})=W(x;k_1,k_2)$ for $x \in \Z_{\geq 0}$. Now we have the following proposition.
\begin{prop} \label{prop:++}
The operator $I_1'$ intertwines $\pi^+_{k_1} \tensor \pi^+_{k_2}$ with $\displaystyle \bigoplus_{x=0}^\infty \pi^+_{k_1+k_2+x}$.
\end{prop}
Using the notation $\pi_1 \cong \pi_2$ if the two representations $\pi_i$, $i=1,2$, are unitarily equivalent, Proposition \ref{prop:++} states that the tensor product representation $\pi^+_{k_1} \tensor \pi^+_{k_2}$ is decomposed into irredubible representations as
\[
\pi_{k_1}^+ \tensor \pi^+_{k_2} \cong \bigoplus_{j=0}^\infty \pi^+_{k_1+k_2+j},
\]
and the corresponding Clebsch-Gordan coefficients are multiples of (dual) Hahn polynomials. See also \cite{KVdJ98}, \cite{Koo81}, \cite{SSS84} for the interpretation of (dual) Hahn polynomials as Clebsch-Gordan coefficients.

\subsection{Racah coefficients for $\pi^+ \tensor \pi^+ \tensor \pi^+$}
To determine the Racah coefficients for the tensor product $\pi^+_{k_1} \tensor \pi^+_{k_2} \tensor \pi^+_{k_3}$ we first determine the eigenvectors of $\Om$.  We consider the operator $\pi^+_{k_1} \tensor \pi^+_{k_2} \tensor \pi^+_{k_3}(\De\tensor 1)(\De(\Om))$. Using the actions of $H, E, F$ on $v^p(x)$, we obtain that $\pi^+_{k_1} \tensor \pi^+_{k_2} \tensor \pi^+_{k_3}(\De \tensor 1)(\De(\Om))$ acts on $v^{p-m}(x) \tensor e_m$ in the same way as $\pi^+_{k_1+k_2+x} \tensor \pi^+_{k_3}$ acts on $e_{p-x-m} \tensor e_m$. So we see that the vector
\[
\begin{split}
w_1^{p-x_1}(x_2)\ &= \sum_{m=0}^{p-x_1} v_{n_3}^{p-x_1}(x_2;k_3,k_1+k_2+x_1)\, v^{p-n_3}(x_1;k_1,k_2) \tensor e_{n_3}\\
& = \sum_{n_3=0}^{p-x_1} \sum_{n_1=0}^{p-n_3}  v_{n_1}^{p-n_3}(x_1;k_1,k_2) v_{n_3}^{p-x_1}(x_2;k_3,k_1+k_2+x_1) \, e_{n_1} \tensor e_{n_2} \tensor e_{n_3},
\end{split}
\]
where $n_1+n_2+n_3=p$, is an eigenvector of $\pi^+_{k_1} \tensor \pi^+_{k_2} \tensor \pi^+_{k_3}(\De \tensor 1)(\De(\Om))$ for the eigenvalue $(k_1+k_2+k_3+x_1+x_2)(1-k_1-k_2-k_3-x_1-x_2)$. Moreover, the actions of $H$, $E$, $F$, on $w_1^{p-x_1}(x_2)$ are the same as the actions of the generators in the representation $\pi^+_{k_1+k_2+k_3+x_1+x_2}$ on $e_{p-x_1-x_2}$. The squared norm of the vector $w_1^{p-x_1}(x_2)$ is $\big(W(x_1;k_1,k_2)W(x_2;k_3,k_1+k_2+x_1,) \big)^{-1}$.

In the same way we find that the vector
\[
\begin{split}
w_2^{p-y_1}(y_2)\ & = \sum_{n_1=0}^{p-y_1} v_{n_1}^{p-y_1}(y_2;k_1,k_2+k_3+y_1)\, e_{n_1} \tensor v^{p-n_1}(y_1;k_3,k_2)\\
&= \sum_{n_1=0}^{p-y_1} \sum_{n_3=0}^{p-n_1} v_{n_3}^{p-n_1}(y_1;k_3,k_2) v_{n_1}^{p-y_1}(y_2;k_1,k_2+k_3+y_1)\, e_{n_1} \tensor e_{n_2} \tensor e_{n_3},
\end{split}
\]
where $p=n_1+n_2+n_3$, is an eigenvector of $\pi^+_{k_1} \tensor \pi^+_{k_2} \tensor \pi^+_{k_3}(1 \tensor \De)(\De(\Om))$ for the eigenvalue $(k_1+k_2+k_3+y_1+y_2)(1-k_1-k_2-k_3-y_1-y_2)$. 

From the associativity of the tensor product, or equivalently the coassociativity of the coproduct, i.e.~$(\De \tensor 1)\De = (1 \tensor \De) \De$, it follows that we have found two bases of eigenvectors for the same operator. Let us consider the overlap coefficients for the transition from one basis to the other, i.e.~the coefficients
\begin{equation} \label{def:Racah coeff}
C_{p_1,p_2}(x_1,x_2;y_1,y_2;k_1,k_2,k_3)= \inprod{w_1^{p_1-x_1}(x_2)}{w_2^{p_2-y_1}(y_2)}.
\end{equation}
Here the $\inprod{ \cdot}{\cdot}$ denotes the usual inner product in $\ell^2(\Z_{\geq 0})^{\tensor 3}$. Since both $w_1$ and $w_2$ are eigenvectors of $\Om$, the coefficient $C$ vanishes if $x_1+x_2 \neq y_1+y_2$. Also, from applying $H$ to $w_1$ and using $H^*=H$, we see that the coefficient vanishes if $p_1 \neq p_2$. Similarly, setting $p_1=p_2=p$, applying $E$ on $w_1$ and using $E^*=-F$, it follows that $C$ is independent of $p$. From here on we denote, $C_{p,p}(x_1,x_2;y_1,y_2;k_1,k_2,k_3)= U(x_1,y_1;k_1,k_2,k_3,t)$, where we assume that $t=x_1+x_2=y_1+y_2$. The coefficient $U$ is called the Racah coefficient.

Let us define a weight function in the variable $x$ as follows;
\[
W(x;k_1,k_2,k_3,t) = W(x;k_1,k_2)W(t-x;k_3,k_1+k_2+x,).
\]
So this weightfunction is the inverse of the squared norm of $w_1^{p-x}(t-x)$. Note that interchanging $k_1$ and $k_3$ gives the inverse of the squared norm of $w_2^{p-x}(t-x)$.
Now we have
\begin{equation} \label{eq:expan1}
\begin{split}
w_2^{p-y}(t-y) &= \sum_{x=0}^{t} U(x,y;k_1,k_2,k_3,t) w_1^{p-x}(t-x) W(x;k_1,k_2,k_3,t),\\
w_1^{p-x}(t-x) &= \sum_{y=0}^{t} U(x,y;k_1,k_2,k_3,t) w_2^{p-y}(t-y) W(y;k_3,k_2,k_1,t).
\end{split}
\end{equation}
This gives us the following orthogonality relations for the Racah coefficients
\[
\begin{split}
\sum_{x=0}^t U(x,y;k_1,k_2,k_3,t) U(x,y';k_1,k_2,k_3,t) W(x;k_1,k_2,k_3,t) = \frac{\de_{yy'} }{W(y;k_3,k_2,k_1,t)},\\
\sum_{y=0}^t U(x,y;k_1,k_2,k_3,t) U(x',y;k_1,k_2,k_3,t) W(y;k_3,k_2,k_1,t) = \frac{\de_{xx'} }{W(x;k_1,k_2,k_3,t)}.
\end{split}
\]
Also note that there is an obvious symmetry property;
\[
U(x,y;k_1,k_2,k_3,t) =U(y,x;k_3,k_2,k_1,t).
\]

In order to find an explicit expression for the Racah coefficient, we take the inner product with $e_{n_1} \tensor e_{n_2} \tensor e_{n_3}$, $n_1+n_2+n_3=p$, in the first identity in \eqref{eq:expan1}. This gives
\begin{equation} \label{eq:expan2}
\begin{split}
&v_{n_3}^{p-n_1}(y;k_3,k_2) v_{n_1}^{p-y}(t-y;k_1,k_2+k_3+y) = \\
&\ \sum_{x=0}^t
U(x,y;k_1,k_2,k_3,t) v_{n_1}^{p-n_3}(x;k_1,k_2) v_{n_3}^{p-x}(t-x;k_3,k_1+k_2+x)W(x;k_1,k_2,k_3,t).
\end{split}
\end{equation}
Using the orthogonality relation for the Clebsch-Gordan coefficient $v_{n_1}^{p-n_3}$ this formula gives an expansion of the Racah coefficient in Clebsch-Gordan coefficients. Now an explicit expression for the Racah coefficient can be found by specializing $n_3$ and $p$ in a suitable way. This leads to expression as a terminating balanced $_4F_3$-series, and then we find that the Racah coefficient is a multiple of a Racah polynomial,
\[
\begin{split}
U(x,y;k_1,k_2,k_3,t) =&(-1)^{x+y-t} \frac{ t!\,(2k_2)_x (2k_2)_y }{ (2k_1)_x (2k_3)_y } \sqrt{ \frac{ (2k_1+2k_2+2k_3+t)_x (2k_1+2k_2+2k_3+t)_y }{ (2k_1+2k_2+x)_x (2k_2+2k_3+y)_y} }\\
& \times R_y(x;2k_2-1, 2k_3-1, -t-1, 2k_1+2k_2+t-1).
\end{split}
\]
See \cite{VdJ03} for details. The above orthogonality relations for the Racah coefficients lead to orthogonality relations for the Racah polynomials on the finite set $\{0,1,\ldots, N\}$.

In the following sections we study other tensor products for $\su(1,1)$ and the corresponding Racah coefficients. In the cases we consider we cannot give eigenvectors of the Casimir operator, but instead we must work with generalized eigenvectors. Therefore the computations will be more complicated, but underlying ideas are much the same as in this section.


\section{Clebsch-Gordan coefficients} \label{sec:CGC}
In this section we decompose various $2$-fold tensor product representations involving discrete series representations, and we give explicit expressions for the corresponding Clebsch-Gordan coefficients. These decompositions and the corresponding Clebsch-Gordan coefficients are well-known in the literature, see e.g.~\cite{KV99}. In this section we emphasize the relation between the Clebsch-Gordan coefficients and orthogonal polynomials.

\subsection{Clebsch-Gordan coefficients for $\pi^+ \tensor \pi^-$} 
First we consider the tensor product ${\pi^+ \tensor \pi^-}$ acting on $\ell^2(\Z_{\geq 0})^{\tensor 2}$.  We define for fixed $p \in \Z$ elements $f_n^p \in \ell^2(\Z_{\geq 0})^{\tensor 2} $ by
\begin{equation} \label{def fn}
f_n^p = 
\begin{cases}
e_n \tensor e_{n-p}, & p \leq 0, \\
e_{n+p} \tensor e_n, & p \geq 0.
\end{cases}
\end{equation}
So $\{f_n^p \ | \ n \in \Z_{\geq 0}, \ p \in \Z \}$ is an orthonormal basis for $\ell^2(\Z_{\geq 0})^{\tensor 2}$. We define a subspace $\mathcal H^p_1 \subset \ell^2(\Z_{\geq 0})^{\tensor 2}$ by
\[
\mathcal H_1^p = \overline{\Csp\{ f_n^p \ | \ n \in \Z_{\geq 0}\} } \cong \ell^2(\Z_{\geq0 }).
\]
From \eqref{eq De Om}, \eqref{def:pos} and \eqref{def:neg} we find that the operator $\pi^+_{k_1} \tensor \pi^-_{k_2}( \De(\Om))$ acting on elements $f_n^p$ extends to an unbounded Jacobi operator acting on $\mathcal H_1^p$. We refer to the book by Akhiezer \cite{Ak65} for more information on Jacobi operators. The Jacobi operator we have here is the Jacobi operator corresponding to the three-term recurrence relation of the continuous dual Hahn polynomials, cf.~ \cite[Prop.2.1]{GK02}.

The continuous dual Hahn polynomials are polynomials in $x^2$ defined by
\begin{equation} \label{def:cont dHahn}
S_n(x;a,b,c) = (a+b)_n (a+c)_n\F{3}{2}{-n,a+ix,a-ix}{a+b,a+c}{1}.
\end{equation}
The continuous dual Hahn polynomials can be obtained as a limit of the Wilson polynomials \cite{Wil80}, see also \cite{AAR99} or \cite{KS98}. 
By Kummer's transformation \cite[Cor.3.3.5]{AAR99} the polynomials $S_n$ are symmetric in $a$, $b$ and $c$. For real parameters $a$, $b$ and $c$, with $a+b$, $a+c$, $b+c$ positive, or for one real parameter and a pair of complex conjugates with a positive real part, the continuous dual Hahn polynomials are orthogonal with respect to a positive measure. The orthonormal continuous dual Hahn polynomials are given by
\[
s_n(x)= s_n(x;a,b,c) = \frac{(-1)^n S_n(x;a,b,c)} {\sqrt{n!(a+b)_n
(a+c)_n (b+c)_n}}.
\]
Assuming without loss of generality that $a$ is the smallest of the real parameters, the polynomials $s_n$ are orthonormal with respect to the measure $\mu(\cdot;a,b,c)$, given by
\[
\int f(x)\, d\mu(x;a,b,c) = 
\frac{1}{2\pi} \int_0^\infty f(x) w(x;a,b,c)\, dx + \sum_{\substack{k \in \Z_{\geq 0}\\ a+k<0 }}f\big(i(a+k)\big)  w_k(a;b,c)  ,
\]
where
\[
\begin{split}
w(x;a,b,c) =& 
 \frac{ \Ga(a\pm ix) \Ga(b\pm ix) \Ga(c\pm ix) } {\Ga(a+b) \Ga(a+c) \Ga(b+c) \Ga(\pm 2ix) }, \\
w_k(a;b,c) =& i\Res{x=i(a+k)} w(x;a,b,c) \\
=&\frac{ \Ga(b-a) \Ga(c-a) }{ \Ga(-2a) \Ga(b+c)}
\frac{(-1)^k (2a)_k (a+1)_k (a+b)_k (a+c)_k  }{k!\, (a)_k (a-b+1)_k (a-c+1)_k }.
\end{split}
\]
The polynomials $s_n$ form an orthonormal basis for the Hilbert space $\mathcal L(\mu)$ consisting of even functions ($\mu$-almost everywhere) that have finite norm with respect to the inner product associated to the measure $\mu$. The orthonormal continuous dual Hahn polynomials satisfy the following three-term recurrence relation,
\begin{equation} \label{eq:S rec}
x^2 s_n(x) = a_n s_{n+1}(x) +b_n s_n(x) +a_{n-1}s_{n-1}(x),
\end{equation}
where
\[
\begin{split}
a_n & =  \sqrt{(n+1)(n+a+b)(n+a+c)(n+b+c)}, \\ b_n & = 
2n^2+2n(a+b+c-\frac{1}{2})+ab+ac+bc.
\end{split}
\]

The Jacobi operator we obtained from the action of the Casimir corresponds to the three-term recurrence relation for the continuous dual Hahn polynomials with parameters
\[
\{a,b,c\} = \begin{cases}
\{ k_1-k_2+\hf, k_1+k_2-\hf, k_2-k_1-p+\hf\}, & p \leq 0,\\
\{ k_2-k_1+\hf, k_1+k_2-\hf, k_1-k_2+p+\hf\}, & p \geq 0.
\end{cases}
\]
The corresponding orthogonality measure depends also on $p$. Using $\Ga(a+p)=(a)_p \Ga(a)$ for the continuous part of the measure, and $(a+p)_k=(a)_k (a+k)_p/(a)_p$ for the discrete part, we split this orthogonality measure in a $p$-dependent part and a $p$-independent part, and we glue the square root of the $p$-dependent part to the continuous dual Hahn polynomials, i.e.~for $\rho \in [0,\infty)$ we define the functions $F_n^p$ by
\[
F_n^p(\rho;k_1,k_2) = 
\begin{cases}
\displaystyle (-1)^{n+p} \sqrt{ \frac{ (k_2-k_1+\hf\pm i\rho)_{-p} }{ (-p)! \, (2k_2)_{-p}}} \\
\qquad \times s_{n}(\rho;k_1-k_2+\hf, k_1+k_2-\hf, k_2-k_1-p+\hf), & p \leq 0,\\ \\
\displaystyle 
(-1)^n \sqrt{ \frac{ (k_1-k_2+\hf \pm i\rho)_{p} }{ p! \, (2k_1)_{p}}}\\
\qquad \times  s_{n}(\rho;k_2-k_1+\hf, k_1+k_2-\hf, k_1-k_2+p+\hf), & p \geq 0,
\end{cases}
\]
and the functions $F_n^p$ are defined similarly for $\rho$ in the discrete part of the support of $\mu_1^p$, where $\mu_1$ denotes the corresponding orthonormality measure. Explicitly, if $k_1$ and $k_2$ are such that the measure is absolutely continuous, $\mu_1^p(\cdot;k_1,k_2)$ is given by
\begin{gather*}
\int f(\rho)\, d\mu_1^p(\rho;k_1,k_2) = \int_0^\infty f(\rho)\, w_1(\rho;k_1,k_2)\, d\rho\\
w_1(\rho;k_1,k_2) =\frac{1}{2\pi}\frac{ \Ga(k_1-k_2+\hf \pm i\rho) \Ga(k_1+k_2-\hf \pm i\rho) \Ga(k_2-k_1+\hf \pm i\rho)}{\Ga(2k_1)\Ga(2k_2)\Ga(\pm 2i\rho)},
\end{gather*}
for all $p \in \Z$. We note that in general the number of discrete mass points of $\mu_1^p$ may still depend on the value of $p$, and that this number is finite for all $p$. Since the weights of the measures $\mu^p_1$ do not depend on $p$, we can define the measure $\mu_1$ to be the measure with support
\[
\bigcup_{p \in \Z} \text{supp}\big(\mu_1^p(\cdot;k_1,k_2) \big),
\]
and with the same weight function as all the other measures $\mu_1^p$.
\begin{rem}
In order to make the situation as simple as possible we assume without much loss of generality that $k_1$ and $k_2$ are such that the measure $\mu_1^p$ does not have discrete mass points. In this case the measures $\mu_1^p$ are the same for all $p$, so $\mu_1^p(\cdot;k_1,k_2)=\mu_1(\cdot;k_1,k_2)$.
\end{rem}
\begin{prop} \label{prop I1}
For $p \in \Z$ the operator $I_2: \mathcal H_1^p \rightarrow \mathcal L\big(\mu_1(\rho;k_1,k_2)\big)$ defined by
\[
I_1: f_n^p \mapsto \Big(\rho \mapsto F_n^p(\rho;k_1,k_2) \Big),
\] 
is unitary and intertwines $\pi^+_{k_1} \tensor \pi^-_{k_2}( \De(\Om))$  with $M_{\rho^2+\frac14}$. 
\end{prop}
Observe that for $\rho \geq 0$, $\rho^2+\frac14$ corresponds to the action of the Casimir in the principal unitary series.
\begin{rem}
Applying the operator $I_1$ to $e_n \tensor e_m$ gives $F_n^{n-m}$ or $F_m^{n-m}$, depending on the sign of $n-m$. To simplify notations, we use both expressions regardless of the sign of $n-m$. This is justified by the following  identity for continuous dual Hahn polynomials, cf.~ \cite[(3.13)]{GK02},
\[
\begin{split}
S_n(\rho;&k_1-k_2+\frac{1}{2},k_1+k_2-\frac{1}{2}, k_2-k_1-p+\frac{1}{2}) = \\
&(-1)^p\, (k_1-k_2+\frac{1}{2}\pm i\rho)_p\,
S_{n-p}(\rho;k_2-k_1+\frac{1}{2},k_1+k_2-\frac{1}{2},
k_1-k_2+p-\frac{1}{2}), \qquad p \geq 0.
\end{split}
\]
\end{rem}
From Proposition \ref{prop I1} we find that
\[
v_1^p(\rho) = \sum_{n=0}^\infty F_n^p(\rho;k_1,k_2)\, f_n^p,
\]
is a generalized eigenvector of $\pi^+_{k_1} \tensor \pi^-_{k_2}( \De(\Om))$ for eigenvalue $\rho^2+\frac14$. We can also determine the actions of $H, E, F$ on $v_1^p(\rho)$, but, since $v_1^p(\rho) \not\in \ell^2(\Z_{\geq 0})$, this will only give formal expressions.

Let $\{e_n\}_{n \in \Z}$ denote the standard orthonormal basis of $\ell^2(\Z)$. We define the unitary operator $I_1'$ by
\[
\begin{split}
I_1': \ell^2(\Z_{\geq 0})^{\tensor 2} &\rightarrow \dirint \ell^2(\Z)\, d\mu_1(\rho;k_1,k_2)\\
f_n^p &\mapsto \int_0^\infty F_n^p(\rho;k_1,k_2)\, e_p\, d\rho
\end{split}
\]
Then, under the assumption that $\mu_1^p$ does not have discrete mass points, we have the following proposition.
\begin{prop} \label{prop:+-}
The operator $I_1'$ intertwines $\pi^+_{k_1} \tensor \pi^-_{k_2}$ with $\dirint \pi^P_{\rho,k_1-k_2} d\rho$.
\end{prop}
Here we use that, for $\eps \in [p, p+1)$ for some $p \in \Z$, the representation $\pi^P_{\rho,\eps}$ is equivalent to the representation $\pi^P_{\rho,\eps-p}$. 

To proof the proposition we only need to check that $I_1'$ indeed intertwines the actions of the generators $H$, $E$ and $F$. The action of $H$ follows from a straightforward computation. The actions of $E$ and $F$ follow from expressing the continuous dual Hahn polynomials as $_3F_2$-series, and using the contiguous relations \eqref{eq:cont.rel}. We refer to \cite[Thm.2.2]{GK02} for a detailed proof.

In case the measure $\mu_1(\cdot;k_1,k_2)$ has discrete mass points, these discrete terms correspond to a finite number of discrete series representations, or exactly one complementary series, cf.~\cite[Thm.2.2]{GK02}.

\subsection{Clebsch-Gordan coefficients for $\pi^+ \tensor \pi^P$}
Next we decompose the tensor product $\pi^+ \tensor \pi^P$ and we compute the corresponding Clebsch-Gordan coefficients. This can be done in the same way as the tensor product $\pi^+ \tensor \pi^-$, i.e.~by diagonalizing the action of the Casimir.

We define for fixed $p \in \Z$ the vector $g_n^p \in \ell^2(\Z_{\geq 0}) \tensor \ell^2(\Z)$ by
\begin{equation} \label{def gn}
g_n^p = e_n \tensor e_{p-n} ,
\end{equation}
and we define the subspace $\mathcal H_2^p$ by
\[
\mathcal H_2^p = \overline{ \Csp \{ g_n^p \ | \ n \in \Z_{\geq 0} \}} \cong \ell^2(\Z_{\geq 0}).
\]
From \eqref{eq De Om}, \eqref{def:pos} and \eqref{def:princ} it follows that the action of the Casimir element $\Om$ in the tensor product $\pi^+_k \tensor \pi^P_{\rho,\eps}$ acting on elements $g_n^p$ extends to an unbounded Jacobi operator acting on $\mathcal H_2^p$. This Jacobi operator corresponds to the recurrence relation for continuous dual Hahn polynomials with parameters
\[
\{a,b,c\} = \{ k+i\rho, k-i\rho, \hf-k-p-\eps \}.
\]
As before, we split the orthonormality measure for these polynomials in a $p$-independent part and a $p$-dependent part, and we glue the square root of the $p$-dependent part to the continuous dual Hahn polynomial, i.e.~for $\tau \in [0,\infty)$, we define
\begin{equation} \label{def:Gn}
G_n^p(\tau;k,\rho,\eps) = 
\begin{cases}
\displaystyle \sqrt{ \frac{ (\hf+\eps \pm i\rho)_{p}}{(\hf+k+\eps \pm i\tau)_{p} }} \ s_n(\tau; k+i\rho, k-i\rho, \hf-k-p-\eps), & \tau\in [0,\infty),\\ \\
\displaystyle \sqrt{ \frac{ (\hf+\eps \pm i\rho)_p }{ (j+1)_p(2k+2\eps-j)_p } } \ s_n(\tau_j; k+i\rho, k-i\rho, \hf-k-p-\eps), & \tau=\tau_j \in i\R_{<0}.
\end{cases}
\end{equation}
Here we denote $\tau_j=i(\hf-k-\eps+j)$, where $j \in \Z$ is such that $\tau_j \in i\R_{<0}$. Observe that it follows from the explicit expression as a $_3F_2$-function, that $\rho\mapsto G_n^p(\tau;k,\rho,\eps)$ is a continuous function for $\rho \in [0,\infty)$.  For fixed $p \in \Z$, the functions $G_n^p$ form an orthonormal basis for the Hilbert space $\mathcal L\big(\mu_2^p(\tau;k,\rho,\eps)\big)$, where $\mu_2^p$ is the measure defined by
\[
\int f(\tau)\, d\mu_2^p(\tau;k,\rho,\eps) = \int_0^\infty f(\tau) w_2(\tau;k,\rho,\eps)\, d\tau + \sum_{\substack{j \in \Z_{\geq -p}  \\ \tau_{j}<0}} f(\tau_j)\, w_{2}(\tau_j;k,\rho,\eps),
\]
where 
\begin{align*}
w_2(\tau;k,\rho,\eps) =& \frac1{2\pi} \frac{ \Ga(k+i\rho \pm i\tau) \Ga(k-i\rho \pm i\tau) \Ga(\hf-k-\eps\pm i\tau)}{ \Ga(2k) \Ga(\hf-\eps \pm i\rho) \Ga(\pm 2i\tau)},& \tau \in [0,\infty),\\
w_{2}(\tau;k,\rho,\eps) =& \frac{\Ga(2k+\eps-\hf \pm i\rho)}{ \Ga( 2k+2\eps-1) \Ga(2k)} \frac{ (-1)^j (1-2k-2\eps)_j (\frac32-k-\eps)_j (\hf-\eps\pm i\rho)_j }{ j!\, (\frac32-2k-\eps \pm i\rho)_j (\hf-k-\eps)_j },& \tau=\tau_j.
\end{align*}
\begin{rem} \label{rem:inf discr}
The number of discrete mass points of $\mu_2^p$ tends to infinity as $p \rightarrow \infty$.  
\end{rem}
Let $\mu_2(\cdot;k,\rho,\eps)$ be the measure correspond to the above defined weight functions $w_2$ and with support
\[
\bigcup_{p \in \Z}\text{supp}\big(\mu_2^p(\cdot;k,\rho,\eps)\big).
\]
So, by Remark \ref{rem:inf discr}, this measure has an infinite number of discrete mass points. 
\begin{prop} \label{prop I2}
For $p \in \Z$ the operator $I_2 : \mathcal H_2^p \rightarrow \mathcal L\big(\mu_2^p(\tau ;k,\rho,\eps)\big)$ defined by
\[
I_2 : g_n^p \mapsto \Big(\tau \mapsto G_n^p(\tau;k,\rho,\eps)\Big),
\]
is unitary and intertwines $\pi^+_k \tensor \pi^P_{\rho,\eps}(\De(\Om))$ with $M_{\tau^2+\frac14}$.
\end{prop}
So we find that
\[
v_2^p(\tau) = \sum_{n=0}^\infty G_n^p(\tau;k,\rho,\eps)\, g_n^p
\]
is a generalized eigenvector of $\pi^+_k \tensor \pi^P_{\rho,\eps}(\De(\Om))$ for eigenvalue $\tau^2+\frac14$. 
We define the unitary operator $I_2'$ by
\[
\begin{split}
I_2': \ell^2(\Z_{\geq 0}) \tensor \ell^2(\Z) &\rightarrow \dirint \ell^2(\Z)\, w(\tau;k,\rho,\eps) d\tau \oplus \bigoplus_j \ell^2(\Z_{\geq 0})\, w(\tau_j;k,\rho,\eps),\\
g_n^p &\mapsto \int_0^\infty G_n^p(\tau;k,\rho,\eps)\, e_p\, d\tau + \sum_{\substack{j \in \Z_{\geq -p}\\ \tau_{j}<0}} G_n^p(\tau_j;k,\rho,\eps) \, e_{p-j}.
\end{split}
\]
Then we have the following proposition.
\begin{prop} \label{prop:+P}
The operator $I_2'$ intertwines $\pi^+_{k} \tensor \pi^P_{\rho,\eps}$ with $\displaystyle \dirint \pi^P_{\tau,\eps'} d\tau \oplus \bigoplus_{j} \pi^+_{\eps'+j}$, where $\eps' = \eps+k$ and the direct sum is over all $j \in \Z$ such that $\eps'+j>0$.
\end{prop}
The proposition is proved by checking that $I_2'$ intertwines the actions of the generators, using the contiguous relations \eqref{eq:cont.rel}. The explicit expression for the Clebsch-Gordan coefficients as $_3F_2$-series can also be found in \cite{KV99}, but the continuous dual Hahn polynomials are not mentioned there.

We remark that the tensor product $\pi^+_k \tensor \pi^C_{\la,\eps}$ can be treated in completely the same way as $\pi^+_k \tensor \pi^P_{\rho,\eps}$. Formally, the Clebsch-Gordan coefficients can be obtained from the substitution $\rho \mapsto -i(\la+\hf)$.

\subsection{Clebsch-Gordan coefficients for $\pi^P \tensor \pi^-$} \label{ssec:CGC P-}
Let $\vartheta :\su(1,1) \rightarrow \su(1,1)$ be the involution defined on the generators of $\su(1,1)$ by
\[
\vartheta(H) = -H, \quad \vartheta(E) = F, \quad \vartheta(F)=E.
\]
This involution can be extended to $\mathcal U(\su(1,1))$ as an algebra morphism. The Casimir element satisfies $\vartheta(\Om) = \Om$. If we identify the representation space for $ \pi^+_k$ with the representation space for $\pi^-_k$, we have $\pi^+_k\big(\vartheta(X)\big) = \pi^-_k(X)$ for $X \in \U(\su(1,1))$.  Furthermore, define a unitary operator $U: \ell^2(\Z) \rightarrow \ell^2(\Z)$ by $Ue_n =  (-1)^n e_{-n}$. Using \eqref{def:princ} it is easy to check that $U \circ \pi^P_{\rho,-\eps}\circ \vartheta=\pi^P_{\rho,\eps}\circ U$, so $\pi^P_{\rho,-\eps}\circ \vartheta$ is equivalent to $\pi^P_{\rho,\eps}$.

For $p\in \Z$, we define
\begin{equation} \label{def:ln}
l_n^p = e_{n+p} \tensor e_n \in \ell^2(\Z) \tensor \ell^2(\Z_{\geq 0}),
\end{equation}
or equivalently $l_n^p = (-1)^{n+p}\, \ Ue_{-p-n} \tensor e_n$. Compare this last expression with \eqref{def gn}. From the above considerations we obtain
\[
\pi^P_{\rho,\eps} \tensor \pi^-_k (\De(X))\,  l_n^p =(-1)^{n+p}\, \pi^P_{\rho,-\eps} \tensor  \pi^+_k \big(\De(\vartheta(X))\big)\, \tilde g_{n}^{-p}, \qquad X \in \U(\su(1,1)),
\]
where $\tilde g_n^p$ is the vector $g_n^p$ with the factors in the tensor product interchanged.
So we can obtain the decomposition of $\pi^P_{\rho,\eps}\tensor \pi^-_k $ from the decomposition of $\pi^+_k\tensor \pi^P_{\rho,-\eps}$ (recall that the coproduct $\De$ is cocommutative). Since $\vartheta(\Om)=\Om$ and since the Clebsch-Gordan coefficients are eigenfunctions of $\Om$ in the tensor product, we obtain from Proposition \ref{prop I2} that the Clebsch-Gordan coefficients for the tensor product $\pi^P_{\rho,\eps} \tensor \pi^-_k $ are the functions $(-1)^{n+p} G_{n}^{-p}(\si;k,\rho,-\eps)$. From the tensor product decomposition for $\pi^+_k \tensor \pi^P_{\rho,\eps}$ we obtain
\[
\pi^P_{\rho,\eps} \tensor \pi^-_k \cong \dirint \pi^P_{\si,\eps'} d\si \oplus \bigoplus_{j} \pi^-_{j-\eps'}, \qquad \eps' = \eps-k,
\]
where the direct sum is over all $j \in \Z$ such that $j-\eps' >0$. \\

Finally let us mention that the tensor product of two principal unitary series can be decomposed as a direct integral over the principal unitary series, plus an infinite sum of discrete series. However, the principal unitary series occur with multiplicity $2$. For this reason we avoid tensor products of two principal unitary series in this paper. See for instance \cite{KV98} for related Clebsch-Gordan coefficients.


\section{Racah coefficients for $\pi^+ \tensor \pi^- \tensor \pi^+$} \label{sec:+-+}
In this section we find explicit expresssions for the Racah coefficients related to the tensor product $\pi^+ \tensor \pi^- \tensor \pi^+$, and we relate the corresponding generalized orthogonality relations to the Wilson function transform of type I.

\subsection{Racah coefficients} \label{ssec:Racah1}
We consider the tensor product representation $\pi^+_{k_1} \tensor \pi^-_{k_2} \tensor \pi^+_{k_3}$, and we want to find the corresponding Racah coefficients.  
First we consider the action of the Casimir $\Om$ in the tensor product $(\pi^+_{k_1} \tensor \pi^-_{k_2}) \tensor \pi^+_{k_3}$, and we want to diagonalize this action. We find that the subspace of $\ell^2(\Z_{\geq 0})^{\tensor 3}$ spanned by elements $e_{n_1} \tensor e_{n_2} \tensor e_{n_3}$ for which $n_1-n_2+n_3$ is constant, is invariant under the action of the Casimir. Therefore we define, for $r \in \Z$,
\[
k_{n,m}^r =f_n^{r-m} \tensor e_m \in \ell^2(\Z_{\geq 0})^{\tensor 3},
\]
where $f_n^p$ is defined by \eqref{def fn}. Then, for $r \in \Z$, the subspace
\[
\mathcal A^r =\overline{ \Csp \{ k_{n,m}^r \ | n,m \in \Z_{\geq 0} \} } \subset \ell^2(\Z_{\geq 0})^{\tensor 3},
\]
is  invariant under the action of $\Om$.  

We define
\[
K_{n,m}^r(\tau,\rho;k_1,k_2,k_3) = F_n^{r-m}(\rho;k_1,k_2) G_m^{r}(\tau;k_3,\rho,k_1-k_2) .
\]
Let $m^r$ be the product measure defined by
\[
\iint f(\tau,\rho)\, dm^r(\tau,\rho;k_1,k_2,k_3) = \iint f(\tau,\rho)\, d\big(\mu_1^{r-m}(\rho;k_1,k_2) \times \mu_2^r(\tau;k_3,\rho,k_1-k_2)\big).
\]
The set $\{K_{n,m}^r\}_{n,m \in \Z_{\geq 0}}$ is an orthonormal basis for $\mathcal H\big(m^r(x,y;k_1,k_2,k_3)\big)$, the Hilbert space consisting of functions in two variables, even in both variables, with inner product associated to the measure $m^r$.
\begin{prop} \label{prop J1}
For $r\in \Z$ the operator $J_1 : \mathcal A^r \rightarrow \mathcal H\big(m^r(\tau,\rho;k_1,k_2,k_3)\big)$ defined by
\[
J_1 : k_{n,m}^r \rightarrow \Big((\tau,\rho) \mapsto K_{n,m}^r(\tau,\rho;k_1,k_2,k_3) \Big), 
\] 
is unitary and intertwines $\pi^+_{k_1} \tensor \pi^-_{k_2} \tensor \pi^+_{k_3} (( \De \tensor 1) \circ \De(\Om))$ with $M_{\tau^2+\frac14}$.
\end{prop}
\begin{proof}
First observe that $J_1$ is unitary since it maps an orthonormal basis onto another orthonormal basis. 

For the second statement we first apply the operator $I_1'$ from Proposition \ref{prop:+-} in the first and second factor of the tensor product:
\[
\begin{split}
I_1' \tensor \mathrm{id} :  k_{n,m}^r= f_{n}^{r-m}\tensor e_m \mapsto \int_0^\infty F_n^{r-m}(\rho;k_1,k_2) e_{r-m}\tensor e_m\, d\rho.
\end{split}
\]
Here $\mathrm{id}$ denotes the identity operator on $\ell^2(\Z_{\geq 0})$.
From the intertwining property of $I_2'$, see Proposition \ref{prop:+-}, we find
\[
I_2'\tensor \mathrm{id}\ \circ\  \pi^+_{k_1} \tensor \pi^-_{k_2} \tensor \pi_{k_3}^+\big(( \De \tensor 1) \circ \De(\Om)\big) =\dirint \pi^P_{\rho,k_1-k_2} \tensor \pi^+_{k_3}(\De(\Om)) \, d\rho.
\]
Next we apply $I_2$ from Proposition \ref{prop I2}, then we find that 
\[
\dirint J_1\, d\rho : k_{n,m}^r \mapsto \Bigg( \tau \mapsto \int_0^\infty F_n^{r-m}(\rho;k_1,k_2) G_m^r(\tau;k_3,\rho,k_1-k_2)\, d\rho \Bigg)
\]
intertwines $\pi^+_{k_1} \tensor \pi^-_{k_2} \tensor \pi^+_{k_3} (( \De \tensor 1) \circ \De(\Om))$ with $M_{\tau^2+\frac14}$. Clearly $M_{\tau^2+\frac14}$ is independent of $\rho$, so the operator $J_1$ also intertwines the action of $\Om$ with $M_{\tau^2+\frac14}$, for almost every $\rho$. Since the function $K_{n,m}^r(\tau,\rho)$ is continuous in $\rho$, the proposition follows.
\end{proof}
From Proposition \ref{prop J1} it follows that 
\[
w_1^r(\rho,\tau) = \sum_{m=0}^\infty G_m^{r}(\tau;k_3,\rho,k_1-k_2)\, v_2^{m-r}(\rho) \tensor e_m = \sum_{n,m=0}^\infty K_{n,m}^r(\tau,\rho;k_1,k_2,k_3)\, k_{n,m}^r
\]
is a generalized eigenvector of $\pi^+_{k_1} \tensor \pi^-_{k_2} \tensor \pi^+_{k_3} (( \De \tensor 1) \circ \De(\Om))$ for eigenvalue $\tau^2+\frac14$.

Combining Propositions \ref{prop I1} and \ref{prop I2} we can construct an explicit intertwiner $J_1'$ to decompose $\pi^+_{k_1} \tensor \pi^-_{k_2} \tensor \pi^+_{k_3}$ into irreducible representations. So $J_1'$ intertwines $\pi^+_{k_1} \tensor \pi^-_{k_2} \tensor \pi^+_{k_3}$ with 
\[
\dirint \left( \dirint \pi^P_{\tau,\eps'} d\tau \oplus \bigoplus_j \pi^+_{\eps'+j} \right) d\rho,
\]
where $\eps'= k_1-k_2+k_3$. We abuse notation and write the intertwiner as
\[
J_1': k_{n,m}^r \mapsto \iint K_{n,m}^r(\tau,\rho;k_1,k_2,k_3)\, e_p\, d\tau d\rho ,
\]
So $e_p$ inside the integral is a basis vector of $\ell^2(\Z)$ for $\tau \in [0,\infty)$ and $e_p$ should be replaced by $e_{j-p} \in \ell^2(\Z_{\geq 0})$ when $\tau$ is in the discrete part of the support of the measure $dm$.\\

Next we consider the action of the Casimir in tensor product $\pi^+_{k_1} \tensor (\pi^-_{k_2} \tensor \pi^+_{k_3})$, and we want to have a unitary operator $J_2$ similar to $J_1$ in Proposition \ref{prop J1}. Recall that the coproduct $\De$ is cocommutative, so we can use Proposition \ref{prop:+-} to decompose $\pi^-_{k_2} \tensor \pi^+_{k_3}$. Then we can find $J_2$ in the same way as above by considering the action of $\Om$ on the element $k_{n,m}^r \in \mathcal A^r$, where we write $k_{n,m}^r$ as
\[
k_{n,m}^r = e_{n} \tensor \tilde f_m^{r-n},
\]
and $\tilde f_n^p$ is the same as $f_n^p$ defined by \eqref{def fn}, but with the factors interchanged, i.e.
\[
\tilde f_n^p = 
\begin{cases}
e_{n-p} \tensor e_{n}, & p \leq 0, \\
e_{n} \tensor e_{n+p}, & p \geq 0.
\end{cases}
\]
By the symmetry in $k_1$ and $k_3$ for the tensor product $\pi^+_{k_1} \tensor \pi^-_{k_2} \tensor \pi^+_{k_3}$, we can find $J_2$ directly from Proposition \ref{prop J1} by interchanging $k_1$ and $k_3$, and $n$ and $m$.
\begin{prop} \label{prop J2}
Let $r\in \Z$ the operator $J_2 : \mathcal A^r \rightarrow \mathcal H\big( m^r(\tau,\si;k_3,k_2,k_1)\big)$ defined by
\[
J_2 : k_{n,m}^r \rightarrow \Big((\tau,\si) \mapsto K_{m,n}^r(\tau,\si;k_3,k_2,k_1) \Big)
\] 
is unitary and intertwines $\pi^+_{k_1} \tensor \pi^-_{k_2} \tensor \pi^+_{k_3} (( 1 \tensor \De) \circ \De(\Om))$ with $M_{\tau^2+\frac14}$.
\end{prop}
So
\[
w_2^r(\si,\tau) = \sum_{n,m}^\infty K_{m,n}^r(\tau,\si;k_3,k_2,k_1) \, k_{n,m}^r
\]
is a generalized eigenvector of $\pi^+_{k_1} \tensor \pi^-_{k_2} \tensor \pi^+_{k_3} (( 1 \tensor \De) \circ \De(\Om))$ for eigenvalue $\tau^2+\frac14$. Again we can make an intertwining operotar $J_2'$ similarly to $J_1'$.\\

From the coassociativity of the coproduct it follows that there exists a unitary operator that intertwines $\pi^+_{k_1} \tensor \pi^-_{k_2} \tensor \pi^+_{k_3} \circ ( 1 \tensor \De) \circ \De$ with $\pi^+_{k_1} \tensor \pi^-_{k_2} \tensor \pi^+_{k_3} \circ ( \De \tensor 1) \circ \De$. Then it follows from Propositions \ref{prop J1} and \ref{prop J2} that there exists a unitary operator
\[
U : \mathcal H\big(m^r(\tau,\rho;k_1,k_2,k_3)\big) \rightarrow \mathcal H\big( m^r(\tau',\si;k_3,k_2,k_1)\big),
\]
such that $U \circ J_1 = J_2$. We write $U$ as an integral operator,
\[
(Uf)(\tau',\si) = \iint f(\tau,\rho)\,U_{\tau',\si}^r(\tau,\rho;k_1,k_2,k_3)\,  dm^r(\tau,\rho;k_1,k_2,k_3),
\]
for some kernel $U_{\tau',\si}^r(\tau,\rho;k_1,k_2,k_3)$. We call this kernel the Racah coefficient for the tensor product $\pi^+_{k_1} \tensor \pi_{k_2}^- \tensor \pi^+_{k_3}$. We justify this terminology in Remark \ref{rem:Racah}. We will determine an explicit expression for the Racah coefficient, using 
\begin{equation} \label{eq:U1}
K_{m,n}^r(\tau',\si;k_3,k_2,k_1) = \iint K_{n,m}^r(\tau,\rho;k_1,k_2,k_3) \, U_{\tau',\si}^r(\tau,\rho;k_1,k_2,k_3)\,  dm^r(\tau,\rho;k_1,k_2,k_3),
\end{equation}
which follows directly from $U \circ J_1 = J_2$. First we determine some useful properties of the Racah coefficient.

Since both $J_1$ and $J_2$ intertwine the action of $\Om$ with $M_{\tau^2+\frac14}$, $U$ leaves functions in the variable $\tau$ fixed, and this gives
\[
U_{\tau',\si}^r(\tau,\rho;k_1,k_2,k_3) = 0, \qquad \text{if }\tau \neq \tau'.
\]
Let $\nu^r(\rho;k_1,k_2,k_3,\tau)$ be the measure $m^r(\tau,\rho;k_1,k_2,k_3)$ restricted to a line on which $\tau=$ constant, i.e.
\[
\int f(\rho)\, d\nu(\rho;k_1,k_2,k_3,\tau) = \int f(\rho)\, w_1(\rho,k_1,k_2) w_2(\tau;k_3,\rho,k_1-k_2) \, d\rho.
\]
For fixed $\tau$ this measure $d\nu^r$ does not depend on $r$, therefore we will omit the superscript $r$. Then, instead of using \eqref{eq:U1}, we can determine an explicit expression for the Racah coefficient from
\begin{equation} \label{eq:U2}
K_{m,n}^r(\tau,\si;k_3,k_2,k_1) = \int_0^\infty K_{n,m}^r(\tau,\rho;k_1,k_2,k_3) \, U_{\tau,\si}^r(\tau,\rho;k_1,k_2,k_3)\, d\nu(y;k_1,k_2,k_3,\tau).
\end{equation}
We use this identity for the following lemma.
\begin{lem} \label{lem:indep}
The Racah coefficient does not depend on $r$.
\end{lem}
\begin{proof}
The functions $K_{m,n}^r(\tau,\si;k_3,k_2,k_1)$ are the ``matrix elements"  of the unitary operator $J_2'$ that intertwines $\pi^+_{k_1} \tensor \pi^-_{k_2} \tensor \pi^+_{k_3}$ with $\iint^\oplus \pi^P_{\tau,\eps}\, d\tau d\si$, where $\eps= k_1-k_2+k_3$. So from the action of $E$ we find 
\begin{equation} \label{eq:rec.rel K}
\begin{split}
A\,& K_{m,n-1}^r(\tau,\si;k_3,k_2,k_1) + B\, K_{m-1,n}^r(\tau,\si;k_3,k_2,k_1)
+ C\, K_{m,n}^r(\tau,\si;k_3,k_2,k_1) \\
&= |r+\eps+\hf+ i\tau|\, K_{m,n}^{r+1}(\tau,\si;k_3,k_2,k_1),
\end{split}
\end{equation}
where we assume for simplicity that $\tau \in [0,\infty)$ and we denoted
\begin{gather*}
A = \sqrt{ n(2k_1+n)}, \quad B = \sqrt{m(2k_3+m-1)},\\  C = \sqrt{ (m+n-r)(2k_2+m+n-r-1)}.
\end{gather*}
Applying \eqref{eq:U2} this gives, 
\[
\begin{split}
&|r+\eps+\hf+ i\tau|\, \int_0^\infty U_{\tau,\si}^{r+1}(\tau,\rho;k_1,k_2,k_3) K_{n,m}^{r+1}(\tau,\rho;k_1,k_2,k_3)\, d\nu(\rho;k_1,k_2,k_3,\tau) \\ 
 =&\int_0^\infty \Big( A\, K_{n-1,m}^r(\tau,\rho;k_1,k_2,k_3) + B\, K_{n,m-1}^r(\tau,\rho;k_1,k_2,k_3)
+ C\, K_{n,m}^r(\tau,\rho;k_1,k_2,k_3) \Big)\\
& \qquad  \times  U^{r}_{\tau,\si}(\tau,\rho;k_1,k_2,k_3)\, d\nu(\rho;k_1,k_2,k_3,\tau) \\
= &|r+\eps+\hf+ i\tau|\, \int_0^\infty U_{\tau,\si}^{r}(\tau,\rho;k_1,k_2,k_3) K_{n,m}^{r+1}(\tau,\rho;k_1,k_2,k_3)\, d\nu(\rho;k_1,k_2,k_3,\tau).
\end{split}
\]
In the last step we applied \eqref{eq:rec.rel K} again with $(n,k_1) \leftrightarrow (m,k_3)$. From the above identity it follows that the Racah coefficient is independent of $r$.
\end{proof}
\begin{rem} \label{rem:Racah}
Let us explain why we call the kernel $U_{\tau',\si}^r(\tau,\rho;k_1,k_2,k_3)$ the Racah coefficient. Formally, the Racah coefficient can be defined in terms of the generalized eigenvectors by
\[
U'(\rho,\si;k_1,k_2,k_3,\tau) \de_{r_1 r_2} \de_{\tau,\tau'} = \langle w_1^{r_1}(\rho,\tau), w_2^{r_2}(\si,\tau') \rangle_{\ell^2(\Z_{\geq 0})^{\tensor 3}}.
\]
Compare this expression with \eqref{def:Racah coeff}. So the Racah coefficients can be considered as ``matrix elements" of the operator $U'$ defined by $U': w_1^{r}(\rho,\tau) \mapsto w_2^{r}(\si,\tau_2)$, and as before the coefficient is independent of $r$. If we ``expand" one eigenvector in terms of the other basis of eigenvectors, we have, for $\rho^2=y$,
\[
w_2^r(\si,\tau_2) = \int U'(\rho,\si;k_1,k_2,k_3,\tau) w_1^r(\rho,\tau) \, d\nu(\rho;k_1,k_2,k_3,\tau),
\]
as an identity for generalized eigenvectors. Now taking inner products with $k_{n,m}^r$ on both sides gives \eqref{eq:U2} for the kernel of $U'$. So the kernel of our operator $U$ is the same as the (formal) Racah coefficient.
\end{rem}

From here on we use the notation
\[
U(\rho,\si;k_1,k_2,k_3,\tau) = U_{\tau,\si}^r(\tau,\rho;k_1,k_2,k_3).
\]

We are now ready to determine an explicit expression for the Racah coefficient. To write down the result it is convenient to introduce the following function:
\begin{equation}\label{def:Wilson function}
\begin{split}
\phi_\la(x;a,b,c,d)
=& \frac{ \Ga(\tilde a+\tilde b+\tilde c+i\la) }{\Ga(a+b) \Ga(a+c) \Ga(1+a-d) \Ga(1-\tilde d -i\la)   \Ga(\tilde b+c+i\la \pm ix) }\\
& \times W(\tilde a+\tilde b+ \tilde c-1 +i\la;a+ix, a-ix, \tilde a+i\la, \tilde b+i\la, \tilde c+i\la ),
\end{split}
\end{equation} 
where Bailey's $W$ notation is used for a very-well-poised $_7F_6$-function, i.e.
\[
\begin{split}
W(a;b,c,d,e,f) &= \F{7}{6}{a,\ 1+\hf a,\ b,\ c,\ d,\ e,\ f\ }{\hf a, 1+a-b, 1+a-c, 1+a-d, 1+a-e, 1+a-f}{1}\\
&= \sum_{n=0}^\infty \frac{ a+2n}{a} \frac{ (a)_n (b)_n (c)_n (d)_n (e)_n (f)_n }{ n! (1+a-b)_n (1+a-c)_n (1+a-d)_n (1+a-e)_n (1+a-f)_n }
\end{split}
\]
The parameters $\tilde a, \tilde b, \tilde c, \tilde d$ are called the dual parameters, and they are defined by
\begin{equation} \label{eq:dual param}
\begin{split}
\tilde a &= \hf(a+b+c+d-1),\quad 
\tilde b=\hf(a+b-c-d+1) , \\
\tilde c&=\hf(a-b+c-d+1), \quad
\tilde d=\hf(a-b-c+d+1).
\end{split}
\end{equation}
The function $\phi_\la(x;a,b,c,d)$ is called the Wilson function. Another useful expression for the Wilson function is obtained from transforming the $_7F_6$-function into a sum of two balanced $_4F_3$-functions, see \cite[(3.5)]{Gr03},
\begin{equation} \label{eq:Wil4F3}
\begin{split}
\phi_{\la}(x;a,b,c,d) =& \frac{ \Ga(1-a-d) }{ \Ga(a+b) \Ga(a+c) \Ga(1-d \pm ix) \Ga(1-\tilde d \pm i\la) } \\
& \times \F{4}{3}{a+ix, a-ix, \tilde a+i\la, \tilde a-i\la}{a+b, a+c, a+d}{1} \\
&+ \frac{ \Ga(a+d-1)} { \Ga(1+b-d) \Ga(1+c-d) \Ga(a \pm ix) \Ga(\tilde a \pm i\la) } \\
& \times \F{4}{3}{1-d+ix, 1-d-ix, 1-\tilde d+i\la, 1-\tilde d-i\la}{1+b-d, 1+c-d, 2-a-d}.
\end{split}
\end{equation}
From this expression it follows that the Wilson function satisfies the duality property
\begin{equation} \label{eq:duality}
\phi_\la(x;a,b,c,d) = \phi_x(\la;\tilde a,\tilde b, \tilde c, \tilde d).
\end{equation}
Also we see from \eqref{eq:Wil4F3} that the Wilson function is invariant under $a \leftrightarrow 1-d$ and $b \leftrightarrow c$. For large $x$ we have \cite{Gr03}
\[
\phi_\la(x;a,b,c,d) \sim x^{d-a-b-c-2i\la}e^{\pi x}\Big( c_1 x^{-2i\la} +c_2 x^{2i\la} \Big),
\]
where $c_1, c_2$ are independent of $x$. 
We show that the Racah coefficient $U(\rho,\si;k_1,k_2,k_3,\tau)$ is a multiple of a Wilson function.
\begin{thm} \label{thm:U+-+}
The Racah coefficient $U(\rho,\si;k_1,k_2,k_3,\tau)$ can explicitly be written as
\[
\frac{2\pi\ \Ga(2k_1) \Ga(2k_2) \Ga(2k_3)\Ga(\pm 2i\tau) }{\Ga(k_2-k_1-k_3+\hf \pm i\tau) } \ \phi_\rho(\si; k_1 + i\tau, k_2+k_3-\hf, k_3-k_2+\hf, 1-k_1+i\tau).
\]
\end{thm}
For the proof we need the following lemma, which is Corollary 6.4 in \cite{Gr03}.
\begin{lem} \label{lem1}
The Wilson function satisfies:
\begin{align*}
\phi_\la(x;a,b,c,d) &= \\
\sum_{n=0}^\infty C_n(x,\la)& \F{3}{2}{-n, f+i\la,f-i\la}{\tilde b+f, \tilde c+f}{1} \F{3}{2}{ -n , f-\tilde a, f+\tilde d-1}{f+\tilde c-c-ix, f+\tilde c-c+ix}{1},
\end{align*}
where
\[
C_n(x,\la) = \frac{(f+\tilde c-c \pm ix)_n   \Ga(f \pm i\la) }{n!\, \Ga(1+a-d+n) \Ga(f+\tilde b) \Ga(f+\tilde c) \Ga(\tilde a \pm i\la) \Ga(1-\tilde d \pm i\la) } .
\]
\end{lem}
\begin{proof}[Proof of Theorem \ref{thm:U+-+}]
We use \eqref{eq:U2}, where we write the functions $K_{n,m}^r$ in terms of continuous dual Hahn polynomials and we put $m=r=0$. Using the orthogonality relations for the continuous dual Hahn polynomials $s_n(\rho;k_2-k_1+\hf, k_1+k_2-\hf, k_1-k_2+\hf)$ we obtain
\[
\begin{split}
V(\rho)= U(\rho,\si;k_1,k_2,k_3,\tau)&   \frac{ \Ga(k_3+i\tau \pm i\rho) \Ga(k_3-i\tau \pm i\rho)\Ga(k_2-k_1-k_3+\hf \pm i\tau) }{2\pi\ \Ga(2k_3)\Ga(k_2-k_1+\hf \pm i\rho) \Ga(\pm 2i\tau)}   \\
=\sum_{n=0}^\infty& \sqrt{\frac{ (k_2-k_3+r+\hf \pm i\si)_n  }{  n!\, (2k_2)_n } } \,
s_n(\tau;k_1+i\si, k_1-i\si, k_2-k_1-k_3+\hf)\\
&\times s_n(\rho;k_1-k_2+\hf, k_1+k_2-\hf, k_2-k_1+\hf)\\
=\sum_{n=0}^\infty& \frac{ (k_2-k_3+\hf+i\si)_{n} (k_2-k_3+\hf-i\si)_{n} }{ n!\, (2k_1)_{n} } \\
& \times  \F{3}{2}{ -n, k_2-k_1+\hf+i\rho, k_2-k_1+\hf-i\rho}{ 2k_2, 1 }{1} \\
& \times \F{3}{2}{ -n, k_2-k_1-k_3+\hf +i\tau, k_2-k_1-k_3+\hf -i\tau}{ k_2-k_3+\hf+i\si, k_2-k_3+\hf-i\si}{1},
\end{split}
\]
provided that this function $V$ is an element of $\mathcal H\big(\mu_1^0(\rho;k_1,k_2)\big)$. Here we assume for simplicity that $\tau \in [0,\infty)$. In case $\tau$ is in the discrete part of the support of the corresponding measure, the proof runs along the same lines.
Now we apply Lemma \ref{lem1} with parameters
\[
(a,b,c,d,f,x,\la) \mapsto (k_1+i\tau,k_2+k_3-\hf,k_3-k_2+\hf, 1-k_1+i\tau, k_2-k_1+\hf,\si,\rho),
\]
then we find the expression given in the theorem for the Racah coefficient. Using Stirling's asymptotic formula for the $\Ga$-function and the asymptotic behaviour of the Wilson function, we verify that indeed $V \in \mathcal H\big(\mu_1^0(\rho;k_1,k_2)\big)$. Now the theorem follows from the fact that  a function in $\mathcal H\big(\mu_1^0(\rho;k_1,k_2)\big)$ is uniquely determined by the coefficients in its expansion in the continuous dual Hahn polynomials $s_n(\rho;k_2-k_1+\hf, k_1+k_2-\hf, k_1-k_2+\hf)$.
\end{proof}
\begin{rem}
We assumed for simplicity that $k_1,k_2,k_3$ are such that discrete terms do not appear in the decompositions of $\pi^+_{k_i} \tensor \pi_{k_2}^-$, $i=1,3$. If there are discrete terms in the decomposition, the Racah coefficients can be determined in exactly the same way, only discrete terms should be added to the integrals. In case the discrete terms correspond to positive discrete series, the Clebsch-Gordan coefficients for the tensor product of two positive discrete series are also needed here.
\end{rem}

\subsection{The Wilson function transform of type I} \label{ssec:WilsonI}
Theorem \ref{thm:U+-+} gives a natural interpretation of the Wilson functions in the representation theory of the Lie algebra $\su(1,1)$. In \cite{Gr03} the Wilson functions are studied as eigenfunctions of a certain difference operator.
Spectral analysis of the difference operator on a specific Hilbert space leads to a unitary integral transform that has the Wilson as a kernel. In this section we show how the integral transform, called the Wilson function transform of type I in \cite{Gr03}, is related to the interpretation of the Wilson functions as Racah coefficients.\\

Let the parameters $a,b,c,d \in \C$ satisfy $a,b,c,1-d>0$, or $b,c>0$, $a=1-\overline{d}$ and $\Re(a)>0$. We define the Hilbert space $\mathcal M=\mathcal M(a,b,c,d)$ to be the Hilbert space consisting of even functions that have finite norm with respect to the inner product 
\[
\inprod{f}{g}_{\mathcal M} = \int_0^\infty f(x) \overline{g(x)} \, dm(x;a,b,c,d),
\]
where $m$ is the measure given by
\[
\int_0^\infty f(x)\, dm(x;a,b,c,d) = \frac{1}{2\pi} \int_0^\infty f(x) \frac{ \Ga(a\pm ix) \Ga(b \pm ix) \Ga(c \pm ix) \Ga(1-d \pm ix) }{ \Ga( \pm 2ix) }dx.
\]
We define an operator $\mathcal F=\mathcal F(a,b,c,d)$ on $\mathcal M$ by 
\[
(\mathcal F f)(\la) = \inprod{f}{\phi_\la(\cdot;a,b,c,d)}_{\mathcal M}.
\]
This operator is initially defined on the domain consisting of functions $f \in \mathcal M(a,b,c,d)$ such that the integral $\inprod{f}{\phi_\la}_{\mathcal M}$ converges. We call this integral operator $\mathcal F$ the Wilson function transform of type I, or just the Wilson transform I. The following theorem is Theorem 4.12 in \cite{Gr03}.
\begin{thm}
The operator $\mathcal F(a,b,c,d) : \mathcal M(a,b,c,d) \rightarrow \mathcal M(\tilde a, \tilde b, \tilde c, \tilde d) $ extends to a unitary operator, and the inverse is given by $\mathcal F(\tilde a, \tilde b, \tilde c, \tilde d)$.
\end{thm}
We sketch a proof of this theorem using the current interpration of the Wilson functions as Racah coefficients.

Let us define the integral transform $\mathcal I=\mathcal I(k_1,k_2,k_3,\tau)$ by
\[
\begin{split}
\mathcal I : \mathcal L\big(\nu(\rho;k_1,k_2,k_3,\tau)\big) \rightarrow \mathcal L\big(\nu(\si;k_3,k_2,k_1,\tau)\big) \\
(\mathcal I f)(\si) = \int_0^\infty f(\rho)\, U(\rho,\si;k_1,k_2,k_3,\tau)\, d\nu(\rho;k_1,k_2,k_3,\tau).
\end{split}
\]
Then from comparing squared norms of the funtions $K_{n,m}^r$ it follows from \eqref{eq:U2} that $\mathcal I$ is a partial isometry.
By the following lemma we find that $\mathcal I$ extends to a unitary operator.
\begin{lem}
The set of functions $\{ K_{n,m}^r(\rho;k_1,k_2,k_3,\tau)\ |\ n,m \in \Z_{\geq 0},\, r\in \Z \}$ spans the space $\mathcal L\big(\nu(\rho;k_1,k_2,k_3,\tau)\big)$.
\end{lem}
\begin{proof}
Let $f \in \mathcal L\big(\nu(\rho;k_1,k_2,k_3,\tau)\big)$ be a function such that 
\[
\int_0^\infty f(\rho) K_{n,m}^r(\tau,\rho;k_1,k_2,k_3)\, d\nu(\rho;k_1,k_2,k_3,\tau) =0,
\]
for all $n,m,r$. We integrate over all $\tau$, then by Fubini's theorem we have
\[
\iint f(\rho) K_{n,m}^r(\tau,\rho;k_1,k_2,k_3)\, dm^r(\tau,\rho;k_1,k_2,k_3)=0,
\]
for all $n,m,r$. Since the functions $K_{n,m}^r$, $n,m \in \Z_{\geq}$ form a basis for $\mathcal H\big( m^r(x,y;k_1,k_2,k_3))$ we find that $f=0$ almost everywhere.
\end{proof}

Writing out the weight functions for the measure $\nu(\rho;k_1,k_3,k_3,\tau)$, we see that this measure is, up to a multiplicative constant, the same as the measure $m(\rho;a,b,c,d)$ for the Wilson transform I, where the parameters $a,b,c,d$ (the Wilson parameters) are given by
\[
(a,b,c,d) = (k_1 + i\tau, k_2+k_3-\hf, k_3-k_2+\hf, 1-k_1+i\tau).
\]
So the above defined integral transform $\mathcal I$ is the Wilson transform I. 

Note that the variables $\rho$ and $\si$ came from decomposing the tensor product $\pi^+ \tensor \pi^-$ into a direct integral over the principal unitary series. So the support of the measure $m$ for the Wilson transform I corresponds to the decomposition of $\pi^+ \tensor \pi^-$ into irreducible representations. Discrete terms can appear in the decomposition, corresponding to discrete series or exactly one complementary series, and in this case the measure $m$ for the Wilson transform I has corresponding discrete mass points. These are exactly the cases that one of the parameters $a,b,c,1-d$ is non-positive.

To define the Racah coefficients we might as well have started with the inverse of \eqref{eq:U2}
\[
K_{n,m}^r(\tau,\rho;k_1,k_2,k_3) = \int_0^\infty K_{m,n}^r(\tau,\si;k_3,k_2,k_1) U(\rho,\si;k_1,k_2,k_3,\tau)\,d\nu(\si;k_3,k_2,k_1,\tau).
\]
Now interchanging $k_1 \leftrightarrow k_3$, $\rho \leftrightarrow \si$, $n \leftrightarrow m$, and comparing the obtained expression with \eqref{eq:U2} gives the following property for the Racah coefficients:
\[
U(\rho,\si;k_1,k_2,k_3,\tau) = U(\si,\rho;k_3,k_2,k_1,\tau).
\]
From \eqref{eq:dual param} we see that interchanging $k_1$ and $k_3$ is the same as interchanging the Wilson parameters with the dual Wilson parameters. So the duality property \eqref{eq:duality} for the Wilson function can be obtained as a consequence of the associativity of the tensor product $\pi^+_{k_1} \tensor \pi^-_{k_2} \tensor \pi^+_{k_3}$. This duality property can be considered as the main reason why the transform $\mathcal F$ is self-dual, up to an involution on the parameters.


\section{Racah coefficients for ${\pi^+ \tensor \pi^P \tensor \pi^-}$} \label{sec:+P-}
In the same way as in the previous section we calculate the Racah coefficients for the tensor product representation ${\pi^+_{k_1} \tensor \pi^P_{\rho,\eps} \tensor \pi^-_{k_3}}$. Also we relate the generalized orthogonality relation for the Racah coefficients to the Wilson function transform of type II.

\subsection{Racah coefficients}
The subspace of $\ell^2(\Z_{\geq 0}) \tensor \ell^2(\Z) \tensor \ell^2(\Z_{\geq 0})$ spanned by elements $e_{n_1} \tensor e_{n_2} \tensor e_{n_3}$ for which $n_1+n_2-n_3$ is constant, is invariant under the action of the Casimir in the tensor product $\pi^+_{k_1} \tensor \pi^P_{\rho,\eps} \tensor \pi^-_{k_3}$. 
For this reason we define, for $r \in \Z$, 
\[
k_{n,m}^r = g_n^{m+r} \tensor e_m=  e_n \tensor  l_m^{r-n}.
\]
Here $g_n^p$ and $l_n^p$ are defined by \eqref{def gn} and \eqref{def:ln}. 
Now the subspace 
\[
\mathcal A^r = \overline{ \Csp \{ k_{n,m}^r\ | \ n,m \in \Z_{\geq 0} \}}
\]
is invariant under the action of $\Om$. 

We define corresponding functions $K_{n,m}^r$ by
\begin{equation} \label{def:K2}
K_{n,m}^r(\tau,\si;k_1,\rho,\eps,k_3) = 
(-1)^{m-r}\,G_n^{m+r}(\si;k_1,\rho,\eps) G_m^{-r} (\tau;k_3,\si,-k_1-\eps), 
\end{equation}
and we denote by $m^r(\tau,\si;k_1,\rho,\eps,k_3)$ the corresponding orthogonality measure;
\[
\iint f(\tau,\si)\, dm^r(\tau,\si;k_1,\rho,\eps,k_3) = \iint f(\tau,\si)\, d\big(\mu_2^{m+r}(\si;k_1,\rho,\eps) \times \mu_2^{-r}(\tau;k_3,\si,-k_1-\eps)\big).
\]
We should be careful with the definition of the function $K_{n,m}^r$ and the corresponding measure $dm^r$ when $\si$ is in the discrete part of the support of $d\mu_2(\si;k_1,\rho,\eps)$, since in this case the second function $G_m^{-r}$ in \eqref{def:K2} must be replace by a function $F_m^{r-j}$, see also Remark \ref{rem:discr}. We leave the details to the reader. Now the following proposition is proved in the same way as Proposition \ref{prop J1}.
\begin{prop} 
For $r\in \Z$ the operator $J_1 : \mathcal A^r \rightarrow \mathcal H\big(m(\tau,\si;k_1,\rho,\eps,k_3)\big)$ defined by
\[
J_1 : k_{n,m}^r \mapsto \Big((\tau,\si) \mapsto K_{n,m}^r(\tau,\si;k_1,\rho,\eps,k_3) \Big), 
\] 
is unitary and intertwines $\pi^+_{k_1} \tensor \pi^P_{\rho,\eps} \tensor \pi^-_{k_3} (( \De \tensor 1) \circ \De(\Om))$ with $M_{\tau^2+\frac14}$.

The operator $J_2 : \mathcal A^r \rightarrow \mathcal H\big(m(\tau,\zeta;k_3,\rho,-\eps,k_3)\big)$ defined by
\[
J_2 : k_{n,m}^r \mapsto \Big( (\tau,\zeta) \mapsto K_{m,n}^{-r}(\tau,\zeta;k_3,\rho,-\eps,k_1)\Big), 
\] 
is unitary and intertwines $\pi^+_{k_1} \tensor \pi^P_{\rho,\eps} \tensor \pi^-_{k_3} (( 1 \tensor \De) \circ \De(\Om))$ with $M_{\tau^2+\frac14}$.
\end{prop}
\begin{rem} \label{rem:discr}
Using Propositions \ref{prop:+-} and \ref{prop:+P} we can decompose the tensor product $\pi^+_{k_1} \tensor \pi^P_{\rho,\eps} \tensor \pi^-_{k_3}$ into irreducible representations in two ways. This gives
\[
\begin{split}
(\pi^+_{k_1} \tensor \pi^P_{\rho,\eps}) \tensor \pi^-_{k_3} \cong&\, \left( \dirint \pi_{\si,\eps+k_1}^P\, d\si \oplus \bigoplus_j \pi^+_{\eps+k_1+j} \right) \tensor \pi^-_{k_3} \\
\cong &\, \dirint \left( \dirint \pi^P_{\tau,\eps+k_1-k_3} \,d\tau \oplus \bigoplus_l \pi^-_{k_3-k_1-\eps+l} \right) d\si \\
&\oplus \bigoplus_j \left( \dirint \pi^P_{\tau,\eps+k_1-k_3+j} \,d\tau \oplus \bigoplus_l \pi^+_{\eps+k_1-k_3+j-l} \right), \\ \\
\pi^+_{k_1} \tensor (\pi^P_{\rho,\eps}  \tensor \pi^-_{k_3}) \cong&\, \pi^+_{k_1} \tensor \left( \dirint \pi^P_{\zeta,\eps-k_3}\, d\zeta \oplus \bigoplus_j \pi^-_{k_3-\eps+j} \right) \\
\cong &\, \dirint \left( \pi^P_{\tau,\eps+k_1-k_3}\, d\tau \oplus \bigoplus_l \pi^+_{\eps+k_1-k_3+l } \right) d\zeta \\
& \oplus \bigoplus_j \left( \dirint \pi^P_{\tau, k_1+\eps-k_3-j} \, d\tau \oplus \bigoplus_l \pi^-_{k_3-k_1-\eps+j-l} \right).
\end{split}
\]
We see that the values for $\tau$ corresponding to positive discrete series representations only occur for discrete values of $\si$ and  for continuous values of $\zeta$ (i.e.~$\zeta \in [0,\infty)$). More precisely, if $\si=i(\hf-k_1-\eps+j)$, then $\tau$ can take the value $\tau= i(\hf-k_1+k_3-\eps+l+j)$, where $l \in \Z_{\geq 0}$ is such that $\tau \in i\R_{<0}$. A similar statement holds for values of $\tau$ corresponding to negative discrete series representations.

So we should be careful when writing down the function $K_{n,m}^r$ defined by \eqref{def:K2} for $\si=i(\hf-k_1-\eps+j)$, since then the second function $G_{m}^{-r}$ in \eqref{def:K2} must be replaced by a function $F_{m}^{r-j}(\tau;k_1+\eps+j,k_3)$, coming from the tensor product $\pi^+_{k_1+\eps+j} \tensor \pi^-_{k_3}$. 
\end{rem}

We first concentrate on the case $\tau \in [0,\infty)$. We define the measure $\nu^r(\si;k_1,\rho,\eps,k_3,\tau)$ to be the measure $m^r(\tau,\si;k_1,\rho,\eps,k_3)$ restricted to a line on which $\tau=$ constant. For $r \rightarrow \infty$ this measure has an infinite number of discrete mass points (cf.~Remark \ref{rem:inf discr}), and we denote the corresponding measure by $\nu$, i.e.~we omit the superscript `r'. As in subsection \ref{ssec:Racah1} we can define the Racah coefficients as the kernel in the integral operator 
$U: \mathcal L\big(\nu(\si;k_1,\rho,\eps,k_3,\tau)\big) \rightarrow \mathcal L\big(\nu(\zeta;k_3,\rho,-\eps,k_1,\tau)\big)$ mapping $K_{n,m}^r(\tau,\si;k_1,\rho,\eps,k_3)$  to $K_{m,n}^{-r}(\tau,\zeta;k_3,\rho,-\eps,k_1)$. So the Racah coefficient $U(\si,\zeta;k_1,\rho,\eps,k_3,\tau)$ is defined by
\begin{equation} \label{eq:U3}
K_{m,n}^{-r}(\tau,\zeta;k_3,\rho,-\eps,k_1) = \int U(\si,\zeta;k_1,\rho,\eps,k_3,\tau) K_{n,m}^r(\tau,\si;k_1,\rho,\eps,k_3) \, d\nu(\si;k_1,\rho,\eps,k_3,\tau).
\end{equation}
In Remark \ref{rem:discr} it is explained that discrete values of $\tau$ only occur when one of the variables $\si,\zeta$ is purely continuous and the other is purely discrete. In \eqref{eq:U3} the meaure $d\nu$ is absolutely continuous (supported on $[0,\infty)$) for $\tau=i(\hf-k_1+k_3-\eps+j) \in i\R_{<0}$, and the measure $d\nu$ is purely discrete (supported on the points $\si_l$ for $l \in \Z$ such that $\si_l=i(\hf-k_1-\eps+l) \in i\R_{<0}$) for $\tau=i(\hf+k_1-k_3+\eps+j) \in i\R_{<0}$.

We use \eqref{eq:U3} and Lemma \ref{lem1} to determine an explicit expression for the Racah coefficient as a multiple of a Wilson function. It turns out that for the discrete values of $\tau$ the Wilson function that we need can be simplified to a multiple of a Wilson polynomial \cite{Wil80}. The Wilson polynomials is a polynomial in $x^2$ which is defined by 
\[
\begin{split}
W_n(x;a,b,c,d)
&= (a+b)_n(a+c)_n (a+d)_n \F{4}{3}{-n, n+a+b+c+d-1, a+ix, a-ix}{ a+b, a+c, a+d}{1}.
\end{split}
\]
By Whipple's transformation \cite[Theorem 3.3.3]{AAR99} this polynomial is symmetric in the parameters $a,b,c,d$.
The following properties of the Wilson function turn out to be useful. 
\begin{lem} \label{lem2} \*
\begin{itemize}
\item[(i)] The Wilson function has the following symmetry property:
\[
\phi_\la(x;a+i\om,b+iy,b-iy,1-a+i\om) = \phi_\om(y;a+i\la,b+ix,b-ix,1-a+i\la).
\]
\item[(ii)] The Wilson function $\phi_\la(x;a,b,c,d)$ is symmetric in $a,b,c,1-d$.
\item[(iii)] For $\la_n=i(\tilde a +n)$, $n \in \Z_{\geq 0}$, the Wilson function is a multiple of a Wilson polynomial;
\[
\phi_{\la_n}(x;a,b,c,d) = \frac{(-1)^n W_n(x;a,b,c,d)}{\Ga(a+b+n) \Ga(a+c+n) \Ga(b+c+n) \Ga(1-d \pm ix) }
\]
\end{itemize}
\end{lem}
\begin{proof}
The first property can be read off directly from the definition \eqref{def:Wilson function} of the Wilson function after the substitution $(a,b,c,d) \mapsto (a+i\om,b+iy,b-iy,1-a+i\om)$. The second property is proved in \cite[Remark 4.5(iii)]{Gr03}. The third property follows directly from \eqref{eq:Wil4F3}, where for $\la=\la_n$ the second $_4F_3$-function vanishes because of the factor $\Ga(\tilde a \pm i \la)^{-1}$. 
\end{proof}
\begin{thm} \label{thm:U+P-}
The Racah coefficient $U(\si,\zeta;k_1,\rho,\eps,k_3,\tau)$ can be expressed as a Wilson polynomial or a Wilson function as follows:
\begin{itemize}
\item for $\tau \in [0,\infty)$,
\[
K_1\,\phi_\zeta(\si;k_1+i\rho,k_3+i\tau, k_3-i\tau, 1-k_1+i\rho),
\]
\item for $\tau=i(\hf-k_1+k_3-\eps+j) \in i\R_{<0}$ and $\zeta=i(\hf-k_3+\eps+l)\in i\R_{<0}$,
\[
K_2 \, W_{j-l}(\si;2k_3-k_1-\hf-\eps-j,k_1+i\rho,k_1-i\rho,\hf-k_1-\eps-j), 
\]
\item for $\tau=i(\hf+k_1-k_3+\eps+j) \in i\R_{<0}$ and $\si=i(\hf-k_1-\eps+l)\in i\R_{<0}$,
\[
K_3\, W_{j-l}(\zeta;2k_1-k_3-\hf+\eps-j,k_3+i\rho,k_3-i\rho,\hf-k_3+\eps-j), 
\]
\end{itemize}
where
\[
\begin{split}
K_1= & 2 \sin\pi(\hf-k_1-\eps+k_3 \pm i\tau) \, \Ga(2k_1) \Ga(2k_3)  \Ga(\pm 2i\tau) \Ga(\hf-\eps+k_3 \pm i\zeta) \Ga(\hf+k_1+\eps \pm i\si),\\
K_2= &\frac{ \Ga(2k_3-2k_1-2\eps-1) \Ga(2k_3) \Ga(2k_3-2\eps-l)  }{ \Ga(2k_1-2k_3+2\eps) \Ga(2k_3-\hf-\eps-l \pm i\rho) } \frac{(-1)^{j-l} l!(\hf-k_3+k_1+\eps)_j }{  (\frac32-k_3+k_1+\eps)_j (\hf+k_1+\eps \pm i\si)_j (2k_1)_{j-l} },\\
K_3= &\frac{ \Ga(2k_1-2k_3+2\eps-1) \Ga(2k_1) \Ga(2k_1+2\eps-l)  }{ \Ga(2k_3-2k_1-2\eps) \Ga(2k_1-\hf+\eps-l \pm i\rho) } \frac{(-1)^{j-l} l!(\hf-k_1+k_3-\eps)_j }{ (\frac32-k_1+k_3-\eps)_j (\hf+k_3-\eps \pm i\si)_j (2k_3)_{j-l} }
\end{split}
\]
\end{thm}
\begin{proof}
First assume that $\tau \in [0,\infty)$. From \eqref{eq:U3} with $m=r=0$ and using the orthogonality of the continuous dual Hahn polynomial $s_n(\si;k_1+i\rho,k_1-i\rho,\hf-\eps-k_1)$, we find
\[
\begin{split}
U_\zeta(\si;k_1,\rho,\eps,k_3,\tau) &= \\
 C \ \sum_{n=0}^\infty &\frac{ (\hf-\eps \pm i\rho)_n }{n!\, (2k_1)_n} \F{3}{2}{ -n, \hf-\eps-k_1 +i\si, \hf-\eps-k_1 -i\si}{ \hf-\eps +i\rho, \hf-\eps -i\rho}{1} \\
&\times \F{3}{2}{ -n, \hf-\eps-k_1+k_3+i\tau, \hf-\eps-k_1+k_3-i\tau}{ \hf-\eps+k_3+i\zeta, \hf-\eps+k_3-i\zeta}{1} .
\end{split}
\]
Using Lemma \ref{lem1} with parameters
\[
(a,b,c,d,f,x,\la)\mapsto (k_1+i\si,k_3+i\zeta,k_3-i\zeta,1-k_1+i\si,\hf-\eps-k_1+k_3,\rho,\tau)
\]
we find that the Racah coefficient is a multiple of the Wilson function $\phi_\tau(\rho;k_1+i\si,k_3+i\zeta,k_3-i\zeta,1-k_1+i\si)$.
Then the theorem follows from applying Lemma \ref{lem2}(i), and writing out explicitly the factor $C$. 

For discrete values of $\tau$ the proof runs along the same lines, using also Lemma \ref{lem2}(ii) and (iii).
\end{proof}
\begin{rem}
Note that the Wilson function and the Wilson polynomials in the theorem do not depend on $\eps$.
\end{rem}

\subsection{The Wilson function transform of type II} \label{ssec:WilsonII}
The Wilson function transform of type II \cite{Gr03} is an integral transform depending on the four Wilson parameters $a,b,c,d$ and an extra parameter $t$. Let the parameters $a,b,c,d,t$ satisfy 
\[
\Re(a),\Re(b)>0,\quad \overline{a} = 1-d, \quad \overline b=c, \quad t\in \R.
\]
For these parameters we define the Hilbert space $\mathcal H=\mathcal H(a,b,c,d;t)$ to be the space consisting of even functions that have finite norm with respect to the inner product 
\[
\inprod{f}{g}_{\mathcal H} = \int f(x)\overline{g(x)}\, dh(x;a,b,c,d;t),
\]
where $h$ is the measure defined by
\[
\int f(x)\, dh(x;a,b,c,d;t) = \frac{C}{2\pi} \int_0^\infty f(x)\, H(x)\, dx + iC \sum_{x \in \mathcal D} f(x)\, \Res{z=x} H(z),
\]  
where
\[
\begin{split}
H(x)&=H(x;a,b,c,d;t) = \frac{ \Ga(a \pm ix) \Ga(b \pm ix) \Ga(c \pm ix) \Ga(1-d \pm ix) }{ \sin \pi(t \pm ix) \Ga( \pm 2ix) },\\
C&=C(a,b,c,d;t) = \sqrt{ \sin \pi(a +t) \sin \pi(b+t) \sin \pi(c+t) \sin \pi(1-d+t) },
\end{split}
\]
and $\mathcal D$ is the infinite discrete set
\[
\mathcal D = \{ i(t-n)\ | \ n \in \Z, t-n<0 \}.
\]
We define the dual parameter $\tilde t$ by
\[
\tilde t = 1- \tilde b-c -t,
\]
and the dual parameters $\tilde a, \tilde b, \tilde c, \tilde d$ are still defined by \eqref{eq:dual param}. 
Now we define an operator $\mathcal G=\mathcal G(a,b,c,d;t)$ by
\[
(\mathcal Gf)(\la) = \inprod{f}{\phi_\la(\cdot;a,b,c,d)}_{\mathcal H},
\]
for all functions $f \in \mathcal H$ such that the integral on the right hand side converges. The operator $\mathcal G$ is called the Wilson transform of type II, or just the Wilson transform II. The following theorem is Theorem 5.10 in \cite{Gr03}
\begin{thm}
The operator $\mathcal G : \mathcal H(a,b,c,d;t) \rightarrow \mathcal H(\tilde a, \tilde b, \tilde c, \tilde d; \tilde t)$ extends to a unitary operator, and its inverse is given by $\mathcal G(\tilde a, \tilde b, \tilde c, \tilde d;\tilde t)$.
\end{thm}

We will sketch how the unitarity of $\mathcal G$ can be obtained from the interpretation of the Wilson functions as Racah coefficients for the tensor product $\pi^+_{k_1} \tensor \pi^P_{\rho,\eps} \tensor \pi^-_{k_3}$.

Let $\tau \in [0,\infty)$. Let us define an operator $\mathcal I=\mathcal I(k_1,\rho,\eps,k_3,\tau):\mathcal L\big(\nu(\si;k_1,\rho,\eps,k_3,\tau)\big) \rightarrow \mathcal L\big(\nu(\zeta;k_3,\rho,-\eps,k_1,\tau)\big)$ by
\[
(\mathcal I f)(\zeta) = \int f(\si) U(\si,\zeta;k_1,\rho,\eps,k_3)\, d\nu(\si;k_1,\rho,\eps,k_3,\tau).
\]
Then, in the same way as in subsection \ref{ssec:WilsonI}, we can obtain from \eqref{eq:U3} that $\mathcal I$ is unitary, and its inverse is the same tranform after $k_1 \leftrightarrow k_3$ and $\eps \mapsto -\eps$. Comparing the measure $\nu$ with the measure $h$ for the Wilson transform II, and using Theorem \ref{thm:U+P-} we see that we have, for an appropriate function $f$,
\[
\begin{split}
(\mathcal I f)(\zeta)= \frac{E}{B(\zeta)} \int \phi_\zeta(\si;a,b,c,d) A(\si) f(\si)\, dh(\si;a,b,c,d;t),
\end{split}
\]
where
\[
(a,b,c,d,t) = (k_1+i\rho,k_3+i\tau, k_3-i\tau, 1-k_1+i\rho, \hf-k_1-\eps),
\]
and
\[
\begin{split}
A(\si) = &\Ga(\hf+k_1+\eps \pm i\si)^{-1}, \qquad 
B(\zeta) =  \Ga(\hf -\eps+ k_3 \pm i\zeta)^{-1}, \\
E =&\sqrt{ \frac{\Ga( \hf - k_3 +\eps+k_1 \pm i\tau) \Ga(\hf +\eps \pm i\rho) }{ \Ga(\hf-k_1-\eps+k_3\pm i\tau)\Ga(\hf-\eps\pm i\rho)} }.
\end{split}
\]
So we have $(\mathcal I f)(\zeta) = E\,\big( M_{B(\zeta)}^{-1} \circ \mathcal G \circ M_{A(\cdot)} f \big) (\zeta)$, where $\mathcal G$ denotes the Wilson transform II with the above parameters. Note that $\tilde A(\cdot) = B(\cdot)$ and $\tilde E = E^{-1}$. Here we use the notation $\tilde O=O(\tilde a,\tilde b,\tilde c,\tilde d, \tilde t)$ for any object $O=O(a,b,c,d,t)$ depending on the parameters $a,b,c,d,t$. Observe here that replacing $(a,b,c,d,t)$ by their dual parameters is the same as replacing $(k_1,k_3,\eps)$ by $(k_3, k_1,-\eps)$. The operator $M_{A(\cdot)}$ is invertible, because $A(\si) \neq 0$ for $\si$ in the support of the measure $\nu(\si;k_1,\rho,\eps,k_3,\tau)$.

The measure $\nu = \nu(\cdot;k_1,\rho,\eps,k_3,\tau)$ is related to the measure $h(\cdot;a,b,c,d;t)$ by
\[
\int f(\si)\, d\nu(\si) = N \, \int f(\si)\,|A(\si)|^2\, dh(\si;a,b,c,d;t), 
\]
where, if $C$ is the normalizing constant in the definition of the measure $h$ satisfying $\tilde C = C$, 
\[
N =  \frac{ \Ga(\hf + k_1+\eps -k_3 \pm i\tau) }{2C\, \Ga( \pm 2i\tau) \Ga(2k_1) \Ga(2k_3) \Ga(\hf -\eps \pm i\rho) }.
\]
So we have $\inprod{f_1}{f_2}_{\mathcal L(\nu)} = N \inprod{A\,f_1}{A\,f_2}_{\mathcal H}$. On the other hand, we find from the unitarity of $\mathcal I$, 
\[
\begin{split}
\inprod{f_1}{f_2}_{\mathcal L(\nu)} &= \inprod{\mathcal I f_1}{\mathcal I f_2}_{\mathcal L(\tilde \nu)} \\
&= |E|^2 \inprod{ \tilde A^{-1} \,\mathcal G (A\,f_1) } {\tilde A^{-1} \,\mathcal G (A\,f_2) }_{\mathcal L(\tilde \nu)} \\
&= |E|^2 \tilde N \inprod{  \mathcal G (A\,f_1) }{ \mathcal G (A\,f_2) }_{\tilde{\mathcal H}}.
\end{split}
\]
Now from  $|E|^2 \tilde N = N$ it follows that the Wilson function transform $\mathcal G$ is unitary.

Let us also remark that the infinite set of discrete mass points of the measure $h$ for the Wilson transform II corresponds to the infinite sum of discrete series in the tensor product decomposition of $\pi^\pm \tensor \pi^P$. The absolutely continuous part of $h$ corresponds to the direct integral over the principal unitary series.\\

For $\tau \not\in [0,\infty)$ the Wilson functions in Theorem \ref{thm:U+P-} reduce to Wilson polynomials. The Wilson polynomials are orthogonal with respect to a positive measure \cite{Wil80}. Let the parameters $a,b,c,d$ be such that that non-real parameters appear in pairs of complex conjugates with positive real part, and such that the pairwise sum of the real parameters is positive. Then the orthogonality relations for the Wilson polynomials read
\[
\frac{1}{2\pi}\int_0^\infty W_n(x) W_m(x) w(x) dx + \sum_j W_n(x_j)W_m(x_j) w_j= h_n\, \de_{nm},
\]
where we denote $W_n(x) = W_n(x;a,b,c,d)$. The sum is over all $j \in \Z_{\geq 0}$ such that $x_j=i(e+j) \in i\R_{<0}$, where $e$ is the smallest of the real parameters. In particular, the sum is empty when there is no negative parameter.
The weights and the squared norm are given by
\[ 
\begin{split}
w(x) = &\frac{ \Ga(a \pm ix) \Ga(b \pm ix) \Ga(c \pm ix) \Ga(d \pm ix) } { \Ga( \pm 2ix) },\\
w_j = & \Res{x=x_j} w(x),\\
h_n =& \frac{ \Ga(a+b+n) \Ga(a+c+n) \Ga(a+d+n) \Ga(b+c+n) \Ga(b+d+n) \Ga(c+d+n) }{ \Ga(a+b+c+d+2n) }\\
&\times  (n+a+b+c+d-1)_n \, n!.
\end{split}
\]
The Wilson polynomials form an orthogonal basis for the Hilbert space consisting of even functions that have finite norm with respect to the inner product
\[
\inprod{f}{g} = \frac{ 1}{2\pi} \int_0^\infty f(x)\overline{g(x)} w(x) dx + \sum_j f(x_j) \overline{g(x_j)} w_j.
\]
So we have the following pair of unitary integral transforms
\[
\begin{cases}
\displaystyle (\mathcal Wf)(n) = \inprod{f}{W_n},\\
\displaystyle f(x)= \sum_{n=0}^\infty (\mathcal W f)(n) W_n(x) h^{-1}_n
\end{cases}
\]
We call $\mathcal W=\mathcal W(a,b,c,d)$ the polynomial Wilson transform.

It can be checked that the unitarity of polynomial Wilson transform in the case $a,d>0$ and $\overline b=c$ can be obtained from the interpretation of the Wilson polynomials as Racah coefficients in Theorem \ref{thm:U+P-}. We leave this to the interested reader. Another interesting fact that we obtain from this interpretation is that the polynomial Wilson transform $\mathcal W(a,b,c,1-d)$ can be considered as a degenerate case of the Wilson function transform of type II $\mathcal G(a,b,c,d;t)$, namely the case that the parameters satisfy
\[
a=1-\overline d, \quad \Re(a)>0, \quad b,c\in\R, \quad b+c>0, \quad t=n-c+1,
\]
for $n \in \Z_{\geq 0}$.


\section{More Racah coefficients} \label{sec:more}
In this section we determine some more Racah coefficients for $3$-fold tensor product representations of $\su(1,1)$. \\

Let us start with the tensor product representation $\pi^+_{k_1} \tensor \pi^+_{k_2} \tensor \pi^-_{k_3}$. We assume for simplicity that there are no discrete terms in the decomposition of $\pi^+_{k_2} \tensor \pi^-_{k_3}$. The $\Om$-invariant subspace in this case is
\[
\mathcal A^r = \overline{\Csp \{ e_{n_1} \tensor e_{n_2} \tensor e_{n_3} \ | \ n_1,n_2,n_3 \in \Z_{\geq 0}, \ n_1+n_2-n_3= r \}} \subset \ell^2(\Z_{\geq 0})^{\tensor 3}.
\] 
Using the results from section \ref{sec:CGC} we obtain operators $J_1$ and $J_2$ as in Propositions \ref{prop J1} and \ref{prop J2}.
Let us define the measure $m_1$  by
\[
\begin{split}
\iint f(\tau,j)\, dm_1(\tau,j;k_1,k_2,k_3) = \sum_{j=0}^\infty \int_\R f(\tau,j)\, W(j;k_1,k_2)\, d\mu_1(\tau;k_1+k_2+j,k_3),
\end{split}
\]
and let $dm_2$ be the product measure defined by
\[
\iint f(\tau,\rho) \,dm_2(\tau,\rho;k_1,k_2,k_3) = \iint f(\tau,\rho)\, d\big(\mu_1(z;k_2,k_3) \times \mu_2(x;k_1,\rho,k_2-k_3)\big).
\]
Now the operators $J_1$ and $J_2$ are given by
\[
\begin{split}
J_1 : \mathcal A^r &\rightarrow \mathcal H\big(m_1(\tau,j;k_1,k_2,k_3)\big)\\
e_{n_1} \tensor e_{n_2} \tensor e_{n_3} & \mapsto \Big((\tau,j) \mapsto v_{n_1}^{n_3+r}(j;k_1,k_2) F_{n_3}^{r-j}(\tau;k_1+k_2+j,k_3) \Big),
\end{split}
\]
and
\[
\begin{split}
J_2 : \mathcal A^r & \rightarrow \mathcal H\big(m_2(\tau,\rho;k_1,k_2,k_3)\big) \\
e_{n_1} \tensor e_{n_2} \tensor e_{n_3} & \mapsto \Big((\tau,\rho) \mapsto  F_{n_3}^{r-n_1}(\rho;k_2,k_3) G_{n_1}^r(\tau;k_1,\rho,k_2-k_3) \Big).
\end{split}
\]
These operators intertwine $\pi^+_{k_1} \tensor \pi^+_{k_2} \tensor \pi^-_{k_3}((\De \tensor 1)\circ \De(\Om))$, respectively $\pi^+_{k_1} \tensor \pi^+_{k_2} \tensor \pi^-_{k_3}((1 \tensor \De)\circ \De(\Om))$, with $M_{\tau^2+\frac14}$. We restrict the measure $m_1(x,y;k_1,k_2,k_3)$ to a line on which $\tau=$ constant, then we obtain a measure which is supported on $\Z_{\geq 0}$. We will denote a mass point of this measure at $j \in \Z_{\geq 0}$ by $w(j;k_1,k_2,k_3,\tau)$. Furthermore, we define 
the measure $\nu(\rho;k_1,k_2,k_3,\tau)$ to be the measure $m_2(\tau,\rho;k_1,k_2,k_3)$ restricted to a line on which $\tau=$ constant. Now the Racah coefficients $U(j,\rho;k_1,k_2,k_3,\tau)$ can be determined from
\begin{equation} \label{eq:U4'}
\begin{split}
F_{n_3}^{r-n_1}&(\rho;k_2,k_3) G_{n_1}^r(\tau;k_1,\rho,k_2-k_3) = \\ &\sum_{j=0}^\infty U(j,\rho;k_1,k_2,k_3,\tau) v_{n_1}^{n_3+r}(j;k_1,k_2) 
F_{n_3}^{r-j}(\tau;k_1+k_2+j,k_3) \, w(j;k_1,k_2,k_3,\tau),
\end{split}
\end{equation}
or equivalently from 
\begin{equation} \label{eq:U4}
\begin{split}
v_{n_1}^{n_3+r}&(j;k_1,k_2) F_{n_3}^{r-j}(\tau;k_1+k_2+j,k_3)  = \\
&\int_0^\infty U(j,\rho;k_1,k_2,k_3,\tau) F_{n_3}^{r-n_1}(\rho;k_2,k_3) G_{n_1}^r(\tau;k_1,\rho,k_2-k_3)\, d\nu(\rho;k_1,k_2,k_3,\tau).
\end{split}
\end{equation}
Then from Lemmas \ref{lem1} and \ref{lem2} we find the following expression for the Racah coefficient.
\begin{thm} \label{thm:U++-}
The Racah coefficient $U(j,\rho;k_1,k_2,k_3,\tau)$ is the Wilson polynomial
\[
K\ W_j(\rho;k_2-k_3+\hf, k_2+k_3-\hf, k_1+i\tau, k_1-i\tau),
\]
where 
\[
K = \frac{2\,(-1)^j \Ga(\pm 2i\tau) \Ga(2k_1+2k_2+2j) \Ga(2k_3) \sin \pi (k_1+k_2-k_3+\hf \pm i\tau)}{(2k_1)_j  \Ga(k_1+k_2+k_3-\hf +j \pm i\tau)    \sqrt{(k_1+k_2-k_3+\hf \pm i\tau)_j (2k_1+2k_2+j)_j}}.
\]
\end{thm}
\begin{proof}
We use formula \eqref{eq:U4}. We use the orthogonality relations for $F_{n_3}^{r-n_1}(\rho)$ to expand the Racah coefficient in terms of the Clebsch-Gordan coefficients. Then by putting $n_1=0$ and $r=j$ and writing out the Clebsch-Gordan coefficients explicitly, we obtain
\[
\begin{split}
U(j,\rho;k_1,&k_2,k_3,\tau) = \\
C' \sum_{n_3=0}^\infty&  \frac{ (j+1)_{n_3} (2k_2+j)_{n_3} }{ n_3!\, (2k_3)_{n_3} } \F{3}{2}{ -n_3, k_2-k_3+j+\hf+i\rho, k_2-k_3+j+\hf-i\rho }{ 2k_2+j, j+1}{1} \\
\times &\F{3}{2}{ -n_3 ,k_1+k_2-k_3+j+\hf+i\tau, k_1+k_2-k_3+j+\hf -i \tau}{ 1, 2k_1+2k_2+2j}{1}.
\end{split}
\]
Applying Lemma \ref{lem1} and Lemma \ref{lem2}(i) shows that the Racah coefficient is equal to the Wilson function $C \, \phi_{i(k_1+k_2-\hf+j)}(\rho;k_2+k_3-\hf,  k_1+i\tau, k_1-i\tau, k_3-k_2+\hf)$, where the factor $C$ is given by
\[
\begin{split}
2& \Ga(\pm 2i\tau) \Ga(2k_1) \Ga(2k_3) \Ga(k_3-k_2+\hf \pm i\rho) \\
& \times \frac{ \Ga(2k_1+2k_2+2j) \sin \pi (k_1+k_2-k_3+\hf \pm i\tau)}{    \sqrt{(k_1+k_2-k_3+\hf\pm i\tau)_j (2k_1+2k_2+j)_j} }.
\end{split}
\]
Finally the theorem follows from applying Lemma \ref{lem2}(ii).
\end{proof}
\begin{rem}
Koornwinder \cite{Koo86} found a similar group theoretic interpretation of the Wilson polynomials by establishing the Wilson polynomials as connection coefficients between two $O(p) \times O(q) \times O(r)$-invariant bases on a hyperboloid in $\R^{p+q+r}$.
\end{rem}

Finally we give the Racah coefficients for the tensor product representations $\pi^+_{k_1} \tensor \pi^P_{\rho,\eps} \tensor \pi^+_{k_3}$ and $\pi^+_{k_1} \tensor \pi^+_{k_2} \tensor \pi^P_{\rho,\eps}$. The derivation of these Racah coefficients is completely the same as before, therefore we omit all details.

The Racah coefficient $U(\si,\zeta;k_1,\rho,\eps,k_3)$ for the tensor product $\pi^+_{k_1} \tensor \pi^P_{\rho,\eps} \tensor \pi^+_{k_3}$ is the kernel in the transform $\mathcal I: \mathcal L\big(\nu(\si;k_1,\rho,\eps,k_3,\tau)\big) \rightarrow \mathcal L\big(\nu(\zeta;k_3,\rho,\eps,k_1,\tau)\big)$, which is defined by
\[
\begin{split}
\Big(\mathcal I\big[\si \mapsto (G_{n_1}^{r-n_3}(\si; k_1,\rho,\eps)& G_{n_3}^r(\tau;k_3,\si, \eps+k_1)\big]\Big)(\zeta)=\\
&  G_{n_3}^{r-n_1}(\zeta; k_3,\rho,\eps) G_{n_1}^r(\tau;k_1,\zeta, \eps+k_3).
\end{split}
\]
Here the measure $\nu(\si;k_1,\rho,\eps,k_3,\tau)$ is the restriction of the product measure
\[
\mu_2(\si;k_1,\rho,\eps,k_3) \times \mu_2(\tau;k_3,\si, \eps+k_1) 
\]
to a line on which $\tau=$ constant. The Racah coefficients can be expressed as a multiple of a Wilson function
\[
\begin{split}
U(\si,\zeta;&k_1,\rho,\eps,k_3) =\\
&\frac{2\pi \Ga(\pm 2i\tau) \Ga(2k_1) \Ga(2k_3) \Ga(\hf-\eps \pm i\rho)}{ \Ga(\hf-k_1-k_3-\eps \pm i\tau)}  \phi_\zeta(\si;k_3+i\tau, k_1+i\rho, k_1-i\rho, 1-k_3+i\tau)
\end{split}
\]
The $\Ga$-functions in the measure $\nu$ depending on $\eps$ are cancelled by the same $\Ga$-functions in the expression of the Racah coefficient, so the corresponding integral transform $\mathcal I$ does not depend on $\eps$. This integral transform $\mathcal I$ is equivalent to the Wilson transform I.\\

The Racah coefficient $U(j,\si;k_1,k_2,\rho,\eps)$ for the tensor product $\pi^+_{k_1} \tensor \pi^+_{k_2} \tensor \pi^P_{\rho,\eps}$ is the kernel in the transform $\mathcal I : \ell^2(\Z_{\geq 0}, w(j;k_1,k_2,\rho,\eps)) \rightarrow \mathcal L\big(\nu(\si;k_1,k_2,\rho,\eps,\tau))$, which is defined by
\[
\begin{split}
\Big(\mathcal I\big[j \mapsto v_{n_1}^{r-n_3}(j;k_1,k_2) &G_{r-n_3-j}^{r-j}(\tau;k_1+k_2+j,\rho,\eps)\big]\Big)(\si)= \\
&G_{r-n_1-n_3}^{r-n_1}(\si;k_2,\rho,\eps) G_{n_1}^r(\tau;k_1,\si, k_2+\eps).
\end{split}
\]
Here $w(j;k_1,k_2,\rho,\eps,\tau)$ is the mass at the points $j \in \Z_{\geq 0}$ of the measure 
\[
W(j;k_1,k_2) \times \mu_2(\tau;k_1+k_2+j,\rho,\eps)
\]
restricted to a line on which $\tau=$ constant, and $\nu(\si;k_1,k_2,\rho,\eps,\tau)$ is the restriction of the product measure
\[
\mu_2(\si;k_2,\rho,\eps) \times \mu_2(\tau;k_1,\si, k_2+\eps).
\]
The Racah coefficient can be expressed as a multiple of a Wilson polynomial,
\[
\begin{split}
U(j,\si&;k_1,k_2,\rho,\eps) = \\
&\frac{(-1)^j\, 2\pi\ \Ga(\pm 2i\tau)  \Ga(\hf-\eps \pm i\rho) \Ga(2k_1+2k_2+2j)}{(2k_1)_j \Ga(\hf-k_1-k_2-\eps \pm i\tau) \Ga(k_1+k_2+j+i\rho \pm i\tau) \Ga(k_1+k_2+j-i\rho \pm i\tau) }\\
\times&  \sqrt{ \frac{  (k_1+k_2+\eps+\hf \pm i\tau)_j }{(2k_1+2k_2+j)_j }} \ W_j(\si;k_2+i\rho,k_2-i\rho,k_1+i\tau, k_1-i\tau).
\end{split}
\]
The corresponding integral transform $\mathcal I$ does not depend on $\eps$, and $\mathcal I$ is equivalent to the polynomial Wilson transform.


\section{Integral and summation formulas} \label{sec:identities}
\subsection{Formulas related to 3-fold tensor products}
In the previous sections we found explicit formulas for the Racah coefficients by expanding the Racah coefficients in terms of the Clebsch-Gordan coefficients. Writing the Clebsch-Gordan coefficients as (continuous) dual Hahn polynomials, we obtain several summation formulas for these polynomials.

Let us start with the tensor product $\pi^+_{k_1} \tensor \pi^+_{k_2} \tensor \pi^-_{k_3}$ from section \ref{sec:more}. We write the Clebsch-Gordan coefficients as (continuous) dual Hahn polynomials and the Racah coefficient as a Wilson polynomial, then formula \eqref{eq:U4'} gives a summation formula involving Wilson, dual Hahn and continuous dual Hahn polynomials. Note that the sum is finite, since the dual Hahn polynomial $T_n(j;\ga, \de,N)$ equals zero for $j > N$. Using the orthogonality of the dual Hahn polynomials, we obtain a summation formula which can be considered as a generalized convolution identity for continuous dual Hahn polynomials.
\begin{thm}  \label{thm:++-}
The continuous dual Hahn polynomials satisfy the following generalized convolution identity: 
\[
\begin{split}
\sum_{l=0}^{n+j} & \binom{n+j}{l} (2k_1)_j \ T_l(j;2k_1-1, 2k_2-1, n+j) \\
& \times S_{n+j-l}(\rho;k_2-k_3+\hf,k_2+k_3-\hf, k_3-k_2+l+r+\hf) \\
&\times S_l(\tau;k_1+i\rho, k_1-i\rho, k_3-k_2-k_1+r+\hf)\,  \\
=&\,W_j(\rho;k_2-k_3+\hf, k_2+k_3-\hf, k_1+i\tau, k_1-i\tau) \\
& \times S_n(\tau; k_1+k_2-k_3+j+\hf, k_1+k_2+k_3+j-\hf, k_3-k_2-k_1+r+\hf).
\end{split}
\]
\end{thm}
\begin{rem}
Several formulas with a similar structure, involving other hypergeometric orthogonal polynomials, are found by Koelink and Van der Jeugt \cite{KVdJ98}. These formulas and Theorem \ref{thm:++-} can be obtained as limit cases of one ``master formula", a generalized convolution identity for Wilson polynomials, which is obtained by Lievens and Van der Jeugt \cite[Thm.1]{LVdJ02}. In Theorem \ref{thm:BE+++-} we show that the ``master formula" itself also has an interpretation in the representation theory of $\su(1,1)$.
\end{rem}

Next we consider the tensor product $\pi^+_{k_1} \tensor \pi^-_{k_2} \tensor \pi^+_{k_3}$ from section \ref{sec:+-+}. Writing the continuous dual Hahn polynomials in \eqref{eq:U2} in their ususal normalisation \eqref{def:cont dHahn}, and using Theorem \ref{thm:U+-+} and orthogonality relations for the continuous dual Hahn polynomials, we obtain an expansion formula for the Wilson funtions.
\begin{thm} \label{thm:+-+}
The continuous dual Hahn polynomials and the Wilson functions satisfy
\[
\begin{split}
\phi_\si(&\rho; k_3 + i\tau, k_1+k_2-\hf, k_1-k_2+\hf, 1-k_3+i\tau) \\
\times S_m&(\tau;k_3+i\rho, k_3-i\rho, k_2-k_1-k_3+r+\hf) = \\
\sum_{n=0}^\infty& C_n\, S_{n}(\rho;k_1-k_2+\hf,k_1+k_2-\hf,k_2-k_1+m+r+\hf)\\ &\times S_{n}(\tau;k_1+i\si, k_1-i\si, k_2-k_1-k_3+r+\hf) \\
&\times S_m(\si;k_3-k_2+\hf, k_2+k_3-\hf, k_2-k_3+n+r+\hf), 
\end{split}
\]
where
\[
\begin{split}
C_n =&
 \frac{ \Ga(k_2-k_1+m+r+\hf \pm i\rho)} {n!\,  \Ga(2k_1+n) \Ga(2k_2+n+r+n) \Ga(n+r+m+1) \Ga(k_3-i\tau\pm i\rho) \Ga(k_3+i\tau \pm i\rho) }.
\end{split}
\]
\end{thm}

The tensor product $\pi^+_{k_1} \tensor \pi^P_{\rho,\eps} \tensor \pi^-_{k_3}$ leads to a summation formula which is equivalent to the formula in Theorem \ref{thm:+-+}. The tensor product $\pi^+_{k_1} \tensor \pi^+_{k_2} \tensor \pi^+_{k_3}$ leads to a generalized convolution identity for (dual) Hahn polynomials, see \cite{KVdJ98}.\\

Since the coproduct $\De$ is cocommutative we have $\pi_1 \tensor \pi_2 \cong \pi_2 \tensor \pi_1$ for representations $\pi_i$, $i=1,2$. We use this symmetry of the tensor product to find relations between Racah coefficients. In order to make the symmetries in the tensor products more transparant, we will denote $\pi_1 \tensor \pi_2 \tensor \pi_3 \circ (\De \tensor 1)\circ \De$ by $[\pi_1 \tensor \pi_2] \tensor \pi_3$, and similarly for the operator related to $(1 \tensor \De)\circ \De$. In order to distinguish between Racah coefficients for different tensor products we add superscripts $+$, $-$ or $P$ to the representation labels. For instance, $U(\si,\zeta;k_1^+,(\rho,\eps)^P,k_3^-, \tau)$ is the Racah coefficient for the tensor product representation $\pi^+_{k_1} \tensor \pi_{\rho,\eps}^P \tensor \pi_{k_3}^-$. We will use similar notation for the corresponding measures. Morover, we denote the corresponding integral transforms by $\mathcal I([k_1^+,(\rho,\eps)^P],k_3^-, \tau)$ and $\mathcal I(k_1^+,[(\rho,\eps)^P,k_3^-], \tau)$, which are each others inverse.

From symmetries in the tensor product $\pi^+_{k_1} \tensor \pi^+_{k_2} \tensor \pi^-_{k_3}$ we find the following theorem. 
\begin{thm} \label{thm:symI}
For $k_1+k_3-\hf>0$ and $k_1-k_3+\hf>0$ we have
\[
\begin{split}
\frac{ 1}{ 2\pi} &\int_0^\infty W_j(\si;k_1-k_3+\hf, k_1+k_3-\hf, k_2+i\tau, k_2-i\tau) \\
& \times \phi_\rho(\si;k_2+i\tau,k_1+k_3-\hf, k_1-k_3+\hf, 1-k_2+i\tau) \\
& \times \frac{ \Ga(k_2 + i\tau \pm i\si) \Ga(k_2-i\tau \pm i\si) \Ga(k_1+k_3- \hf \pm i\si) \Ga(k_1-k_3+\hf \pm i\si)}{ \Ga( \pm 2i\si) }\, d\si \\
=\,& (-1)^j W_j(\rho;\,k_2-k_3+\hf, k_2+k_3-\hf, k_1 + i\tau, k_1-i\tau). 
\end{split}
\]
\end{thm}
\begin{proof}
We can factorize the transform $\mathcal I(k_1^+,[k_2^+,k_3^-],\tau)$ using the following chain of equivalent tensor products:
\[
\begin{split}
\pi^+_{k_1} \tensor [\pi^+_{k_2} \tensor \pi^-_{k_3}] &\rightarrow 
\pi^+_{k_1} \tensor [\pi^-_{k_3} \tensor \pi^+_{k_2}] \xrightarrow{\mathcal I}
[\pi^+_{k_1} \tensor \pi^-_{k_3}] \tensor \pi^+_{k_2} \\
& \rightarrow \pi^+_{k_2} \tensor [\pi^+_{k_1} \tensor \pi^-_{k_3}] \xrightarrow{\mathcal I} [\pi^+_{k_2} \tensor \pi^+_{k_1}] \tensor \pi^-_{k_3} \rightarrow  [\pi^+_{k_1} \tensor \pi^+_{k_2}] \tensor \pi^-_{k_3}.
\end{split}
\]
The $\mathcal I$ above an arrow indicates that the arrow corresponds to an integral transform with a Racah coefficient as a kernel. 
This gives us the following factorization,
\begin{equation} \label{eq:factorI}
\mathcal I(k_1^+, [k_2^+, k_3^-],\tau) =  C\, \mathcal I(k_2^+,[k_1^+,k_3^-],\tau) \circ \mathcal I(k_1^+,[k_3^-, k_2^+],\tau),
\end{equation}
for a (still to be determined) non-zero factor $C$. This factor can be found from applying both sides to the same function. We denote $\mathcal J= \mathcal I(k_2^+,[k_1^+,k_3^-],\tau) \circ \mathcal I(k_1^+,[k_3^-, k_2^+],\tau)$, then from applying \eqref{eq:U2} and \eqref{eq:U4} we obtain 
\[
\begin{split}
\Big(\mathcal J \big[\rho \mapsto F_{n_3}^{r-n_1}(\rho);k_2,k_3) G_{n_1}^r(\tau;k_1,\rho,k_2-k_3)\big] \Big)(j)= v_{n_2}^{r+n_3}(j;k_2,k_1) F_{n_3}^{r-j}(\tau;k_1+k_2+j,k_3),
\end{split}
\]
where $r=n_1+n_2-n_3$. Also, from \eqref{eq:U4} we obtain for $\mathcal I = \mathcal I(k_1^+,[k_2^+, k_3^-],\tau)$
\[
\Big(\mathcal I \big[\rho \mapsto F_{n_3}^{r-n_1}(\rho);k_2,k_3) G_{n_1}^r(\tau;k_1,\rho,k_2-k_3)\big]\Big)(j) = v_{n_1}^{r+n_3}(j;k_1,k_2) F_{n_3}^{r-j}(\tau;k_1+k_2+j,k_3) .
\]
Now $C$ is obtained from
\[
v_{n_1}^{n_1+n_2}(j;k_1,k_2) = (-1)^j \frac{ (2k_2)_j }{(2k_1)_j} v_{n_2}^{n_1+n_2}(j;k_2,k_1),
\]
which follows directly from \eqref{eq:symT}.

We apply both sides of \eqref{eq:factorI} to a suitable function $f$. If we choose $f$ smooth enough we may interchange the order of integration on the right hand side, and since $f$ was choosen arbitrary we then obtain
\[
\begin{split}
(-1)^j\frac{(2k_1)_j}{(2k_2)_j} &U(j,\rho;k_1^+, k_2^+, k_3^-, \tau) = \\ &\int_0^\infty U(j,\si;k_2^+, k_1^+, k_3^-, \tau) U(\rho,\si;k_1^+, k_3^-, k_2^+, \tau)\, d\nu(\si;k_2^+,  k_3^-,k_1^+, \tau),
\end{split}
\]
where we use
\[
\int_0^\infty g(\si)\, d\nu(\si;k_1^+, k_3^-, k_2^+,\tau) = \int_0^\infty g(\si)\, d\nu(\si;k_1^+,k_2^+,k_3^-,\tau).
\]
Writing the Racah coefficients as Wilson functions or Wilson polynomials, we obtain the theorem.
\end{proof}
We remark that Theorem \ref{thm:symI} is valid for all $k_1,k_3>0$ if appropriate discrete terms are added to the integral, corresponding to the discrete terms in the tensor product decomposition of $\pi^+_{k_1} \tensor \pi^-_{k_3}$.

Let us relabel
\[
(a,b,c,d) = (k_2+i\tau, k_1+k_3-\hf, k_1-k_3-\hf, 1-k_2+i\tau),
\]
then Theorem \ref{thm:symI} can be written as 
\[
\big(\mathcal F W_j(\cdot\,;a,b,c,1-d)\big)(\rho) = (-1)^j W_j(\rho;\tilde a, \tilde b, \tilde c, 1-\tilde d).
\]
This is \cite[Theorem 6.7]{Gr03}, which states that the Wilson transform I maps a basis of orthogonal polynomials in $\mathcal M(a,b,c,d)$ to a basis of orthogonal polynomials in $\mathcal M(\tilde a, \tilde b, \tilde c, \tilde d)$.

\subsection{Formulas related to 4-fold tensor products} \label{ssec:4fold}
In general the Racah coefficients, or $6j$-symbols, are the recoupling coefficients that connect the coupled basis vectors for $(\pi_1 \tensor \pi_2) \tensor \pi_3$ with the coupled basis vectors for $\pi_1 \tensor (\pi_2 \tensor \pi_3)$. The $9j$-symbols are the recoupling coefficients related to a $4$-fold tensor product. These $9j$-symbols can be expressed in terms of the $6j$-symbols. In this section we consider $9j$-symbols that connect the coupled basis vectors for $[(\pi_1 \tensor \pi_2) \tensor \pi_3] \tensor \pi_4$ with the coupled basis vectors for $\pi_1 \tensor [\pi_2 \tensor (\pi_3 \tensor \pi_4)]$. There are several possible transitions from the first coupled basis to the second coupled bases. We consider the following two transitions:
\begin{equation} \label{eq:BEid}
\begin{split}
[(\pi_1 \tensor \pi_2) \tensor \pi_3] \tensor \pi_4 &\xrightarrow{\mathcal I} (\pi_1 \tensor \pi_2) \tensor [\pi_3 \tensor \pi_4] \xrightarrow{\mathcal I} \pi_1 \tensor (\pi_2 \tensor [\pi_3 \tensor \pi_4]),\\
[(\pi_1 \tensor \pi_2) \tensor \pi_3] \tensor \pi_4 &\xrightarrow{\mathcal I} [\pi_1 \tensor (\pi_2 \tensor \pi_3)] \tensor \pi_4 \\
&\xrightarrow{\mathcal I} \pi_1 \tensor [(\pi_2 \tensor \pi_3) \tensor \pi_4] \xrightarrow{\mathcal I} \pi_1 \tensor [\pi_2 \tensor (\pi_3 \tensor \pi_4)].
\end{split}
\end{equation}
These possibilities give rise to two different expressions for the $9j$-symbols, so we obtain a summation or integral formula relating different Racah coefficients. Biedenharn \cite{Bie53} and Elliott \cite{Ell53} used a similar method to find the famous Biedenharn-Elliott identity, or pentagon identity,  for the $\su(2)$-Racah coefficients.

Let us consider the tensor product $\pi^+_{k_1} \tensor \pi^+_{k_2} \tensor \pi^-_{k_3} \tensor \pi^+_{k_4}$. The two transitions from \eqref{eq:BEid} lead to the following identity for operators:
\[
\begin{split}
\mathcal I\big([k_1^+, k_2^+]&,(\si,k_4-k_3)^P,\tau\big)\,\circ  \mathcal I\big([(k_1+k_2+j)^+, k_3^-], k_4^+,\tau\big) \\
= &\,C\, \mathcal I\big([k_2^+, k_3^-], k_4^+,\zeta\big)\circ \mathcal I\big([k_1^+, (\om,k_2-k_3)^P],k_4^+,\tau\big) \circ \mathcal I\big([k_1^+,k_2^+],k_3^-,\rho\big),
\end{split}
\]
where the factor $C$ is still to be determined.  

This factor can be determined in the same way as in the proof of Theorem \ref{thm:symI}. Then we find that the factor we are looking for can be determined from
\[
C\, G_n^{r-j}(\tau;k_4,\rho,k_1+k_2+j-k_3) = G_n^{r}(\tau;k_4,\rho,k_1+k_2-k_3).
\]
and this gives
\[
C = \sqrt{ \frac{(k_1+k_2-k_3+\hf \pm i\rho)_j}{ (k_1+k_2-k_3+k_4+\hf \pm i\tau)_j} }.
\]
Then, similar as in the proof of Theorem \ref{thm:symI}, we obtain the following identity.
\begin{thm} \label{thm:BE++-+}
For $ \zeta,\rho, \si, \tau \in \R$, $k_1,k_2,k_3,k_4>0$, $j \in \Z_{\geq 0}$, $k_2+k_3-\hf>0$, and $k_2-k_3+\hf>0$, 
\[
\begin{split}
\int_0^\infty &B(\om) W_j(\om;k_2-k_3+\hf, k_2+k_3-\hf, k_1+i\rho, k_1-i\rho)\\
& \times \phi_{\si}(\om; k_4+i\zeta,k_2+k_3-\hf, k_2-k_3+\hf, 1-k_4+i\zeta) \\
& \times \phi_{\tau}(\om;k_1+i\rho, k_4+i\zeta, k_4-i\zeta, 1-k_1+i\rho) d\om \\
= &\ W_j(\zeta ; k_2+i\si, k_2-i\si, k_1+i\tau, k_1-i\tau)\\
&\times \phi_{\si}(\rho;k_4+i\tau, k_1+k_2+k_3+j-\hf, k_1+k_2-k_3+j+\hf, 1-k_4+i\tau), \end{split}
\]
where
\[
\begin{split}
B(\om) = \frac1{2\pi} & \frac{   \Ga(k_2+k_3-\hf\pm i\om) \Ga(k_2-k_3+\hf\pm i\om) \Ga(k_1+i\rho\pm i\om)}{ \Ga(k_1+k_2+k_3+j-\hf\pm i\rho) \Ga(k_1+k_2-k_3+j+\hf\pm i\rho)} \\
\times\ & \frac{  \Ga(k_1-i\rho\pm i\om) \Ga(k_4+i\zeta \pm i\om) \Ga(k_4-i\zeta\pm i\om) }{ \Ga(\pm 2i\om) }  .
\end{split}
\]
\end{thm}
Note that we used the symmetry relation Lemma \ref{lem2}(i) for the Wilson functions. If discrete mass points with the appropriate residues at the poles of the weight $B(\om)$ are added to the integral, Theorem \ref{thm:BE++-+} is valid for $k_2,k_3>0$. 

Using the orthogonality relations for the Wilson polynomials, it follows that Theorem \ref{thm:BE++-+} is equivalent to the following bilinear summation formula for the Wilson polynomials:
\[
\begin{split}
\sum_{j=0}^\infty&  \frac{2k_1+2k_2+2j-1}{2k_1+2k_2-1}\, \frac{ (2k_1+2k_2-1)_j }{j!\,(2k_1)_j (2k_2)_j}\\
& \times \phi_{\si}(\rho;k_4+i\tau, k_1+k_2+k_3+j-\hf, k_1+k_2-k_3+j+\hf, 1-k_4+i\tau) \\
& \times W_j(\om;k_2-k_3+\hf, k_2+k_3-\hf, k_1+i\rho, k_1-i\rho) \\
& \times W_j(\zeta ; k_2+i\si, k_2-i\si, k_1+i\tau, k_1-i\tau)\\
=&K\,\phi_{\si}(\om; k_4+i\zeta,k_2+k_3-\hf, k_2-k_3+\hf, 1-k_4+i\zeta) \\
& \times \phi_{\tau}(\om;k_1+i\rho, k_4+i\zeta, k_4-i\zeta, 1-k_1+i\rho),
\end{split}
\]
where
\[
K = \frac{\Ga(2k_1) \Ga(2k_2) \Ga(k_4+i\zeta \pm i\om) \Ga(k_4-i\zeta \pm i\om) }{ \Ga(2k_1+2k_2)}.
\]
This sum can be considered as the $q=1$ analogue of the bilinear summation formula for Askey-Wilson polynomials \cite[Theorem 3.3]{KVdJ99} which is obtained by a different method.

In the same way as we obtained Theorem \ref{thm:BE++-+}, we can find Biedenharn-Elliott identities for other tensor products.  From  $\pi^+_{k_1} \tensor \pi^+_{k_2} \tensor \pi^+_{k_3} \tensor \pi^-_{k_4}$ we obtain a generalized convolution identity for the Wilson polynomials, also involving Racah polynomials. 
\begin{thm} \label{thm:BE+++-}
Let $k_1,k_2,k_3,k_4>0$, $j_1,j_2 \in \Z_{\geq 0}$ and $\rho,\zeta \in \R$, then the Wilson polynomials satisfy the following generalized convolution identity
\[
\begin{split}
\sum_{j_3=0}^{j_1+j_2}& c_{j_3} R_{j_3}(j_1;2k_2-1, 2k_3-1, -j_1-j_2-1, 2k_1+2k_2+j_1+j_2-1) \\
&\times W_{j_3}(\rho;k_3-k_4+\hf, k_3+k_4-\hf, k_2+i\zeta, k_2-i\zeta)\\
&\times W_{j_1+j_2-j_3}(\zeta; k_2+k_3-k_4+j_3+\hf, k_2+k_3+k_4+j_3-\hf, k_1+i\tau, k_1-i\tau) \\
=&\,W_{j_2}(\rho;k_3-k_4+\hf, k_3+k_4-\hf, k_1+k_2+j_1+i\tau, k_1+k_2+j_1-i\tau)\\
&\times W_{j_1}(\zeta; k_2+i\rho, k_2-i\rho, k_1+i\tau, k_1-i\tau),\\
\end{split}
\]
where
\[
c_{j_3} = \binom{j_1+j_2}{j_3} \frac{ (2k_2)_{j_1} (2k_2+2k_3)_{2j_3} (2k_1+2k_2+2k_3+j_1+j_2-1)_{j_3}} { (2k_3)_{j_3} (j_3+2k_2+2k_3-1)_{j_3} (2k_2+2k_3)_{j_1+j_2+j_3} }.
\]
\end{thm}
\begin{rem}
This convolution identity is obtained in \cite{LVdJ02} by analytic continuation of the classical Biedenharn-Elliott identity for $\su(2)$ Racah coefficients, which can be considered as a generalized convolution identity for the Racah polynomials. Many generalized convolution identities for hypergeometric orthogonal polynomials can be obtained from Theorem \ref{thm:BE+++-} as limit cases. For instance, if we put $i\tau= k_1+k_2+k_3-k_4-\hf$ in Theorem \ref{thm:BE+++-} and we let $k_1 \rightarrow \infty$, then we obtain a generalized convolution identity for continuous dual Hahn polynomials which is equivalent to Theorem \ref{thm:++-}. See \cite{LVdJ02} and \cite{KVdJ98} for other limit cases.
\end{rem}

The Biedenharn-Elliott identity for the tensor product  $\pi^+_{k_1} \tensor \pi^+_{k_2} \tensor \pi^-_{k_3} \tensor \pi^-_{k_4}$ leads to the following theorem.
\begin{thm}  \label{thm:BE ++--}
For $\rho,\zeta,\tau \in \R$, $k_1,k_2,k_3,k_4>0$, $j_1,j_2 \in \Z_{\geq 0}$, and $k_2+k_3-\hf>0$,
\[
\begin{split}
\int_0^\infty &B(\om) \ W_{j_2}(\om;k_3-k_2+\hf, k_2+k_3-\hf, k_4+i\zeta, k_4-i\zeta)\\
&\times W_{j_1}(\om;k_2-k_3+\hf, k_2+k_3-\hf, k_1+i\rho, k_1-i\rho) \\
&\times \phi_{\tau}(\om;k_1+i\rho, k_4+i\zeta, k_4-i\zeta, 1-k_1+i\rho)\, d\om\\
=&\ W_{j_1}(\zeta;k_2-k_3-k_4-j_2+\hf, k_2+k_3+k_4+j_{2}-\hf, k_1+i\tau, k_1-i\tau) \\
&\times W_{j_2}(\rho;k_3-k_1-k_2-j_1+\hf, k_1+k_2+k_3+j_1-\hf, k_4+i\tau, k_4-i\tau),
\end{split}
\]
where 
\[
\begin{split}
B(\om) =& \frac{1}{2 \pi}  \frac{ \Ga(k_1+k_2+k_3+k_4+j_1+j_2-\hf\pm i\tau)  \Ga(k_2+k_3-\hf \pm i\om)}{ \Ga(k_1+k_2+k_3+j_1-\hf\pm i\rho) \Ga(k_2+k_3+k_4+j_2-\hf\pm i\zeta)} \\ &\times  \frac{\Ga(k_1+i\rho\pm i\om) \Ga(k_1-i\rho\pm i\om) \Ga(k_4+i\zeta\pm i\om) \Ga(k_4-i\zeta\pm i\om)  } {\Ga(\pm 2i\om) }.
\end{split}
\]
\end{thm}
Theorem \ref{thm:BE ++--} can be obtained from Theorem \ref{thm:BE++-+} by setting $\si = i(k_3+k_4+j_2-\hf)$ and using Lemma \ref{lem2}. Note that the Racah coefficient for the tensor product $\pi^+_{k_2} \tensor \pi^-_{k_3} \tensor \pi^-_{k_3}$ is needed to obtain Theorem \ref{thm:BE ++--}. It can be checked, using the algebra involution $\vartheta$, that these Racah coefficients are the same as the Racah coefficients for the tensor product $\pi^+_{k_4} \tensor \pi^+_{k_3} \tensor \pi^-_{k_2}$, which are given in Theorem \ref{thm:++-}.

The Biedenharn-Elliott identities for the tensor products $\pi^+ \tensor \pi^+ \tensor \pi^+ \tensor \pi^P$, $\pi^+ \tensor \pi^+ \tensor \pi^P \tensor \pi^+$, and $\pi^+ \tensor \pi^+ \tensor \pi^P \tensor \pi^-$ lead to identities equivalent to those in Theorems \ref{thm:BE++-+} and \ref{thm:BE+++-}. Finally let us remark that Theorems \ref{thm:BE++-+}, \ref{thm:BE+++-}, \ref{thm:BE ++--} can all be considered as connection coefficient formulas between two orthogonal bases in certain Hilbert spaces.


\section{Racah coefficients for $\U_q(\su(1,1))$} \label{sec:Uq}
Throughout this section we assume that $0 <q <1$. Racah coefficients for the quantized universal enveloping algebra $\U_q(\su(1,1))$ can be calculated in the same way as for $\su(1,1)$. We use notations which make it easy for the reader to compare with the Lie algebra case.

\subsection{The quantized universal enveloping algebra $\U_q(\su(1,1))$}
The quantized universal enveloping algebra $\mathcal U_q = \U_q\big(\su(1,1)\big)$ is the unital, associative, complex algebra generated by $K$, $K^{-1}$, $E$, and $F$, subject to the relations
\[
K K^{-1} = 1 = K^{-1}K, \quad KE = qEK, \quad KF= q^{-1}FK, \quad EF-FE= \frac{K^2-K^{-2}}{q-q^{-1}}.
\]
The Casimir element
\begin{equation} \label{eq:Casimir}
\Om = \frac{q^{-1} K^2 +qK^{-2}-2}{(q^{-1}-q)^2} + EF= \frac{q^{-1}K^{-2}+qK^2-2}{(q^{-1}-q)^2} +FE
\end{equation}
is a central element of $\U_q$. The algebra
$\U_q$ is a Hopf $*$-algebra with comultiplication $\De$ defined on the generators by
\begin{equation} \label{eq:comult}
\begin{split}
\De(K) = K \tensor K,&\quad  \De(E)= K \tensor E + E \tensor K^{-1}, \\
\De(K^{-1}) = K^{-1} \tensor K^{-1},&\quad  \De(F) = K \tensor F + F \tensor K^{-1} ,
\end{split}
\end{equation}
and $\De$ can be extended to $\mathcal U_q$ as an algebra homomorphism. Observe here that the coproduct is not cocommutative. The $*$-structure is defined by
\[
K^*=K, \quad E^*=-F, \quad F^* = -E, \quad (K^{-1})^* = K^{-1}.
\]
The irreducible $*$-representations of $\U_q$ have been determined by Burban and Klimyk \cite{BK93}. There are, besides the trivial representation, five classes of irreducible $*$-representations:\\

\emph{Positive discrete series.} The positive discrete series $\pi^+_k$
are labelled by $k>0$. The representation space is $\ell^2(\Z_{\geq 0})$
with orthonormal basis $\{e_n\}_{n \in \Z_{\geq 0} }$. The action is
given by
\begin{equation} \label{eq:qpos}
\begin{split}
\pi^+_k(K)\, e_n &= q^{k+n}\, e_n, \quad \pi^+_k(K^{-1})\, e_n =q^{-(k+n)}\,e_n,  \\ 
(q^{-1}-q)\pi^+_k(E)\, e_n &= q^{-\hf-k-n} \sqrt{(1-q^{2n+2})(1-q^{4k+2n})}\, e_{n+1},  \\
(q^{-1}-q)\pi^+_k(F)\, e_n &= -q^{\hf-k-n}
\sqrt{(1-q^{2n})(1-q^{4k+2n-2})}\, e_{n-1},
\\ (q^{-1}-q)^2 \pi^+_k(\Om)\, e_n &= (q^{2k-1}+q^{1-2k}-2)\, e_n.
\end{split}
\end{equation}

\emph{Negative discrete series.} The negative discrete series
$\pi^-_k$ are labelled by $k>0$. The representation space is
$\ell^2(\Z_{\geq 0})$ with orthonormal basis $\{e_n\}_{n \in \Z_{\geq
0}}$. The action is given by
\begin{equation} \label{eq:qneg}
\begin{split}
\pi^-_k(K)\, e_n &= q^{-(k+n)}\, e_n, \quad \pi^-_k(K^{-1})\, e_n =q^{k+n}\,e_n,
 \\  (q^{-1}-q)\pi^-_k(E)\, e_n &= -q^{\hf-k-n}
\sqrt{(1-q^{2n})(1-q^{4k+2n-2})}\, e_{n-1},
\\(q^{-1}-q)\pi^-_k(F)\, e_n &= q^{-\hf-k-n}
\sqrt{(1-q^{2n+2})(1-q^{4k+2n})}\, e_{n+1}, \\ (q^{-1}-q)^2
\pi^-_k(\Om)\, e_n &= (q^{2k-1}+q^{1-2k}-2)\,e_n.
\end{split}
\end{equation}

\emph{Principal unitary series.} The principal unitary series
representations $\pi^P_{\rho,\eps}$ are labelled by $0\leq\rho \leq
-\frac{\pi}{2 \ln q}$ and $\eps \in [0,1)$, where $(\rho,\eps)
\neq (0,\hf)$. The representation space is $\ell^2(\Z)$ with
orthonormal basis $\{e_n\}_{n \in \Z}$. The action is given by
\begin{equation} \label{eq:qprinc}
\begin{split}
\pi^P_{\rho,\eps}(K)\, e_n &= q^{n+\eps}\, e_n, \quad
\pi^P_{\rho,\eps} (K^{-1})\, e_n = q^{-(n+\eps)}\,e_n,\\  
(q^{-1}-q)\pi^P_{\rho,\eps}(E)\, e_n
&= q^{-\hf-n-\eps} \sqrt{(1-q^{2n+2\eps+2i\rho+1})
(1-q^{2n+2\eps-2i\rho+1})}\, e_{n+1},
\\(q^{-1}-q)\pi^P_{\rho,\eps}(F)\, e_n &= -q^{\hf-n-\eps}
\sqrt{(1-q^{2n+2\eps+2i\rho-1}) (1-q^{2n+2\eps-2i\rho-1})}\, e_{n-1},\\
(q^{-1}-q)^2 \pi^P_{\rho,\eps} (\Om)\, e_n &=(q^{2i\rho}+q^{-2i\rho} -2)\,e_n.
\end{split}
\end{equation}
For $(\rho,\eps)= (0,\hf)$ the representation $\pi^P_{0,\hf}$ splits into a direct
sum of a positive and a negative discrete series representation: $\pi^P_{0,\hf}= \pi^+_\hf
\oplus \pi^-_\hf$. The representation space splits into two invariant subspaces: $\{e_n \, |\,
n<0 \} \oplus \{ e_n\, | \, n \geq 0 \}$.\\

\emph{Complementary series.} The complementary series
representations $\pi^C_{\la,\eps}$ are labelled by $\la$ and
$\eps$, where $\eps \in [0,\hf)$ and $\la \in (-\hf,-\eps)$, or
$\eps \in (\hf,1)$ and $\la \in (-\hf, \eps-1)$. The
representation space is $\ell^2(\Z)$ with orthonormal basis
$\{e_n\}_{n \in \Z}$. The action of the generators is the same as for the principal unitary series, with $i\rho$ replaced by $\la+\hf$.\\

The fifth class consists of the strange series representations. The actions of the generators in the strange series are the same as in principal unitary series with $i\rho$ replaced by $- \frac{ i\pi}{2\ln q}+\si$, $\si>0$. We do not need the strange series representation in this paper.\\

\subsection{Clebsch-Gordan coefficients}
The Clebsch-Gordan coefficients for $\U_q$ that we need in this paper are $q$-analogues of the (continuous) dual Hahn polynomials, which we needed to describe the Clebsch-Gordan coefficients for $\su(1,1)$.

The $q$-Racah polynomial \cite{AW79} are defined by
\[
R_n(x;\al,\be,\ga,\de;q) = \ph{4}{3}{q^{-n}, \al \be q^{n+1}, q^{-x}, \ga \de q^{x+1}}{\al q, \be \de q, \ga q}{q,q},
\]
where one of the lower parameters is equal to $q^{-N}$ for a non-negative integer $N$. These are polynomials of degree $n$ in in $q^{-x}+\ga \de q^{x+1}$.

The dual $q$-Hahn polynomials are $q$-Racah polynomials with $\al q=q^{-N}$ and $\be=0$. So the dual $q$-Hahn polynomial are defined by
\begin{equation} \label{def:dual qHahn}
T_n(x;\ga,\de,N;q)=  \ph{3}{2}{q^{-n}, q^{-x}, \ga \de q^{x+1}}{\ga q ,q^{-N}}{q,q},
\end{equation}
where $N$ is a non-negative integer, and $n \in \{0,1, \ldots, N\}$. Under certain condition on $\ga$ and $\de$ these polynomials are orthogonal on the finite set $\{0,1,\ldots,N\}$ with respect to a positive measure. The orthonormal dual $q$-Hahn polynomials are defined by
\[
t_n(x) = t_n(x;\ga,\de,N;q) =  (\ga \de q)^{-\hf n}\sqrt{ \frac{ (\ga q, q^{-N};q)_n }{ (q,q^{-N}/\de;q)_n} }\, T_n(x;\ga,\de,N;q),
\]
and the orthogonality relations read, for $0<\ga<q^{-1}$ and $0<\de<q^{-1}$ or for $\ga>q^{-N}$ and $\de>q^{-N}$, 
\[
\begin{split}
&\sum_{x=0}^N t_m(x;\ga,\de,N;q)t_n(x;\ga,\de,N;q) w(x;\ga,\de,N;q) = \de_{mn},\\
&w(x;\ga,\de,N;q) = (-1)^x(\ga q)^{N-x}q^{Nx-\hf x(x-1)}\frac{ (1-\ga \de q^{2x+1}) (\ga q, q^{-N};q)_x  (\de q;q)_N}{(q, \de q ;q)_x (\ga \de q^{x+1};q)_{N+1}  }.
\end{split}
\]
The dual $q$-Hahn polynomials satisfy the following three-term recurrence relation:
\[
(q^{-x}+\ga \de q^{x+1})\, t_n(x)= a_n t_{n+1}(x) + b_n t_n(x) + a_{n-1} t_{n-1}(x),
\]
where
\[
\begin{split}
a_n &= \sqrt{ \ga \de q(1-q^{n+1}) (1-\ga q^{n+1}) (1-q^{n-N}/\de)(1-q^{n-N}) },\\
b_n &= q^{n-N}+ \ga q^{n+1} + \ga q^{n-N} + \ga\de q^{n+1}- \ga q^{2n-N}(1+q).\\
\end{split}
\]

The Askey-Wilson polynomials \cite{AW85}, see also \cite{GRa04}, \cite{KS98}, are polynomials in $x+x^{-1}$ defined by
\begin{equation} \label{def:AW pol}
p_n(x)=p_n(x;a,b,c,d|q)= a^{-n} (ab,ac,ad;q)_n \ph{4}{3}{q^{-n}, abcdq^{n-1}, ax, a/x } {ab, ac, ad} {q, q}.
\end{equation}
By Sear's $_4\varphi_3$ transformation formula \cite[(III.15)]{GRa04} the polynomials $p_n$ are symmetric in the parameters $a$, $b$, $c$ and $d$. Let the parameters $a,b,c,d$ be such that non-real parameters occur in pairs of complex conjugates, and the pairwise product of two real parameters is smaller than $1$. Then the Askey-Wilson polynomials satisfy the following orthogonality relations,
\[
\frac{1}{4\pi i}\int_\T p_m(x)p_n(x) w(x) \frac{dx}{x} + \sum_{\substack{ k \in \Z_{\geq 0}\\ |eq^k|>1}} p_m(eq^k) p_n(eq^k) \Res{x=eq^k} \frac{w(x)}{x} = \de_{mn}\, h_n.
\]
Here $\T$ is the counter-clockwise oriented unit circle in the complex plane, $e$ is any of the real parameters, and
\[
\begin{split}
w(x)&= \frac{ (x^{\pm 2};q)_\infty }{ (ax^{\pm 1}, bx^{\pm 1}, cx^{\pm 1}, dx^{\pm 1};q)_\infty },\\
h_n &= \frac{ (abcdq^{n-1};q)_n (abcdq^{2n};q)_\infty }{ (q^{n+1}, abq^n , acq^n, adq^n, bcq^n, bdq^n, cdq^n;q)_n }.
\end{split}
\]
Observe that the Askey-Wilson polynomials can be considered as $q$-analogues of the Wilson polynomials.

The continuous dual $q$-Hahn polynomials $S_n$ are Askey-Wilson polynomials with $d=0$;
\begin{equation} \label{def:cdqHahn}
S_n(x;a,b,c|q)= p_n(x;a,b,c,0|q)= a^{-n} (ab,ac;q)_n \ph{3}{2}{q^{-n}, ax, a/x} {ab,ac} {q,q}.
\end{equation}
These polynomials are symmetric in the parameters $a,b,c$. Let us assume that either $a,b,c$ are real and $|ab|, |ac|, |bc|<1$, or one of the parameters is real and the other two parameters are complex conjugates with positive real part. In this case the polynomials
\[
s_n(x)=s_n(x;a,b,c|q) = \frac{S_n(x;a,b,c|q)}{\sqrt{(q,ab,ac;q)_n}}
\]
are orthonormal with respect to the measure $\mu(\cdot;a,b,c|q)$ defined by
\begin{gather}
\int f(x)\,d\mu(x;a,b,c|q) =\frac{1}{4\pi i}\int_\T f(x) w(x;a,b,c|q)\, \frac{dx}{x} + \sum_{\substack{k \in \Z_{\geq 0}\\ |eq^k|>1}} f(eq^k) w_k(a;b,c;q), \label{eq:meas cdqH}\\ 
w(x;a,b,c|q) =  \frac{(q,ab,ac,bc;q)_\infty  (x^{\pm 2};q)_\infty } { (ax^{\pm 1}, bx^{\pm 1}, cx^{\pm 1};q)_\infty }, \nonumber
\end{gather}
where $e$ is any of the real parameters $a$, $b$, $c$. If we assume $e=a$, we have
\[
\begin{split}
w_k(a;b,c;q) &= \Res{x=aq^k} \frac{w(x;a,b,c|q)}{x} \\
&=\frac{1-a^2q^{2k}}{1-a^2} \frac{ (a^{-2}, bc;q)_\infty }{(b/a, c/a;q)_\infty} \frac{(a^2,ab,ac;q)_k} {(q, aq/b, aq/c ; q)_k} (-1)^k q^{-\hf k(k-1)} (a^2bc)^{-k}.
\end{split}
\]
The polynomials $s_n$ form an orthonormal basis for the Hilbert space $\mathcal L(\mu)$ consisting of functions $f$ that satisfy $f(x)=f(x^{-1})$ ($\mu$-almost everywhere), with the inner product associated to the measure $\mu$.
The polynomials $s_n$ satisfy the three-term recurrence relation
\begin{equation} \label{eq:ttr cdqH}
(x+x^{-1} s_n(x) = a_n s_{n+1}(x) + b_n s_n(x) + a_{n-1} s_{n-1}(x),
\end{equation}
where
\begin{align*}
a_n & = \sqrt{ (1-q^{n+1})(1-abq^{n})(1-acq^n)(1-bcq^n) }, \\ b_n &= a+a^{-1}-a^{-1}(1-abq^n)(1-acq^n) - a(1-q^n)(1-bcq^{n-1}).
\end{align*}

\* 
First we consider the tensor product $\pi^+_{k_1} \tensor \pi^+_{k_2}$. This representation can be decomposed into irreducible representations as
\begin{equation} \label{eq:q++}
\pi^+_{k_1} \tensor \pi^+_{k_2} \cong \bigoplus_{j=0}^\infty \pi^+_{k_1+k_2+j}.
\end{equation}
The Clebsch-Gordan coefficients can be obtained by identifying the action of the Casimir operator on a suitable subspace of $\ell^2(\Z_{\geq 0})^{\tensor 2}$ with the Jacobi operator corresponding to the dual $q^2$-Hahn polynomials, analogously to Proposition \ref{prop I}. See also \cite{KR89}, \cite{Va89}, \cite{KK89}, \cite{KVdJ98} for the interpretation of the (dual) $q$-Hahn polynomials as  Clebsch-Gordan coefficients. Let us define the functions
\[
v_n^p(x;k_1,k_2) =(-1)^n q^{p(2k_1+x)}\sqrt{ \frac{ (-1)^x (q^{-2p};q^2)_x (q^{4k_2};q^2)_p}{ (q^{4k_1+4k_2+2x};q^2)_p}}\ t_n(x;q^{4k_1-2}, q^{4k_2-2}, p;q^2).
\]
These functions are orthonormal on $\{0,1, \ldots, p\}$ with respect to the weight
\[
W(x;k_1,k_2) = q^{4k_1 x} q^{-x(x-1)} \frac{ (1-q^{4k_1+4k_2+4x-2}) (q^{4k_1};q^2)_x }{(1-q^{4k_1+4k_2+2x-2}) (q^2,q^{4k_2};q^2)_x}.
\]
The intertwiner for \eqref{eq:q++} can now be defined by
\[
\begin{split}
I': \ell^2(\Z_{\geq 0})^{\tensor 2} &\rightarrow \bigoplus_{x=0}^\infty \ell^2(\Z_{\geq 0})\, W(x;k_1,k_2),\\
e_{n_1} \tensor e_{n_2} &\mapsto \sum_{x=0}^{n_1+n_2} v_{n_1}^{n_1+n_2}(x;k_1,k_2)\,e_{n_1+n_2-x} ,
\end{split} 
\]
The intertwining property of $I'$ can be checked by writing the Clebsch-Gordan coefficients as $_3\varphi_2$-functions, and using the following contiguous relations:
\begin{equation} \label{eq:cont rel 3phi2}
\begin{split}
&\frac{ (1-b/q)(1-c/q) }{ (1-d/q)(1-e/q) q^n } \ph{3}{2}{q^{-n}, b,c}{d,e}{q,q} =\\ & \qquad 
\ph{3}{2}{q^{-n}, b/q,c/q}{d/q,e/q}{q,q}   - \ph{3}{2}{q^{-n-1}, b/q,c/q}{d/q,e/q}{q,q}, \\& \\
&\ph{3}{2}{q^{-n}, b,c}{d,e}{q,q} =  q^{-n}\frac{(1-dq^{n})(1-eq^n)}{(1-d)(1-e)} \ph{3}{2}{q^{-n}, bq,cq}{dq,eq}{q,q} \\
&\qquad + \frac{(1-q^n)(1-de q^{n-1}/bc) }{(1-d)(1-e)}\ph{3}{2}{q^{1-n}, bq,cq}{dq,eq}{q,q}.
\end{split}
\end{equation}
These relations follow from \cite[(2.4)]{GIM96}, respectively \cite[(2.3)]{GIM96}, and the transformation
\[
\ph{3}{2}{q^{-n}, b,c}{d,e}{q,q} = \frac{ (q^{-n},b,c;q)_n }{ (q,d,e;q)_n} \ph{3}{2}{q^{-n}, q^{1-n}/d,q^{1-n}/e}{q^{1-n}/b,q^{1-n}/c}{q,\frac{ deq^n}{bc}},
\]
which is obtained from reversing the summation for the $_3\varphi_2$-series.\\

Next we consider the tensor product representation $\pi^+_{k_1} \tensor \pi^-_{k_2}$. The action of $\pi^+_{k_1} \tensor \pi^-_{k_2}( \De(\Om))$ on a suitable subspace of $\ell^2(\Z_{\geq 0})^{\tensor 2}$ can be identified with a Jacobi operator for continuous dual $q^2$-Hahn polynomials. We define for $\rho \in [0,-\pi/2\ln q]$,
\[
F_n^p(\rho;k_1,k_2) = 
\begin{cases}
\displaystyle (-1)^{p} \sqrt{ \frac{ (q^{2k_2-2k_1+1 \pm 2i\rho};q^2)_{-p} }{ (q^2, q^{4k_2};q^2)_{-p} }} \\
\quad \times s_{n}(q^{2i\rho};q^{2k_2-2k_1-2p+1} ,q^{2k_1+2k_2-1} , q^{2k_1-2k_2+1} |q^2), & p \leq 0,\\
\displaystyle \sqrt{ \frac{ (q^{2k_1-2k_2+1 \pm 2i\rho};q^2)_p }{ (q^2, q^{4k_1};q^2)_p }} \\
\quad \times s_{n}(q^{2i\rho};q^{2k_1-2k_2+2p+1}, q^{2k_1+2k_2-1}, q^{2k_2-2k_1+1}|q^2), & p \geq 0, 
\end{cases}
\]
and similarly for $q^{2i\rho}$ in the discrete part of the support of the corresponding orthogonality measure. We denote the orthonormality measure for the functions $F_n^p$ by $\mu_1^p(\cdot;k_1,k_2)$. So, assuming that $k_1$ and $k_2$ are such that $\mu_1^p$ is absolutely continuous, we have
\[
\begin{split}
\int f(\rho)\,& d\mu_1^p(\rho;k_1,k_2) =\\ 
&\frac{1}{2\pi} \int_0^{-\frac{\pi}{2\ln q}} f(\rho) \frac{ (q^2,q^{4k_1}, q^{4k_2}, q^{\pm 4i\rho};q^2)_\infty} {(q^{2k_1-2k_2+1 \pm 2i\rho}, q^{2k_2-2k_1+1 \pm 2i\rho}, q^{2k_1+2k_2-1 \pm 2i\rho};q^2)_\infty} \, d\rho
\end{split}
\]
As in the Lie algebra case we will define the measure $\mu_1(\cdot;k_1,k_2)$ to be the measure $\mu_1^p(\cdot;k_1,k_2)$ with the largest possible support. We define the unitary operator $I_1'$ by
\[
\begin{split}
I_1' : \ell^2(\Z_{\geq 0})^{\tensor 2} &\rightarrow \qdirint \ell^2(\Z) d\mu_1(\rho;k_1,k_2) \\
e_{n_1} \tensor e_{n_2} & \mapsto \int_0^{-\frac{\pi}{2\ln q}} F_n^{n_1-n_2}(\rho;k_1,k_2) \, e_{n_1-n_2} \, d\rho,
\end{split}
\]
where $n=\min\{n_1,n_2\}$. With this intertwining operator we find the following decomposition, cf.~\cite{KM96}, \cite{Gr02},
\[
\pi^+_{k_1} \tensor \pi^-_{k_2} \cong \qdirint \pi^P_{\rho,k_1-k_2} d\rho.
\]
This can be proved by checking the actions of $K$, $K^{-1}$, $E$, $F$, using the contiguous relations \eqref{eq:cont rel 3phi2}. If $d\mu_1$ has discrete mass points, then we must replace $\pi^P_{\rho,k_1-k_2}$ in the discrete mass points by a discrete series representation or a complementary series representation. See \cite[Theorem 2.4]{Gr02} for a precise formulation.
\begin{rem} \label{rem:q -+}
Different from the Lie algebra case, we cannot obtain the tensor product decomposition of $\pi^-_{k_1} \tensor \pi^+_{k_2}$ from the decomposition of $\pi^+_{k_1}\tensor \pi^-_{k_2}$, because the coproduct is not cocommutative. When we would try our approach to find the decomposition for $\pi^-_{k_1} \tensor \pi^+_{k_2}$, we encounter the following problems. Let $\mathcal U_q^{\mathrm{opp}}$ be the Hopf-algebra $\U_q$ with opposite comultiplication $\De^{\text{opp}}$. Then $\U_q^{\mathrm{opp}}$ is isomorphic to $\U_{q^{-1}}$ as a Hopf-algebra, where the isomorphism is the interchanging of $E$ and $F$. So $\pi^-_{k_1} \tensor \pi^+_{k_2} (\De(\Om))$ can be identified with the Jacobi operator for the continuous dual $q^{-2}$-Hahn polynomials. This Jacobi operator corresponds to an indeterminate moment problem, therefore $\pi^-_{k_1} \tensor \pi^+_{k_2}(\De(\Om))$ is \emph{not} essentially self-adjoint. We refer to the book by Akhiezer \cite{Ak65} for more information about moment problems. 
\end{rem}

The action of the Casimir in the tensor product representation $\pi^+_k \tensor \pi^P_{\rho,\eps}$ on a suitable subspace of $\ell^2(\Z_{\geq 0}) \tensor \ell^2(\Z)$ can again be identified with a Jacobi operator corresponding to continuous dual $q^2$-Hahn polynomials. We define 
\[
\begin{split}
&G_n^p(\tau;k,\rho,\eps) = \\
&\begin{cases}
\displaystyle (-1)^n \sqrt{ \frac{q^{4pk}\,(q^{1+2\eps \pm 2i\rho};q^2)_{p} }{ (q^{1+2k+2\eps  \pm 2i\tau};q^2)_{p}}}\ s_n(q^{2i\tau}; q^{2k+2i\rho}, q^{2k-2i\rho}, q^{1-2k-2p-2\eps}|q^2),& \tau \in [0,-\frac{\pi}{2\ln q}),\\
\displaystyle (-1)^n \sqrt{ \frac{q^{4pk}\,(q^{1+2\eps \pm 2i\rho};q^2)_{p} }{ (q^{2+2j},q^{4k+4\eps-2j};q^2)_{p}}}\ s_n(q^{2\tau_j}; q^{2k+2i\rho}, q^{2k-2i\rho}, q^{1-2k-2p-2\eps}|q^2),& \tau_j <0,
\end{cases}
\end{split}
\]
where $\tau_j= \hf-k-\eps+j$ for $j \in \Z$. The functions $G_n^p$ are orthonormal with respect to a measure $d\mu_2^p(\cdot;k,\rho,\eps)$, given by
\[
\int f(\tau)\, d\mu_2(\tau;k,\rho,\eps) =  \int_0^{-\frac{\pi}{2\ln q}} f(\tau) w_2(\tau;k,\rho,\eps)\, d\tau + \sum_{\substack{j \in \Z_{\geq -p}\\\tau_j<0}} f(\tau) w_2(\tau_j;k,\rho,\eps),
\]
where
\begin{align*}
w_2(\tau;k,\rho,\eps) =& \frac{1}{2 \pi} \frac{ (q^2, q^{4k}, q^{1-2\eps \pm 2i\rho}, q^{\pm 4i\tau};q^2)_\infty }{ (q^{2k+2i\rho \pm 2i\tau}, q^{2k-2i\rho \pm 2i\tau}, q^{1-2k-2\eps \pm 2i\tau};q^2)_\infty }, & \textstyle \tau \in [0, -\frac{\pi}{2 \ln q} ),\\
w_2(\tau_j;k,\rho,\eps) =& \frac{ 1-q^{2-4k-4\eps+4j}}{1-q^{2-4k-4\eps}}\, \frac{ (q^{4k+4\eps-2}, q^{4k};q^2)_\infty }{ (q^{4k+2\eps-1 \pm 2i\rho};q^2)_\infty } \frac{ (q^{2-4k-4\eps}, q^{1-2\eps \pm 2i\rho};q^2)_j }{ (q^2, q^{3-4k-2\eps \pm 2i\rho};q^2)_j }&\\
& \times (-1)^j q^{-j(j-1)} q^{-j(2-4\eps)}, & \tau=\tau_j \in \R_{<0}.
\end{align*}
We see that the number of discrete mass points of the measure $d\mu_2^p$ tends to inifinity as $p \rightarrow \infty$. We denote the measure we obtain in this way, i.e.~with an infinite number of discrete mass points, by $d\mu_2$.
The intertwining operator $I_2'$ is now given by
\[
\begin{split}
I_2' : \ell^2(\Z_{\geq 0}) \tensor \ell^2(\Z) &\rightarrow \qdirint \ell^2(\Z) \, w_2(\tau;k,\rho,\eps)\, d\tau \oplus \bigoplus_j \ell^2(\Z_{\geq 0})\, w_2(\tau_j;k,\rho,\eps)\\
e_{n_1} \tensor e_{n_2} &\mapsto \int_0^{-\frac{ \pi}{2\ln q} } G_{n_1}^{n_1+n_2}(\tau;k,\rho,\eps)\, e_{n_1+n_2}\, d\tau + \sum_j G_{n_1}^{n_1+n_2}(\tau_j;k,\rho,\eps) \, e_{n_1+n_2-j} .
\end{split}
\]
This leads to the following decomposition of the tensor product $\pi^+_k \tensor \pi^P_{\rho,\eps}$,
\[
\pi^+_k \tensor \pi^P_{\rho,\eps} \cong \displaystyle \qdirint \pi^P_{\si,\eps'}\, d\si \oplus \bigoplus_{j} \pi^+_{\eps'+j}, \qquad \eps' = \eps+k.
\]
This is proved by checking the intertwining property of $I_2'$ for the actions of the generators of $\U_q$, using the contiguous relations \eqref{eq:cont rel 3phi2}.  

The Clebsch-Gordan coefficients for the tensor products $\pi^+_k \tensor \pi^C_{\la,\eps}$ and $\pi^+_k \tensor\pi^S_{a,\eps}$ can be obtained in the same way as for $\pi^+_k \tensor \pi^P_{\rho,\eps}$. Formally the results follow from the substitutions $\rho \mapsto -i(\la+\hf)$, respectively $\rho\mapsto -i(- \frac{ i\pi}{2\ln q}+a)$.\\

In order to find the decomposition of the tensor product representation $\pi^P_{\rho,\eps} \tensor \pi^-_k$ we introduce the algebra involution $\vartheta$, defined on the generators by
\begin{equation} \label{eq:map theta}
\vartheta(K)=K^{-1}, \quad \vartheta(E)=F, \quad \vartheta(F)=E, \quad \vartheta(K^{-1})=K,
\end{equation}
and $\vartheta$ is extended to $\U_q$ as an algebra mophism. From \eqref{eq:qpos} and \eqref{eq:qneg} it follows that 
$\pi^+_k\big(\vartheta(X)\big) = \pi^-_k(X)$ for $X \in \mathcal U_q$. Here we identify the representation space for $ \pi^+_k$ with the representation space for $\pi^-_k$. Furthermore, if we define a unitary operator $U: \ell^2(\Z) \rightarrow \ell^2(\Z)$ by $Ue_n =  (-1)^n e_{-n}$. Then by \eqref{eq:qprinc} we have $U \circ \pi^P_{\rho,-\eps}\circ \vartheta=\pi^P_{\rho,\eps}\circ U $.

A direct verification shows that $\De \circ \vartheta = \vartheta \tensor \vartheta \circ \De^{\text{opp}}$, so we have
\[
\pi^P_{\rho,\eps} \tensor \pi^-_{k} \circ \De \cong \pi^+_k \tensor \pi^P_{\rho,-\eps} \circ \De \circ \vartheta
\]
Since $\vartheta(\Om)= \Om$ and the Clebsch-Gordan coefficients are eigenfunctions of $\De(\Om)$, the Clebsch-Gordan coefficients for the tensor product $\pi^P_{\rho,\eps}\tensor \pi^-_k $ can be obtained from $\pi^+_k \tensor \pi^P_{\rho,\eps}$.  Then in the same way as for the Lie algebra case, see subsection \ref{ssec:CGC P-}, we find that the Clebsch-Gordan coefficients for $\pi^P_{\rho,\eps} \tensor \pi^-_{k} $ are the functions $(-1)^{n+p} G_{n}^{-p}(\si;k,\rho,-\eps)$. Furthermore, from the tensor product decomposition for $\pi^+_k \tensor \pi^P_{\rho,\eps}$ we obtain
\[
\pi^P_{\rho,\eps} \tensor \pi^-_k \cong \qdirint \pi^P_{\si,\eps'}\, d\si \oplus \bigoplus_{j} \pi^-_{j-\eps'}, \qquad \eps' = \eps-k.
\]
Here the direct sum is over all $j \in \Z$ such that $j-\eps' >0$.

\subsection{Racah coefficients} \label{ssec:qRacah}
We give explicit expressions for the Racah coefficients for various tensor products. The method is the same as in the Lie algebra case, therefore we omit the details. In this section we need the $q$-analogues of the Wilson functions and the Wilson polynomials, the Askey-Wilson functions and the Askey-Wilson polynomials. \\

The Askey-Wilson function transform is obtained by Koelink and Stokman \cite{KSt01A}. The Askey-Wilson function transform can be considered as a $q$-analogue of the the Wilson transform II. To stress the similarity between both integral transform, we use notations for the Askey-Wilson transform similar to the notations we used for the Wilson transform II.

Let the parameters $a,b,c,d,t$ satisfy the following conditions 
\[
t>0, \quad  0 < \Re(a), \Re(b)<1, \quad a = q/\overline d, \quad \overline b =c.
\]
We define dual parameters by
\[
\tilde a=|b|e^{i\arg a}, \quad \tilde b=\frac{ab}{\tilde a}, \quad \tilde c= \frac{ac}{\tilde a}, \quad \tilde d=\frac{ad}{\tilde a}, \quad \tilde t= \frac{\tilde d}{at}.
\]
Note that if $a$ and $d$ are real, we have $\tilde a =\sqrt{abcd/q}$.

We define the Askey-Wilson function by
\begin{equation} \label{def:AWf}
\begin{split}
\phi_\la&(x;a,b,c,d|q) =\\
&\frac{ (qa\la x^{\pm 1} /\tilde d,  q/\tilde d \la, ab, ac, qa/d ;q)_\infty }{ (\tilde a \tilde b \tilde c \la;q)_\infty}\ {}_8W_7(\tilde a \tilde b \tilde c \la/ q; ax, a/x, \tilde a \la, \tilde b \la, \tilde c \la;q, q/\tilde d \la).
\end{split}
\end{equation}
This is a slightly different function than the function used by Koelink and Stokman in \cite{KSt01A}. To identify $\phi_\la$ with the Askey-Wilson function in \cite{KSt01A}, let $\psi_\la$ be the function defined by \cite[(3.2)]{KSt01A}, then 
\[
\phi_\la(x;a,b,c,d|q) = G(x,\la) \psi_\la(x;a,b,c,d|q),
\]
where
\[
G(x,\la) =(qx^{\pm 1}/d, q\la^{\pm 1}/\tilde d, ab, ac, bc, qa/d ;q)_\infty.
\]
The normalization we use here corresponds to the normalization of the Wilson function defined by \eqref{def:Wilson function}. The reason for choosing this normalization is that the Askey-Wilson function $\phi_\la(x;a,b,c,d|q)$ is symmetric in the parameters $a,b,c,q/d$, while the function $\psi_\la(x;a,b,c,d|q)$ is not. This allows us to define the Askey-Wilson function transform that is unitary for a parameter domain that is different from the parameter domain used in \cite{KSt01A}, namely the parameters $a$, $b$, $c$, $d$ may be non-real.

By Bailey's transformation \cite[(III.36)]{GRa04} we have
\begin{equation} \label{eq:AWf 4ph3}
\begin{split}
\phi_\la(x;&a,b,c,d|q) = \frac{ (qx^{\pm 1}/d, q\la^{\pm 1}/\tilde d, ab,ac ;q)_\infty}{ (q/ad ;q)_\infty } \ph{4}{3}{ ax, a/x, \tilde a \la , \tilde a /\la }{ ab, ac , ad}{q;q} \\
&+  \frac{ (ax^{\pm 1}, \tilde a \la^{\pm 1}, qb/d, qc/d ;q)_\infty }{ (ad/q;q)_\infty} \ph{4}{3}{ qx/d, q/dx, q\la/\tilde d, q/\tilde d \la }{ qb/d, qc/d, q^2/ad}{q,q} .
\end{split}
\end{equation}
From this identity it follows that the Askey-Wilson function satisfies the duality property
\[
\phi_\la(x;a,b,c,d|q) =  \phi_x(\la;\tilde a,\tilde b, \tilde c, \tilde d|q),
\]
and that $\phi_\la(x;a,b,c,d|q)$ is symmetric in $a, q/d$ and in $b,c$. So, from $\overline{a}=q/d$ and $\overline{b}=c$ it follows that for $x,\la \in \R \cup \T$, the Askey-Wilson functions are real-valued. From Bailey's transformation \cite[(III.36)]{GRa04} with parameters $(a,b,c,d,e,f) \mapsto (\tilde a \tilde b \tilde c\la/q, ax, a/x, \tilde a \la, \tilde b \la, \tilde c \la)$ it even follows that $\phi_\la$ is symmetric in $a,b,c,q/d$, but we do not need this here.
 
Let $H$ be the weight function
\[
H(x)=H(x;a,b,c,d|q)  = \frac{ (x^{\pm 2};q)_\infty }{ (ax^{\pm 1}, bx^{\pm 1}, cx^{\pm 1}, qx^{\pm 1}/d;q)_\infty \te(tx^{\pm 1};q) }.
\] 
Here $\te(y;q) = (y, q/y ;q)_\infty$ denotes the Jacobi theta function. Let $\mathcal D$ be the infinite discrete set defined by
\[
\mathcal D = \{ tq^{-k}\ |\ k \in \Z,\ tq^{-k}>1 \}.
\]
We define a positive measure $h(\cdot)=h(\cdot;a,b,c,d;t|q)$ by
\[
\begin{split}
\int f(x)\, dh(x) = \frac{ C}{4\pi i} \int_\T f(x) H(x)  \frac{dx}{x} 
+ C \sum_{x \in \mathcal D} f(x) \Res{y=x}\left( \frac{H(y)}{y}\right)
\end{split}
\]
where $C$ is the normalizing constant given by
\[
C = (q;q)_\infty \sqrt{ \te(at;q) \te(bt;q) \te(ct;q) \te(qt/d;q) }.
\]
Explicitly we have, for $x= tq^{-k} \in \mathcal D$,
\[
\begin{split}
\Res{y=x}\left( \frac{H(y)}{y}\right) =& \frac{1}{(q,q,at,a/t,bt,b/t, ct, c/t qt/d, q/dt;q)_\infty } \\
\times& \frac{ (a/t, b/t, c/t, q/dt;q)_k }{ (q/at, q/bt, q/ct, d/t;q)_k } \left(1-q^{2k}/t^2 \right) (abct^2q/d)^{-2k} q^{k^2+k}.
\end{split}
\]
Let $\mathcal H(a,b,c,d;t)$ be the Hilbert space consisting of functions $f$ that satisfy $f(x) = f(x^{-1})$ ($h$-almost everywhere), with inner product $\langle \cdot, \cdot \rangle_{\mathcal H}$ defined by
\[
\langle f, g \rangle_{\mathcal H} = \int f(x) \overline{g(x)} dh(x).
\]
The Askey-Wilson function transform $\mathcal G : \mathcal H(a,b,c,d;t) \rightarrow  \mathcal H(\tilde a, \tilde b, \tilde c, \tilde d; \tilde t)$ is defined by
\[
(\mathcal G f)(\la) = \langle f, \phi_\la \rangle_{\mathcal H}.
\]
The Askey-Wilson function transform is unitary, with inverse $\mathcal G^{-1}= \tilde{\mathcal G}$, where $\tilde{\mathcal G}$ is the Askey-Wilson function transform with dual parameters. For real parameters $a,b,c,d$ this is Theorem 1 in \cite{KSt01A}, for non-real parameters this can be proved using the same method as in \cite{KSt01A}. In this section we show that the Askey-Wilson function with non-real parameters has an interpretation as a Racah coefficient. Then a similar reasoning as in subsection \ref{ssec:WilsonII} can be used to proof unitarity of $\mathcal G$ in this case. \\

To determine explicit expressions for the Racah coefficients we need the following $q$-analogue of Lemma \ref{lem1} and Lemma \ref{lem2}(iii). 
\begin{lem} \label{lem:qhulp}
The Askey-Wilson function satisfies:
\begin{align*}
(i)&&& \phi_\la(x;a\om,by,b/y,q\om/a|q) = \\
&&& \quad D\, \sum_{j=0}^\infty \frac{ (at\om, at/\om;q)_n }{ (q,a^2;q)_n } \,t^{-2n}\ph{3}{2}{q^{-n}, tx, t/x}{at\om,at/\om}{q,q} \ph{3}{2}{ q^{-n}, bty, bt/y}{abt\la, abt/\la}{q,q},\\
&&& D = \frac{ (a^2,byx,by/x,bx/y,b/yx,abt\la,abt/\la)_\infty }{ (bty, bt/y ;q)_\infty }, \\ \\
(ii) &&& \phi_{\tilde a q^n}(x;a,b,c,d|q) = \\
&&& \qquad (-d)^n q^{\hf n(n-1)} (qx^{\pm 1}/d,  abq^n, acq^n, bcq^n;q)_\infty\, p_n(\mu(x);a,b,c,d|q), \qquad n \in \Z_{\geq 0}.
\end{align*}
\end{lem}
\begin{proof}
The first statement is a special case of \cite[Cor.3.4]{KVdJ99}: Let $ct=c'$ in \cite[(3.21)]{KVdJ99}, then the $_3\varphi_2$-function on the right hand side of \cite[(3.21)]{KVdJ99} and the $_3\varphi_2$-function in the coefficient $G_j$ both become $_2\varphi_1$-functions that can be summed by the $q$-Gauss sum \cite[(II.8)]{GRa04}. This gives the following formula
\[
\begin{split}
\sum_{j=0}^\infty &\frac{ (ac,bc;q)_j }{ (q,ab;q)_j } \,c^{-2j}\ph{3}{2}{q^{-j}, cx, c/x}{ac,bc}{q,q} \ph{3}{2}{ q^{-j}, cty, ct/y}{a'ct, abct/a'}{q,q}  = \\
&\frac{ (bty^{\pm 1}, cty^{\pm 1}, abt x^{\pm 1}/a', a't/b;q)_\infty} { (tyx^{\pm 1}, tx^{\pm 1}/y, a'ct, abct/a', ab^2t/a' ;q)_\infty} \\
\times&\ {}_8W_7( ab^2t/a'q;bx, b/x, aby/a', ab/a'y , bt/a';q, a't/b) .
\end{split}
\]
Then the result follows from substituting 
\[
(a,b,c,a',t) \mapsto (a/\om, a\om, t, a/\la, b),
\]
and the definition \eqref{def:AWf} of the Askey-Wilson function.

The second statement follows directly from substituting $\la=\tilde a q^n$, $n \in \Z_{\geq 0}$, in \eqref{eq:AWf 4ph3}, and from the definition \eqref{def:AW pol} of the Askey-Wilson polynomial.
\end{proof}
There is also a $q$-analogue of Lemma \ref{lem2}(i),
\[
\phi_\la(x;a\rho,by,b/y, q \rho/a|q) = \phi_\rho(y;a\la, bx, b/x, q\la/a |q),
\]
which follows directly from the definition of the Askey-Wilson function, but we do not need this identity here.\\

First we consider the tensor product representation $\pi^+_{k_1} \tensor \pi^+_{k_2} \tensor \pi^-_{k_3}$. The Racah coefficients can be defined as the kernel in the (possibly discrete) integral transform that maps the eigenfuntions of the product of Clebsch-Gordan coefficients corresponding to $(\pi^+_{k_1} \tensor \pi^+_{k_2}) \tensor \pi^-_{k_3}$ to the Clebsch-Gordan coefficients corresponding to $\pi^+_{k_1} \tensor (\pi^+_{k_2} \tensor \pi^-_{k_3})$. 
We define the weight function $w(j;k_1,k_2,k_3,\tau)$ on $\Z_{\geq 0}$ as the mass of the measure 
\[
W(j;k_1,k_2) \times \mu_1(\tau;k_1+k_2+j,k_3),
\]
on a line $\tau=$ constant. Also we define the measure $\nu(\rho;k_1,k_2,k_3,\tau)$ to be the measure defined by
\[
\mu_2(\rho;k_2,k_3) \times \mu_3(\tau;k_1,\rho,k_2-k_3), 
\]
on the line $\tau=$ constant. The Racah coefficients for the tensor product $\pi^+_{k_1} \tensor \pi^+_{k_2} \tensor \pi^-_{k_3}$ can now be determined from \[
\begin{split}
v_n^{m+r}&(j;k_1,k_2) F_m^{r-j}(\tau;k_1+k_2+j,k_3) = \\
&\int  U(j,\rho;k_1,k_2,k_3,\tau)F_m^{r-n}(\rho;k_2,k_3) G_n^r(\tau;k_1,\rho,k_2-k_3) \, d\nu(\rho;k_1,k_2,k_3,\tau),
\end{split}
\]
or equivalently from 
\begin{equation} \label{eq:U++-q}
\begin{split}
&F_m^{r-n}(\rho;k_2,k_3) G_n^r(\tau;k_1,\rho,k_2-k_3) =\\
&\ \sum_{j=0}^\infty U(j,\rho;k_1,k_2,k_3,\tau) v_n^{m+r}(j;k_1,k_2) F_m^{r-j}(\tau;k_1+k_2+j,k_3)\,  w(j;k_1,k_2,k_3,\tau).
\end{split}
\end{equation}
Using the orthogonality relations for the Clebsch-Gordan coefficients, their explicit expressions as $_3\varphi_2$-series, and Lemma \ref{lem:qhulp}, we can express the Racah coefficient as a multiple of an Askey-Wilson polynomial.
\begin{thm} 
The Racah coefficient $U_j(\rho;k_1,k_2,k_3,\tau)$ is the Askey-Wilson polynomial
\[
N\ p_j(q^{2i\rho};q^{2k_2-2k_3+1}, q^{2k_2+2k_3-1}, q^{2k_1+2i\tau}, q^{2k_1-2i\tau}|q^2),
\]
where 
\[
\begin{split}
N =& \frac{2\pi\,q^{-2jk_1} q^{\hf j(j-1)}  \, (q^{2k_1+2k_2+2k_3-1+2j \pm 2i\tau};q^2)_\infty }{(q^{4k_1}; q^2)_j\, (q^2,q^2, q^{4k_3}, q^{4k_1+4k_2+4j}, q^{\pm 4i\tau} ;q^2)_\infty  } \\
&\times \frac{\te(q^{2k_1+2k_2-2k_3+1+2j \pm 2i\tau};q^2) } {\sqrt{ (q^{4k_1+4k_2+2j}, q^{1-2k_1-2k_2+2k_3-2j \pm 2i\tau} ;q^2)_j}}
\end{split}
\]
\end{thm}
\begin{rem}
There are not many group theoretic interpretations of the Wilson polynomials. In contrast with this, there are many interpretations in the context of quantum groups of the $q$-analogues of the Wilson polynomials, the Askey-Wilson polynomials. For example, Koornwinder \cite{Koo93} showed that they occur as spherical functions for $\mathrm{SU}_q(2)$, see also \cite{NM90} and \cite{Ko96}. In \cite{Ro00} the Askey-Wilson polynomials occur as matrix elements for representations of $\mathcal U_q(\su(1,1))$, in \cite{KVdJ98}, \cite{Gr02} and \cite{BR00} they occur as Clebsch-Gordan coefficients. 
\end{rem}

Next let us consider the tensor product representation $\pi^+_{k_1} \tensor \pi^P_{\rho,\eps} \tensor \pi^-_{k_3}$. Her we need the measure $\nu(\cdot;k_1,\rho,\eps,k_3,\tau)$, which is the restriction of the product measure
\[
\mu_2(\si;k_1,\rho,\eps) \times \mu_2(\tau;k_3,\si,-k_1-\eps),
\]
for $\tau=$ constant.
The Racah coefficients for the tensor product $\pi^+_{k_1} \tensor \pi^P_{\rho,\eps} \tensor \pi^-_{k_3}$ can be determined from 
\begin{equation} \label{eq:U+P-q}
\begin{split}
&(-1)^{n}\,G_m^{n-r}(\zeta;k_3,\rho,-\eps) G_n^{r}(\tau;k_1,\zeta,\eps-k_3) = \\
&\int U(\si,\zeta;k_1,\rho,\eps,k_3,\tau) G_n^{m+r}(\si;k_1,\rho,\eps) G_m^{-r} (\tau;k_3,\si,-k_1-\eps)\, d\nu(\si;k_1,\rho,\eps,k_3,\tau).
\end{split}
\end{equation}
The inverse for this formula is obtained by substituting $k_1 \leftrightarrow k_3$ and $\eps \mapsto -\eps$. This also leads to the following symmetry property for the Racah coefficients:
\begin{equation} \label{eq:symm U q}
U(\si,\zeta;k_1,\rho,\eps,k_3,\tau) = U(\zeta,\si;k_3,\rho,-\eps,k_1,\tau).
\end{equation}
Using the orthogonality of continuous dual $q^2$-Hahn polynomials and Lemma \ref{lem:qhulp}, we can express the Racah coefficient as a multiple of an Askey-Wilson function.
\begin{thm} \label{thm U+P-q}
For $\tau \in [0,-\pi/2\ln q)$ the Racah coefficient $U(\si,\zeta;k_1,\rho,\eps,k_3,\tau)$ is the Askey-Wilson function
\[
N\,\phi_{q^{2i\zeta}}(q^{2i\si};q^{2k_1+2i\rho},q^{2k_3+2i\tau}, q^{2k_3-2i\tau}, q^{2-2k_1+2i\rho}|q^2),
\]
where
\[
N= \frac{ 2\pi\ \te(q^{1+2k_1-2k_3+2\eps \pm 2i\tau};q^2) }{ (q^2,q^{4k_1}, q^{4k_3}, q^{\pm 4i\tau}, q^{1+2\eps+2k_1\pm 2i\si}, q^{1-2\eps+2k_3\pm 2i\zeta};q^2)_\infty}.
\]
\end{thm}
Observe that \eqref{eq:symm U q} is the same as the duality property for the Askey-Wilson funtions. In the same way as for the Lie algebra $\su(1,1)$ we can show that the Racah coefficients are the kernel in a unitary integral transform. This integral transform is equivalent to the Askey-Wilson transform. Let $\mathcal I$ denote the unitary integral transform
\[
\begin{split}
\mathcal I & \, : \mathcal L\big(\nu(\si;k_1,\rho,\eps,k_3,\tau)\big) \rightarrow \mathcal L\big(\nu(\zeta;k_3,\rho,-\eps,k_1,\tau),\\
(\mathcal I f)(\zeta) &= \int f(\si) U(\si,\zeta;k_1,\rho,\eps,k_3,\tau) d\nu(\si;k_1,\rho,\eps,k_3,\tau).
\end{split}
\]
This integral transform is related to the Askey-Wilson transform $\mathcal G$ in base $q^2$, with parameters
\[
(a,b,c,d,t) = (q^{2k_1+2i\rho},q^{2k_3+2i\tau}, q^{2k_3-2i\tau}, q^{2-2k_1+2i\rho},q^{1-2\eps-2k_1}),
\]
by
\[
(\mathcal If)(\zeta) = E\ \big( M^{-1}_{B(\zeta)} \circ \mathcal G \circ M_{A(\cdot)} \, f)(\zeta),
\]
where
\begin{gather*}
A(\si) = (q^{1+2k_1+2\eps \pm 2i\si};q^2)_\infty, \qquad B(\zeta) = (q^{1+2k_3-2\eps \pm 2i\zeta};q^2)_\infty, \\
E = \sqrt{ \frac{ (q^{1-2\eps \pm 2i\rho}, q^{1-2k_1+2k_3-2\eps \pm 2i\tau};q^2)_\infty } { (q^{1+2\eps \pm 2i\rho}, q^{1+2k_1-2k_3+2\eps \pm 2i\tau};q^2)_\infty } }.
\end{gather*}
\begin{rem}
(i) Similar to the Lie algebra case, for the values of $\tau$ corresponding to discrete series representations, the Askey-Wilson function in Theorem \ref{thm U+P-q} reduces to an Askey-Wilson polynomial, and the corresponding integral transform is the integral transform associated to these Askey-Wilson polynomials.

(ii) Koelink and Stokman \cite{KSt01B} found the Askey-Wilson function transform from the interpretation of the Askey-Wilson function as a spherical function on the quantum group $\mathrm{SU}_q(1,1)$. The duality property of the Askey-Wilson functions comes from the interpretation of the Askey-Wilson transform as a diffence Fourier transform related to a rank $1$ double affine Hecke algebra \cite{St03}. The interpretation as Racah coefficients given in Theorem \ref{thm U+P-q} gives a different explanation for the duality property.
\end{rem}

Let us also give the Racah coefficient for the tensor product representation 
$\pi^+_{k_1} \tensor \pi^+_{k_2} \tensor \pi^P_{\rho,\eps}$. Let the weight function $w(j;k_1,k_2,\rho,\eps,\tau)$ be the mass at the points $j \in \Z_{\geq 0}$ of the measure
\[
W(j;k_1,k_2) \times \mu_3(\tau;k_1+k_2+j,\rho,\eps)
\]
on the line $\tau=$ constant. Moreover, let $\nu(\si;k_1,k_2,\rho,\eps,\tau)$ be the restriction of the product measure
\[
\mu_3(\si;k_2,\rho,\eps) \times \mu_3(\tau;k_1,\si, k_2+\eps)
\]
to a line one which $\tau=$ constant. Then the Racah coefficient for the tensor product $\pi^+_{k_1} \tensor \pi^+_{k_2} \tensor \pi^P_{\rho,\eps}$ is determined from  
\[
\begin{split}
G_{r-m-n}^{r-n}&(\si;k_2,\rho,\eps) G_n^r(\tau;k_1,\si, k_2+\eps)= \\ \sum_{j=0}^\infty &U(j,\si;k_1,k_2,\rho,\eps,\tau) v_n^{r-m}(j;k_1,k_2) G_{r-m-j}^{r-j}(\tau;k_1+k_2+j,\rho,\eps) w(j;k_1,k_2,\rho,\eps,\tau),
\end{split}
\]
or from
\[
\begin{split}
v_n^{r-m}&(j;k_1,k_2) G_{r-m-j}^{r-j}(\tau;k_1+k_2+j,\rho,\eps)=\\
 \int &U(j,\si;k_1,k_2,\rho,\eps,\tau) G_{r-m-n}^{r-n}(\si;k_2,\rho,\eps) G_n^r(\tau;k_1,\si, k_2+\eps)\, d\nu(\si;k_1,k_2,\rho,\eps,\tau).
\end{split}
\]
Using orthogonality relations for the Clebsch-Gordan coefficients and Lemma \ref{lem:qhulp}, we can write the Racach coefficient as an Askey-Wilson polynomial. 
\begin{thm}
The Racah coefficient for the tensor product $\pi^+_{k_1} \tensor \pi^+_{k_2} \tensor \pi^P_{\rho,\eps}$ is given by
\[
U_j(\si;k_1,k_2,\rho,\eps) =K\  p_j(q^{2i\si};q^{2k_2+2i\rho}, q^{2k_2-2i\rho},q^{2k_1+2i\tau}, q^{2k_1-2i\tau}\,|q^2),
\]
where
\[
\begin{split}
N=& 2\pi\,\frac{ (q^{1-2k_1-2k_2-2\eps \pm 2i\tau}, q^{2k_1+2k_2+2j+2i\rho \pm 2i\tau}, q^{2k_1+2k_2+2j-2i\rho \pm 2i\tau} ;q^2)_\infty }{(q^{4k_1};q^2)_j\, (q^2,q^{4k_1+4k_2+4j},q^{1-2\eps \pm 2i\rho}, q^{\pm 4i\tau};q^2)_\infty} \\
&\times q^{-2jk_1} q^{\hf j(j-1)} \sqrt{ \frac{ (q^{1-2k_1-2k_2-2\eps-2j \pm 2i\tau};q^2)_j}{ (q^{4k_1+4k_2+2j};q^2)_j }} .
\end{split}
\] 
\end{thm}

Recall that in Lie algebra case the tensor product $\pi^+ \tensor \pi^- \tensor \pi^+$ corresponds to the Wilson function transform of type I. So it would be interesting to find the Racah coefficients for $\U_q$ for the tensor product $\pi^+ \tensor \pi^- \tensor \pi^+$, since the corresponding integral transform would give a $q$-analogue of the Wilson function transform of type I. The problem for finding these Racah coefficients is that we do not know the decomposition of $\pi^- \tensor \pi^+$, cf.~Remark \ref{rem:q -+}.

\subsection{Biedenharn-Elliott identities} 
In the same way as in subsection \ref{ssec:4fold} we find integral formulas or summation formulas for the Clebsch-Gordan coefficients and the Racah coefficients. These do not give essentially new formulas but only new interpretations of known formulas. Therefore we only write down two of them, which contain many other known formulas as special cases or as limit cases. 
Both these formulas are obtained from Biedenharn-Elliott identities. 
To simplify notations, we replace $q^2$ by $q$, $q^{k_i}$ by $k_i$ for $i=1,2,3,4$, and we replace $q^{\pm i\rho}$ by $\rho^{\pm 1}$ (and similarly for $\si,\zeta$ etc.). Moreover, in the formulas we will assume for simplicity $0 \leq k_i <1$, $j_i \in Z_{\geq 0}$ and $\rho,\si,\tau,\zeta, \xi \in \T$. 

We write out the Biedenharn-Elliott identity for $\pi^+_{k_1} \tensor \pi^+_{k_2} \tensor \pi^+_{k_3} \tensor \pi^-_{k_4}$. To find this identity we also need the Racah coefficient for the tensor product $\pi^+_{k_1} \tensor \pi^+_{k_2} \tensor \pi^+_{k_3}$. This Racah coefficient can be expressed as a multiple of the $q$-Racah polynomial, cf.~\cite{KVdJ98},
\[
R_{j_3}(j_2;q^{4k_2-2}, q^{4k_3-2}, q^{-2j_1-2j_2-2}, q^{4k_1+ 4k_2+2j_1+2j_2-2} ;q^2).
\]
The Biedenharn-Elliott identity leads the following convolution identity for Askey-Wilson polynomials.
\begin{thm} \label{thm:conv. AW}
The Askey-Wilson polynomials satisfy the following generalized convolution identity
\[
\begin{split}
\sum_{j_3=0}^{j_1+j_2}& \genfrac{[}{]}{0pt}{}{j_1+j_2}{j_3}_q \frac{k_2^{j_1-j_3}\ (k_3^2;q)_{j_1} (k_2^2;q)_{j_2} (k_1^2 k_2^2  k_3^2q^{j_1+j_2-1};q)_{j_3} }{ (k_3^2, k_2^2 k_3^2q^{j_3-1};q)_{j_3} (k_2^2 k_3^2 q^{2j_3};q)_{j_1+j_2-j_3} } \\
& \times R_{j_3}(j_2;k_2^2/q, k_3^2/q, q^{-j_1-j_2-1}, k_1^2 k_2^2 q^{j_1+j_2-1 } ;q) \\
& \times p_{j_1+j_2-j_3}(\si;k_2 k_3 q^{j_3+\hf}/k_4 , k_2 k_3 k_4 q^{j_3-\hf}, k_1\tau, k_1/\tau \,|q)\\
& \times p_{j_3}(\rho;k_3q^{\hf}/k_4 , k_3 k_4q^{-\hf}, k_2\si, k_2/\si\, |q) \\
= &p_{j_1}(\rho;k_3q^{\hf}/k_4 , k_3 k_4q^{-\hf}, k_1 k_2 \tau q^{j_2 }, k_1 k_2 q^{j_2}/\tau \,|q)\\
& \times  p_{j_2}(\si; k_2 \rho, k_2/\rho, k_1 \tau, k_1/\tau \, |q).
\end{split}
\]
\end{thm}
\begin{rem} Koelink and Van der Jeugt \cite[Thm.4.10]{KVdJ98} found this identity using the interpretation of the Askey-Wilson polynomials as Clebsch-Gordan coefficients for continuous basis vectors for positive discrete series representations of $\mathcal U_q\big(\su(1,1)\big)$. If we replace $k_4$ by $k_3k_4$ in Theorem \ref{thm:conv. AW} and then we put $k_3=0$, we obtain a generalized convolution identity for continuous dual $q$-Hahn polynomials which can also be obtained from writing out \eqref{eq:U++-q} and using the orthogonality relation for the dual $q$-Hahn polynomials. Theorem \ref{thm:conv. AW} contains many other generalized convolution identities for orthogonal polynomials as special cases, see \cite{KVdJ98}.
\end{rem}

From the Biedenharn-Elliott identity related to the tensor product $\pi^+_{k_1} \tensor \pi^+_{k_2} \tensor \pi^P_{\rho,\eps} \tensor \pi^-_{k_4}$ we obtain the following theorem.
\begin{thm} \label{thm:BE++P-q}
The Askey-Wilson functions and Askey-Wilson polynomials satisfy
\[
\begin{split}
\phi_{\si}(\rho;&\,k_1 k_2 \xi q^{j}, k_4\tau, k_4/\tau, \xi q^{1-j}/k_1 k_2\, |q)\ p_j(\zeta; k_2 \si, k_2/\si, k_1 \tau,k_1/\tau\, |q)\\ \frac{1}{4\pi i}&\int_\T  p_j(\om; k_2 \xi, k_2/\xi, k_1 \rho, k_1/\rho\, |q) \ \phi_{\tau}(\om; k_1\rho,  k_4 \zeta, k_4/\zeta, q\rho/k_1 \, |q) \\
 &\times \phi_{\si}(\om; k_2\xi, k_4 \zeta, k_4/\zeta, \xi q/k_2\, |q)  \\
& \times \frac{k_4^{-j}\, (\om^{\pm 2}, k_1 k_2 \xi q^{j}\rho^{\pm 1}, k_1 k_2 q^{j}\rho^{\pm 1}/\xi;q)_\infty }{ ( k_2 \xi \om^{\pm 1}, k_2 \om^{\pm 1}/ \xi,  k_4 \zeta  \om^{\pm 1}, k_4 \om^{\pm 1}/\zeta,  k_1 \rho \om^{\pm 1}, k_1  \om^{\pm 1}/\rho;q)_\infty} \frac{d\om}{\om}.
\end{split}
\]
\end{thm}
\begin{rem}
If $\si = k_4\xi q^{j_2}$ in Theorem \ref{thm:BE++P-q} and we use Lemma \ref{lem:qhulp}(ii), two Askey-Wilson functions in Theorem \ref{thm:BE++P-q} reduce to Askey-Wilson polynomials, and thus we obtain a formula with the following structure
\[
p_{j_1}(\zeta) p_{j_2}(\rho) = \int \phi_\tau(\om) p_{j_1}(\om) p_{j_2}(\om) d\om.
\]
This identity can also be obtained from the Biedenharn-Elliott identity for the tensor product $\pi^+_{k_1} \tensor \pi^+_{k_2} \tensor \pi^-_{k_3} \tensor \pi^-_{k_4}$. The Askey-Wilson polynomials in this identity can be further specialized to continuous dual $q$-Hahn polynomials, and this gives an identity that is equivalent to \eqref{eq:U+P-q}. Using the orthogonality for one of the continuous dual $q$-Hahn polynomials inside the integral then leads to the bilinear generating function for continuous dual $q$-Hahn polynomials \cite[Cor.3.4]{KVdJ99}. This is the generating function from which Lemma \ref{lem:qhulp}(ii) is derived as a special case.

Also, using the orthogonality relations for the Askey-Wilson polynomial inside the integral in Theorem \ref{thm:BE++P-q}, we find the bilinear generating function for Askey-Wilson polynomials \cite[Thm.3.3]{KVdJ99}. It is remarkable that in \cite{KVdJ99} the generating function is proved using only matrix elements and Clebsch-Gordan coefficients for positive discrete series representations of $\mathcal U_q$. 
\end{rem}


\end{document}